\documentclass[11pt]{article}
\usepackage[T1]{fontenc}
\usepackage{lmodern}
\usepackage[margin=1in]{geometry}
\usepackage{amsmath,amssymb,amsthm,mathtools,microtype,indentfirst}
\usepackage{tikz}
\usetikzlibrary{arrows.meta}
\definecolor{ofInk}{RGB}{35,39,43}
\definecolor{ofRed}{RGB}{177,48,58}
\definecolor{ofBlue}{RGB}{43,91,151}
\definecolor{ofRegionFill}{RGB}{244,247,250}
\definecolor{ofRegionBorder}{RGB}{151,165,179}
\tikzset{
  of diagram/.style={x=1cm,y=1cm,font=\small,
    line cap=round,line join=round,
    every node/.style={text=ofInk,inner sep=2pt,outer sep=0pt}},
  of vertex/.style={circle,draw=none,fill=ofInk,
    inner sep=0pt,minimum size=3.6pt,outer sep=0pt,label distance=2pt},
  of edge/.style={draw=ofInk,line width=.8pt},
  of red/.style={of edge,draw=ofRed,line width=.95pt},
  of blue/.style={of edge,draw=ofBlue,line width=.95pt,
    dash pattern=on 3pt off 2pt},
  of gap/.style={preaction={draw=white,line width=3.4pt,solid}},
  of bundle/.style={of edge,line width=1.6pt},
  of region/.style={draw=ofRegionBorder,fill=ofRegionFill,
    line width=.55pt,rounded corners=3pt},
  of set/.style={draw=ofRegionBorder,fill=white,line width=.55pt,
    rounded corners=2pt,minimum width=.85cm,minimum height=.48cm},
  of panel/.style={font=\small,inner sep=0pt},
  of arrow/.style={of edge,line width=.65pt,
    -{Stealth[length=4.5pt,width=4pt]}}
}
\usepackage{flafter}
\usepackage{multicol}
\usepackage[colorlinks=true,linkcolor=blue,citecolor=blue,urlcolor=blue]{hyperref}
\theoremstyle{definition}
\newtheorem{theorem}{Theorem}[section]
\newtheorem{lemma}[theorem]{Lemma}
\newtheorem{proposition}[theorem]{Proposition}
\newtheorem{corollary}[theorem]{Corollary}
\numberwithin{equation}{section}
\hypersetup{pdftitle={A Sharp Form of the Overfull Conjecture},
pdfauthor={Guantao Chen, Frederic Folz and Songling Shan}}
\title{A Sharp Form of the Overfull Conjecture}
	\author{Guantao Chen \thanks{Georgia State
		University, Department of Mathematics and Statistics,  Atlanta, GA 30303, USA. Email: \href{mailto:gchen@gsu.edu}{\texttt gchen@gsu.edu}.  Partially supported by NSF grant DMS-2154331.}
	\quad
	Frederic Folz\thanks{Altifina Discovery Lab,
		66802 \"Uberherrn, Germany. Email:
		\href{mailto:research@altifina.com}{\texttt{research@altifina.com}}.}
	\quad
	Songling Shan\thanks{Auburn University, Department of Mathematics
		and Statistics, Auburn, AL 36849, USA. Email:
		\href{mailto:szs0398@auburn.edu}{\texttt{szs0398@auburn.edu}}.
		Partially supported by NSF grant DMS-2451895.}}
\date{November 25, 2026}
\begin{document}
\maketitle
\emph{\textbf{Abstract}.}
Determining whether a graph of maximum degree $\Delta$ has an edge
$\Delta$-coloring is a fundamental problem in graph theory.
For a simple graph $G$ of order $n$ and maximum degree $\Delta$,
a subgraph $H$ is \emph{$\Delta$-overfull} if
$e(H)>\Delta\lfloor |V(H)|/2\rfloor$.
The Overfull Conjecture, proposed by Chetwynd and Hilton in 1986,
asserts that such a subgraph is the only obstruction to an edge
$\Delta$-coloring when $\Delta>n/3$.
In this paper, we prove that, for all sufficiently large $n$,
a graph $G$ of order $n$ with maximum degree
$\Delta\ge2\lfloor(n-1)/8\rfloor$ has an edge $\Delta$-coloring
if and only if it contains no $\Delta$-overfull subgraph.
Thus we confirm the Overfull Conjecture for sufficiently large
graphs and extend its conclusion to a sharp degree bound near
$n/4$.
This characterization gives a polynomial-time algorithm for determining
the chromatic index in the stated range and strengthens the
$1$-Factorization Conjecture for sufficiently large graphs.
We also obtain the Just-overfull and Vertex-splitting Conjectures
for sufficiently large graphs, and Vizing's Independence Number,
Average Degree, and $2$-Factor Conjectures in the corresponding
degree ranges. For regular graphs, we further prove the
characterization at $\Delta\ge2\lfloor n/10\rfloor$.
Both bounds are best possible for each parity of $\Delta$.
For every integer $\Delta\ge3$, we construct a connected
class~2 graph with maximum degree $\Delta$ and no $\Delta$-overfull
subgraph. Its order is $4\Delta+9$ when $\Delta$ is even and
$4\Delta+5$ when $\Delta$ is odd. We also construct
connected $\Delta$-regular counterexamples on $5\Delta+10$
vertices when $\Delta$ is even and $5\Delta+5$ vertices when
$\Delta$ is odd.


\begingroup
\small
\setlength{\columnsep}{1.5em}
\hyphenpenalty=10000
\begin{multicols}{2}[\section*{\normalsize Contents}]
\makeatletter
\@starttoc{toc}
\makeatother
\end{multicols}
\endgroup
\clearpage

\section{Introduction}\label{sec:introduction}
In this paper, a \emph{graph} is finite and simple, and a
\emph{multigraph} may have multiple edges but has no loops.
We write $e(G)$ for the number of edges of $G$.
Degrees and edge counts in a multigraph include multiplicity.
For a nonnegative integer $k$, write $[k]=\{1,\ldots,k\}$, with
$[0]=\emptyset$. An
\emph{edge $k$-coloring} of a multigraph $G$ is a mapping
$\varphi:E(G)\to[k]$ such that adjacent edges receive different
colors. The edges receiving any one color form a \emph{color class}.
The \emph{chromatic index} of $G$, denoted by $\chi'(G)$, is the
smallest integer $k$ for which $G$ has an edge $k$-coloring. A
\emph{$1$-factorization} is a partition of the edge set into
perfect matchings.

Let $G$ be a graph of order $n$ and maximum degree $\Delta$.
By Vizing's Theorem~\cite{Vizing64}, $\chi'(G)$ is either $\Delta$
or $\Delta+1$. The graph $G$ is of \emph{class~1} in the first
case and of \emph{class~2} in the second. Distinguishing these two
classes is a central problem in edge-coloring. Holyer~\cite{Holyer1981}
showed that determining whether a graph is of class~1 is NP-complete,
even when the graph is cubic.

An elementary obstruction to an edge
$\Delta$-coloring is an overfull subgraph. A subgraph $H$ of $G$
is \emph{$\Delta$-overfull} if
$e(H)>\Delta\lfloor |V(H)|/2\rfloor$.
Since every color class restricted to $H$ is a matching of size at
most $\lfloor |V(H)|/2\rfloor$, the presence of such a subgraph
forces $G$ to be of class~2.

In 1986, Chetwynd and Hilton~\cite{ChetwyndHilton} proposed the
\emph{Overfull Conjecture}, which states that, if $\Delta>n/3$,
then $G$ is of class~2 if and only if it contains a
$\Delta$-overfull subgraph. Thus the conjecture would make an
elementary density condition a complete criterion for determining
which of Vizing's two values occurs. Since the existence of a
$\Delta$-overfull subgraph can be decided in polynomial time by an
odd-cut algorithm~\cite{PadbergRao1982}
(see Corollary~\ref{pre:overfull-detection}), this characterization
would give a polynomial-time algorithm for determining the
chromatic index of graphs with $\Delta>n/3$.
The conjecture also connects this edge-coloring problem with
several longstanding conjectures concerning critical graphs, including
the Just-overfull and Vertex-splitting Conjectures and Vizing's
Independence Number, Average Degree, and $2$-Factor Conjectures;
see~\cite{SSTF,HiltonZhao1997}.

For $X\subseteq V(G)$, let $e_G(X)=e(G[X])$. For a positive integer $k$,
an odd set $X$ is \emph{$k$-overfull} if
$e_G(X)>k(|X|-1)/2$. A graph contains a $\Delta$-overfull subgraph
if and only if it contains a $\Delta$-overfull set. Indeed, a
$\Delta$-overfull subgraph must have odd order, and replacing it
by the subgraph induced on the same vertex set preserves the
overfull inequality. We will use the set formulation throughout
the proofs.

The $1$-Factorization Conjecture is a special case of the Overfull
Conjecture. It asserts that a $\Delta$-regular graph of even order
$n$ has a $1$-factorization whenever
$\Delta\ge2\lceil n/4\rceil-1$. Csaba, K\"uhn, Lo, Osthus, and
Treglown~\cite{CKLOT2016} proved this conjecture for sufficiently
large graphs. Much of the other progress on the Overfull
Conjecture has also required a minimum-degree assumption; see,
for example,~\cite{Plantholt2022,PlantholtShan2023,ShanOdd2024}.

When only the maximum degree is restricted, Chetwynd and
Hilton~\cite{ChetwyndHilton1989} proved the conjecture for
$\Delta\ge n-3$. For sufficiently large graphs, Shan~\cite{Shan2026}
proved it for $\Delta\ge13n/14$, and Chen, McDonald, and
Shan~\cite{OverfullII} proved it for
$\Delta\ge(1+\varepsilon)n/2$ for every fixed $0<\varepsilon<1$.
In this paper, we prove the following sharp form of this
characterization under a maximum-degree condition alone.

\begin{theorem}\label{main:quarter}
There exists a positive integer $n_0$ for which the following holds.
If $G$ is a graph of order $n\ge n_0$ with maximum degree
$\Delta\ge2\lfloor(n-1)/8\rfloor$, then
$\chi'(G)=\Delta$ if and only if $G$ contains no
$\Delta$-overfull subgraph.
\end{theorem}

The degree condition in Theorem~\ref{main:quarter} is equivalent to
$\lfloor(n-1)/8\rfloor\le\lfloor\Delta/2\rfloor$, or
$n\le8\lfloor\Delta/2\rfloor+8$. Thus it becomes
$n\le4\Delta+8$ when $\Delta$ is even and
$n\le4\Delta+4$ when $\Delta$ is odd.
Theorem~\ref{main:quarter} therefore proves the Overfull Conjecture for
sufficiently large graphs and extends its characterization to
a maximum-degree bound near $n/4$.

The bound is best possible for both parities of $\Delta$.
In Proposition~\ref{sharp:thm:family}, we construct a connected
class~2 graph with maximum degree $\Delta$ and no $\Delta$-overfull
subgraph, of order $4\Delta+9$ when $\Delta$ is even and
$4\Delta+5$ when $\Delta$ is odd, for every integer $\Delta\ge3$.
Each counterexample has the first order excluded by the
corresponding parity-dependent bound.
Together with Corollary~\ref{pre:overfull-detection},
Theorem~\ref{main:quarter} gives a polynomial-time algorithm for
determining the chromatic index throughout the stated range.

Throughout the following corollaries, $n_0$ is as in
Theorem~\ref{main:quarter}, enlarged if necessary.
Every color class in an edge $\Delta$-coloring of a $\Delta$-regular
graph is a perfect matching. Thus Theorem~\ref{main:quarter} gives the
following corollary. It is stronger than the $1$-Factorization
Conjecture for sufficiently large graphs, since a $\Delta$-regular
graph of even order $n$ with $\Delta\ge2\lceil n/4\rceil-1$
contains no $\Delta$-overfull set. Indeed, the smaller side of an
odd cut has odd order $s\le2\lceil n/4\rceil-1\le\Delta$, so the
cut has at least $s(\Delta-s+1)\ge\Delta$ edges, since
$s(\Delta-s+1)-\Delta=(s-1)(\Delta-s)\ge0$. Regularity then
gives $e_G(X)\le\Delta(|X|-1)/2$ for every odd set $X$.

\begin{corollary}\label{intro:cor:one-factorization}
Every $\Delta$-regular graph $G$ of order $n\ge n_0$ with
$\Delta\ge2\lfloor(n-1)/8\rfloor$ and no
$\Delta$-overfull set has a $1$-factorization.
\end{corollary}

For regular graphs, the same characterization holds at the sharp
degree bound near $n/5$. We prove the following strengthening in
Section~\ref{sec:proof-regular} by deleting the edges of one dense
component and adapting the proof of Theorem~\ref{main:quarter} to a
critical subgraph of the remaining graph.

\begin{theorem}\label{intro:thm:regular}
There exists a positive integer $n_0$ for which the following holds.
A simple $\Delta$-regular graph $G$ of order $n\ge n_0$ with
$\Delta\ge2\lfloor n/10\rfloor$ has a $1$-factorization if and
only if it contains no $\Delta$-overfull set.
\end{theorem}

We take $n_0$ large enough for both theorems.
A regular graph of positive degree with no $\Delta$-overfull set has even order,
since otherwise its whole vertex set is overfull. To express the
degree condition in Theorem~\ref{intro:thm:regular} as an order bound,
observe that $\Delta\ge2\lfloor n/10\rfloor$ is equivalent to
$\lfloor n/10\rfloor\le\lfloor\Delta/2\rfloor$, and hence to
$n<10(\lfloor\Delta/2\rfloor+1)$. Since $n$ is even, the last
inequality is equivalent to $n\le10\lfloor\Delta/2\rfloor+8$.
When $\Delta$ is even, $\lfloor\Delta/2\rfloor=\Delta/2$,
giving $n\le5\Delta+8$. When $\Delta$ is odd,
$\lfloor\Delta/2\rfloor=(\Delta-1)/2$, giving $n\le5\Delta+3$.
Proposition~\ref{sharp:prop:regular} gives connected regular
counterexamples at the first even order beyond each bound,
namely $5\Delta+10$ and $5\Delta+5$, respectively.

The characterization in Theorem~\ref{main:quarter} also has
consequences for critical graphs.
A graph $G$ of order $n$ and maximum degree $\Delta$ is
\emph{$\Delta$-critical} if its chromatic index
is $\Delta+1$ and every proper subgraph is $\Delta$-edge-colorable.
It is \emph{just overfull} if
$e(G)=\Delta\lfloor n/2\rfloor+1$.
The Just-overfull Conjecture states that every $\Delta$-critical
graph of order $n$ with $\Delta\ge n/2$ is just overfull;
see~\cite[Conjecture~4.23]{SSTF} and~\cite{CaoChenShan2022Core}.
As a consequence of Theorem~\ref{main:quarter}, we obtain this
conclusion for $n\ge n_0$ under the weaker degree condition
$\Delta\ge2\lfloor(n-1)/8\rfloor$.
\begin{corollary}\label{intro:cor:just-overfull}
Every $\Delta$-critical graph $G$ of order $n\ge n_0$ with
$\Delta\ge2\lfloor(n-1)/8\rfloor$
has odd order and is just overfull.
In particular, $e(G)=\Delta(n-1)/2+1$.
\end{corollary}
\begin{proof}
By Theorem~\ref{main:quarter}, $G$ has a $\Delta$-overfull set $X$.
Criticality forces $X=V(G)$, since otherwise $G[X]$ would be a
proper subgraph that is not $\Delta$-edge-colorable.
Thus $n$ is odd and $e(G)>\Delta(n-1)/2$.
For any edge $e$, a $\Delta$-edge-coloring of $G-e$ gives
$e(G)-1\le\Delta(n-1)/2$, proving the equality.
\end{proof}

Another consequence concerns vertex-splitting. A
\emph{vertex-splitting} in a graph $G$ replaces a vertex $v$ by two adjacent
vertices $x,y$ and partitions the neighbors of $v$ into two nonempty sets,
assigned respectively to $x$ and $y$.
The Vertex-splitting Conjecture of Hilton and
Zhao~\cite{HiltonZhao1997} concerns a connected $\Delta$-regular
class~1 graph $D$ with $\Delta>|V(D)|/3$ and predicts that
the resulting graph is $\Delta$-critical.
A revised formulation excludes a specific nine-vertex
counterexample~\cite{QiTianGutin2026}. The following corollary establishes
the conclusion in the degree range of Theorem~\ref{main:quarter},
where $n$ denotes the order after splitting.

\begin{corollary}\label{intro:cor:vertex-splitting}
If $D$ is a connected $\Delta$-regular class~1 graph of order $n-1$,
where $n\ge n_0$ and $\Delta\ge2\lfloor(n-1)/8\rfloor$,
then every graph obtained from $D$ by a vertex-splitting is
$\Delta$-critical.
\end{corollary}
\begin{proof}
Let $H$ be obtained by splitting $v\in V(D)$ into $x,y$.
A $\Delta$-edge-coloring of $D$ decomposes it into perfect
matchings, so $n$ is odd.
The graph $H$ has order $n$, maximum degree $\Delta$, and
$e(H)=\Delta(n-1)/2+1$. Hence $H$ is of class~2.

For every $e\in E(H)$, we have $\Delta(H-e)=\Delta$: among
the $n-2$ vertices inherited from $D-v$, at most two are incident
with $e$, and every other one still has degree $\Delta$. By Theorem~\ref{main:quarter},
it therefore suffices to show that $H-e$ has no
$\Delta$-overfull set. Its whole vertex set has exactly
$\Delta(n-1)/2$ edges, so we need only exclude proper
$\Delta$-overfull sets of $H$.

Let $S$ be a proper odd subset of $V(H)$. If $S$ contains at most
one of $x,y$, then
$H[S]$ is isomorphic to a subgraph of $D$ and is therefore
$\Delta$-edge-colorable.
If $S$ contains both $x,y$, contracting $xy$ gives a nonempty
proper set $T\subset V(D)$ of even size $|S|-1$, with
$e_H(S)=e_D(T)+1$. Since $D$ is connected and $\Delta$-regular,
$e_D(T)<\Delta|T|/2$. The latter bound is an integer, so
$e_H(S)\le\Delta|T|/2=\Delta(|S|-1)/2$.
Thus no $H-e$ has a $\Delta$-overfull set. Every $H-e$ is therefore $\Delta$-edge-colorable. Since $H$ is
connected and has no isolated vertices, each proper subgraph is
contained in some $H-e$. Hence $H$ is $\Delta$-critical.
\end{proof}

Finally, we consider three conjectures of Vizing. Let
$\alpha(G)$ and $\overline d(G)$ denote the independence number
and average degree of $G$, respectively.
For a $\Delta$-critical graph $G$ of order $n$, Vizing's
Independence Number and Average Degree Conjectures assert that
$\alpha(G)\le n/2$ and $\overline d(G)\ge\Delta-1+3/n$,
respectively; see~\cite{Vizing1968,SSTF,CaoChen2021}.
His $2$-Factor Conjecture asserts that $G$ has a spanning
$2$-regular subgraph; see~\cite{KannoShan2019}.
Corollary~\ref{intro:cor:just-overfull} yields all three conclusions
in the degree range of Theorem~\ref{main:quarter}.

\begin{corollary}\label{intro:cor:independence}
Every $\Delta$-critical graph $G$ of order $n\ge n_0$ with
$\Delta\ge2\lfloor(n-1)/8\rfloor$ satisfies $\alpha(G)\le(n-1)/2$.
\end{corollary}
\begin{proof}
By Corollary~\ref{intro:cor:just-overfull}, $n$ is odd and
$e(G)>\Delta(n-1)/2$.
For any independent set $I$, every edge has an end outside $I$,
so $e(G)\le\sum_{v\notin I}d_G(v)\le\Delta(n-|I|)$.
Consequently, $|I|<(n+1)/2$, giving the claimed bound.
\end{proof}

\begin{corollary}\label{intro:cor:average-degree}
Every $\Delta$-critical graph $G$ of order $n\ge n_0$ with
$\Delta\ge2\lfloor(n-1)/8\rfloor$ satisfies
\[
 \overline d(G)=\Delta-\frac{\Delta-2}{n}
 \ge\Delta-1+\frac{3}{n}.
\]
\end{corollary}
\begin{proof}
The equality follows from Corollary~\ref{intro:cor:just-overfull}.
The inequality follows from $\Delta\le n-1$.
\end{proof}

\begin{corollary}\label{intro:cor:two-factor}
Let $G$ be a $\Delta$-critical graph of order $n\ge n_0$ with
$\Delta\ge2\lfloor(n-1)/8\rfloor$. Then every edge of $G$
belongs to a $2$-factor of $G$.
In particular, $G$ has a $2$-factor.
\end{corollary}
\begin{proof}
Fix an edge $xy$. Criticality gives a $\Delta$-edge-coloring
of $G-xy$. By Corollary~\ref{intro:cor:just-overfull}, $n$ is odd
and $e(G-xy)=\Delta(n-1)/2$. Since every color class has at
most $(n-1)/2$ edges, each class has exactly this size and is
missing at exactly one vertex.
Choose a color $c_1$ missing at $x$ and a color $c_2$
missing at $y$; these exist because the edge $xy$ was deleted.
Since each color is missing at exactly one vertex, $c_1\ne c_2$,
$c_1$ is present at $y$, and $c_2$ is present at $x$.
The union of their color classes has degree one at $x,y$ and
degree two at every other vertex. Adding $xy$ gives a $2$-factor
containing the prescribed edge.
\end{proof}

The remainder of this paper is organized as follows. In
Section~\ref{sec:cases}, we partition the bipartite double cover and
prove Theorems~\ref{main:quarter} and~\ref{intro:thm:regular} using
boundary colorings and prescribed palettes.
Section~\ref{sec:preliminaries} contains the preliminary results.
The boundary-coloring theorem and its regular-graph variant are
proved in Section~\ref{sec:boundary}, and the prescribed-palette
theorem is proved in Section~\ref{sec:palettes}.
Section~\ref{sharp:sec} contains the sharpness constructions for
both main theorems.

The main results of this paper have been formally verified in Lean~4 (v4.12.0). 
The complete formalization is  available on GitHub at \url{https://github.com/Altifina-Dev/overfull-lean}.

\section{Proof of the main theorems}\label{sec:cases}
We prove Theorems~\ref{main:quarter} and~\ref{intro:thm:regular}
using a common partition and coloring argument.
The average-degree bound of Cao and Chen makes the total degree
deficiency $\sum_v(\Delta-d_H(v))$ of a critical counterexample $H$
small. Sparse cuts in its
bipartite double cover then yield dense ordinary and bipartite components
with many walks of a common odd length. We color the edges outside
these components, balance the missing colors, and extend the coloring
inside each component. For regular graphs, we apply this argument
after removing the edges of one component, using the boundary variant in
Lemma~\ref{boundary:lem:regular-ambient}.

The exact parity-dependent order bound in
Theorem~\ref{main:quarter} is used when correcting the
boundary coloring. In the only remaining obstruction, four disjoint
sets of original vertices have even orders, and at least one original
vertex lies outside their union. The bound guarantees that one of
these sets has at most $\Delta$ vertices. This supplies the additional
edges needed to remove the obstruction; see Step~4 of the proof of
Proposition~\ref{boundary:prop:recolor}. In the sharpness construction,
all four even sets have more than $\Delta$ vertices.

We use $G,H,J,D,R$ for graphs and multigraphs, with subscripts
or primes for related graphs. The letters $C,P,M$ denote cycles
or circuits, paths, and matchings, respectively; $K_r$ denotes a
complete graph and $Q_k$ a binary cube. Sets of vertices or edges
use uppercase letters such as $A,B,U,V,W,X,Y$, while $x,y,z,u,v,w$ and their
indexed or primed forms are reserved for vertices. Colors are denoted
by $c$ and its indexed or primed forms. Calligraphic letters denote
families or sets of colors; in particular, $\mathcal P(v)$ is the
palette at $v$.

Numerical parameters use other lowercase letters or Greek letters. We use $n,m,N$ for orders or cardinalities,
$L$ for a walk length, and $i,j,k,\ell,q,r,s,t$ for indices or
integer counts. We retain the standard degree notation
$d_G(v)$, $\delta(G)$ and $\Delta(G)$, and specify additional
notation when it is introduced. In probabilistic arguments, $Z$
denotes a random variable and $\xi$ a centered increment.

All unions that describe partitions are disjoint. For a vertex $v$ and a set
$X\subseteq V(G)$, let $d_G(v,X)$ be the number of edges
from $v$ to $X$, counted with
multiplicity in a multigraph. For disjoint sets $X,Y$, let
$G[X,Y]$ have vertex set $X\cup Y$ and all edges of $G$ between
$X$ and $Y$, and set $e_G(X,Y)=e(G[X,Y])$.
We write $\partial_G(X)$ for the set of edges with exactly one
endpoint in $X$.
For sets $A,B$, write $A\oplus B=(A\setminus B)\cup(B\setminus A)$
for their symmetric difference.
A matching covers a vertex if the vertex is an end of one of
its edges. For a matching $M$, write $V(M)$ for its set of
covered vertices. For an edge set $B$, let $d_B(v)$ count its
incidences at $v$, and set $d_B(Y)=\sum_{v\in Y}d_B(v)$ for
a vertex set $Y$.
The constants in the auxiliary results
are independent of the order and maximum degree of the graph.
A decomposition of a graph into subgraphs means that their edge sets
partition the edge set of the graph; the subgraphs may share vertices.

Let $u,v$ be vertices of a simple graph $H$, not necessarily distinct.
A \emph{$(u,v)$-walk of length $L$ in $H$} is a sequence
$u=x_0,x_1,\ldots,x_L=v$ such that
$x_{i-1}x_i\in E(H)$ for every $i\in[L]$.
Vertices and edges may repeat. Two such walks are distinct if their
vertex sequences differ. All walk counts below count distinct walks;
the walks need not be edge-disjoint or vertex-disjoint.
A walk is \emph{closed} if its endpoints coincide, a \emph{trail}
if no edge repeats, and a \emph{path} if no vertex repeats.
For walks and trails in a multigraph, the sequence also specifies
the edge copy used at each step.
Subgraphs are \emph{edge-disjoint} if their edge sets are disjoint,
and \emph{vertex-disjoint} if their vertex sets are disjoint.

\subsection{Sparse cuts in the bipartite double cover}
For a vertex partition of a graph $H$ with exceptional set $U$,
we associate a subgraph with each remaining part $V_i$.
An \emph{ordinary component} is the induced graph $H[V_i]$.
A \emph{bipartite component} is the graph $H[A_i,B_i]$ for a
specified partition $V_i=A_i\cup B_i$; it contains only the edges
between its two sides. The term component refers to these specified
subgraphs. The graph $H$ may also have edges between different
components or within a side of a bipartite component.
Degrees and walks in a component are taken within its specified
subgraph. The walk conditions below make each component connected
and each ordinary component nonbipartite.

The components supplied below retain all but a small number of
incident edges at each vertex and contain many walks of a common
odd length between the relevant pairs of vertices. These are the
properties needed for the boundary-coloring and prescribed-palette
arguments.

The robust-partition theorem of K\"uhn, Lo, Osthus, and
Staden~\cite{KLOS2015} gives a more general structural description.
A sparse cut in the bipartite double cover can separate two ordinary
components or distinguish the sides of a bipartite component. The following
argument supplies the degree and walk bounds in both cases.
We state it in the range $m\le5\Delta+8$ so that it also applies
to Theorem~\ref{intro:thm:regular}. In the narrower range of
Theorem~\ref{main:quarter}, the minimum-degree bound will force
at most four ordinary components and individual bipartite sides altogether.
We state the conclusion as four properties of a vertex partition.
All its components use the same odd walk length $L$ and the same
coefficient $b$ in the lower bound $bm^{L-1}$. The tolerance
$\varepsilon(b,L)$ specifies how small the exceptional set and
the degree loss must be for these walk parameters.
\begin{lemma}\label{elem:partition}
For each $0<b\le1$ and each odd integer $L\ge3$, let
$\varepsilon(b,L)\in(0,1)$ be prescribed. There exist a constant
$\zeta>0$ and a positive integer $m_0$, depending only on this
function $\varepsilon$, for which the following holds.
Let $H$ be a simple graph of order $m\ge m_0$ and maximum degree
$\Delta$, with $m\le5\Delta+8$ and
$\sum_{v\in V(H)}(\Delta-d_H(v))\le\zeta m^2$.

Then $H$ has a vertex partition
\[
 V(H)=U\cup V_1\cup\cdots\cup V_{q+\ell},
\]
where $U$ may be empty and $q,\ell$ are nonnegative integers,
together with numbers
$0<b\le1$ and $h\ge0$ and an odd integer $L\ge3$, satisfying
the following properties.
\begin{enumerate}
\item $1\le q+2\ell\le5$ and $|U|\le h\le\varepsilon(b,L)m$.

\item Each part $V_{q+j}$, for $j\in[\ell]$, has a partition
$V_{q+j}=A_j\cup B_j$. The ordinary components
$H_i=H[V_i]$ for $i\in[q]$ and the bipartite components
$H_{q+j}=H[A_j,B_j]$ for $j\in[\ell]$ satisfy
$\delta(H_i)\ge\Delta-h>m/6$ for every $i\in[q+\ell]$.

\item For each $i\in[q]$ and every $u,v\in V_i$, including
$u=v$, the graph $H_i$ contains at least $bm^{L-1}$ distinct
$(u,v)$-walks of length $L$.

\item For each $j\in[\ell]$ and every $(u,v)\in A_j\times B_j$,
the graph $H_{q+j}$ contains at least $bm^{L-1}$ distinct
$(u,v)$-walks of length $L$.
\end{enumerate}
Moreover, there is a finite list of pairs $(b,L)$, depending only
on $\varepsilon$ and independent of $H$, from which the common
pair $(b,L)$ in the conclusion can always be chosen.
\end{lemma}

This last assertion allows the order thresholds in later lemmas
to be chosen simultaneously for all possible walk parameters,
before the graph $H$ and its partition are known.

\begin{proof}
Let $D$ be the bipartite double cover of $H$, with one vertex $x^+$
and one vertex $x^-$ for each $x$, and edges $x^+y^-,y^+x^-$ for
each edge $xy$. Write $\sigma$ for the involution interchanging the two copies.
We partition the cover so that $\sigma$ permutes its parts. We call these
parts cells until we project them back to $H$.

We use the increasing cut thresholds
$t_1m^2,\ldots,t_{20}m^2$. Fix $0<t_{20}<10^{-6}$.
For $i=20,\ldots,1$, given $t_i$, set $g_i=t_i/16$, let $L_i$ be
the least odd integer at least $4/t_i+5$, and set $b_i=g_i^{L_i-1}$.
Choose $0<t_{i-1}<t_i$ so small that
\[
 18\sqrt{t_{i-1}}<
 \min\{t_i/8,\ \varepsilon(b_i,L_i),\ 1/100\}.
\]
The finite list in the last assertion is
$\{(b_i,L_i):i\in[20]\}$. Set $\zeta=t_0/2$, and choose $m_0$
large enough for the inequalities below. These choices depend
only on $\varepsilon$.

Start with one cell. After $j$ refinements, split a cell if it has a cut
with both sides of size at least $m/40$ and fewer than $t_{j+1}m^2$
edges between them inside the cell. Split its image under $\sigma$ by
the image cut as well. If the cell is fixed by $\sigma$, take the common
refinement by both cuts: intersect each side of one cut with each
side of the other, and omit empty intersections.
At most $2t_{j+1}m^2$ additional edges have their ends in
distinct cells after this refinement.
Consequently, writing $W(z)$ for the cell containing a cover vertex $z$,
\[
 \lambda_j=\sum_z\bigl(\Delta-d_D(z,W(z))\bigr)
 \le t_0m^2+4m^2\sum_{i=1}^j t_i \le 81t_jm^2
 \qquad(0\le j\le20),
\]
where degrees in cells are taken in the induced graph of the cover.
The first term counts the deficiency at both copies of every vertex. Each newly
separated edge increases the sum by two. The two scales have
different roles: $t_j$ bounds the degree loss accumulated so far,
whereas $t_{j+1}$ is the cut bound used to test the current cells.
The choice of $t_j$ ensures that deleting vertices with large
degree loss leaves enough edges across every cut that passes the
next test.

Call a cell large if it has more than $m/10$ vertices. Since
$\Delta\ge(m-8)/5\ge m/6$ for large $m$, every vertex in a cell of
size at most $m/10$ contributes at least $m/15$ to $\lambda_j$.
All such cells together therefore contain at most $1215t_jm<m/40$
vertices. This bound also holds after the next refinement, with
$t_{j+1}$ in place of $t_j$. Each side of the chosen cut is a
union of new cells and has at least $m/40$ vertices. It cannot
consist entirely of small new cells, so it contains a large one.
The same holds for the image cut. An unaffected large cell remains large,
and each cell being split is replaced by at least two large cells.
Thus every refinement strictly increases the number of large cells.
There are fewer than twenty in a cover of order $2m$, so the process
stops at a stage $j<20$.

Choose the common walk parameters $(b,L)=(b_{j+1},L_{j+1})$
from the fixed list. Set $t=t_{j+1}$, $g=t/16$,
$s=9\sqrt{t_j}$ and $h=2sm$. Delete every
vertex with $\Delta-d_D(z,W(z))>sm$. At most $\lambda_j/(sm)\le sm$ vertices are deleted.
Each surviving vertex loses at most another $sm$ neighbors,
so its degree in the surviving cell is at least $\Delta-2sm=\Delta-h$. The deficiency of a vertex equals that of its image under $\sigma$,
so deletion respects $\sigma$. Our parameter choice gives
$s<1/200<1/15$. Every vertex in a small cell has deficiency
at least $m/15>sm$ and is therefore deleted. Fewer than $m/10$
vertices are deleted altogether. Thus a vertex survives in every large cell. Each surviving
cell has minimum degree at least $\Delta-h>m/6$ for large $m$.
Every cut in a surviving cell with both sides of size at least $m/40$
still has at least $tm^2/2$ edges. Indeed, assign all vertices deleted from that cell to either side
of the cut. Both sides still have at least $m/40$ vertices, so
the stopping rule applies to the resulting cut, and
at most $sm^2$ of its edges meet a deleted vertex, since a cover
vertex has degree at most $m$.

Fix one surviving cell $W$, set $J=D[W]$, and let $A,B$ be its sides. For $S\subseteq A$, let
$T$ consist of the vertices of $B$ with at least $gm$ neighbors in $S$.
Suppose $|S|,|A\setminus S|\ge m/40$. The cut at $S\cup T$ meets the
preceding size condition. The degree bounds in $J$ give
\[
 \frac{tm^2}{2}
 \le \sum_{y\in T}d_J(y)-\sum_{x\in S}d_J(x)+2e_J(S,B\setminus T)
 \le \Delta(|T|-|S|)+hm+2gm^2.
\]
Since $h/m<t/8$ and $\Delta\le m$, it follows that
$|T|\ge |S|+tm/4$. The same calculation holds from the other side.

From any vertex $x$, start with its neighbors and
iterate the preceding neighborhood operation, alternating sides.
After $k$ steps, for every vertex $y$ in the resulting set there
are at least $(gm)^{k-1}$ distinct $(x,y)$-walks of length $k$ in $J$:
each new endpoint has at least $gm$
neighbors in the preceding set, and extending the walks through
these different neighbors gives distinct vertex sequences. As long as the current set misses at least
$m/40$ vertices of its side, it grows by at least $tm/4$ at the next
step. The starting set has more than $m/6$ vertices and neither side
has more than $m$, so there can be fewer than $4/t$ such growth
steps. Once fewer than $m/40$ vertices are missing, every vertex
on the opposite side has at least $\Delta-h-m/40>gm$ neighbors
in the set. The next reached set is therefore the whole opposite
side, and every subsequent reached set is also a whole side.
This happens within $\lceil4/t\rceil+2$ steps. Since the chosen
odd length satisfies $L\ge4/t+5$, every vertex $y$ on the side
opposite $x$ is reached at step $L$, with at least
$(gm)^{L-1}=bm^{L-1}$ distinct $(x,y)$-walks in $J$.

The deleted pairs project to $U$. A cell fixed by $\sigma$ is the double
cover of the graph induced on its projection. In two distinct cells
interchanged by $\sigma$, the two sides of either cell project to
disjoint sets, and that cell gives precisely the crossing graph on
these sets. Group each cell with its image under $\sigma$, using
one cell when it is fixed. Different groups project to disjoint
vertex sets of $H$. Denote the ordinary projections by
$V_1,\ldots,V_q$, and the bipartite projections by
$V_{q+j}=A_j\cup B_j$ for $j\in[\ell]$. Together with $U$,
these sets form the required vertex partition. Their component
graphs have the minimum degrees in property~(2).
A walk in $H$ determines a unique walk in the cover once its
starting copy is fixed. The same walk counts therefore hold within the corresponding
subgraph of $H$, including from $x$ to itself in an ordinary
component. This proves properties~(3) and~(4).
The vertex set of every ordinary component and each individual
bipartite side has size greater than $m/6$, by the minimum-degree
bound. These $q+2\ell$ sets are disjoint, so $1\le q+2\ell\le5$.
Finally, $|U|\le sm/2\le h$ and $h\le\varepsilon(b,L)m$,
which completes property~(1).
\end{proof}

\subsection{Proof of Theorem~\ref{main:quarter}}\label{bip:part:section}
We first state the two edge-coloring theorems used in the proof.
The first theorem colors the edges outside a collection of dense
subgraphs. A \emph{private leaf} at $v$ is a new vertex adjacent
only to $v$. We attach one such leaf for each unit of degree
deficiency $\Delta-d_G(v)$; its edge records a color missing at $v$.
We also call these vertices \emph{deficiency leaves} and their
incident edges \emph{deficiency edges}.
\begin{theorem}\label{boundary:thm:main}
There exist a constant $0<\alpha<1$ and a positive integer
$\Delta_1$ for which the following holds. Let $\Delta$ be an
integer and let $G$ be a graph of order $n$ and maximum degree
at most $\Delta$, where
$\Delta\ge\max\{\Delta_1,2\lfloor(n-1)/8\rfloor\}$.
Suppose that $G$ has a vertex partition
$V(G)=U\cup V_1\cup\cdots\cup V_k$, with $1\le k\le4$ and
$U$ possibly empty, satisfying the following conditions.
\begin{enumerate}
\item $G$ has no $\Delta$-overfull set.
\item $|U|\le\alpha\Delta$ and
$\delta(G[V_i])\ge(1-\alpha)\Delta$ for every $i\in[k]$.
\item For some $i\in[k]$ and some $e\in E(G[V_i])$, the graph
$G-e$ has an edge $\Delta$-coloring.
\end{enumerate}
Let $D$ be obtained from $G$ by attaching $\Delta-d_G(v)$
distinct private leaves to each $v\in V(G)$, and set
\[
 J=D-\bigcup_{i=1}^k E(G[V_i]).
\]
Then $J$ has an edge $\Delta$-coloring in which every color is
missing at an even number of vertices of each $V_i$.
\end{theorem}

The second theorem extends this coloring inside each part.
At each vertex, its internal edges must receive exactly the colors
left unused by the exterior coloring. This required set is its
\emph{palette}. We use two settings. In the ordinary case, $H$ is
a graph with no specified bipartition, and we set $m=|V(H)|$.
In the bipartite case, $H$ has specified sides $A,B$, and we set
$m=\max\{|A|,|B|\}$. The proof uses odd walks, whose ends
must lie on opposite sides in the bipartite case.
An ordered pair $(u,v)$ is \emph{admissible} if
$(u,v)\in V(H)\times V(H)$ in the ordinary case, and if
$(u,v)\in A\times B$ in the bipartite case. In the ordinary case,
$u=v$ is permitted.

\begin{theorem}\label{palette:thm:main}
For every $b>0$ and every odd integer $L\ge3$, there exist a
constant $\eta=\eta(b,L)>0$ and a positive integer
$m_0=m_0(b,L)$ for which the following holds.
Let $\Delta$ be a positive integer and let $H$ be a simple graph.
In the ordinary case, set $m=|V(H)|$. In the bipartite case,
let $A,B$ be the specified sides of $H$ and set
$m=\max\{|A|,|B|\}$. For each $v\in V(H)$, let
$\mathcal P(v)\subseteq[\Delta]$. Suppose that $m\ge m_0$ and that
the following conditions hold.
\begin{enumerate}
\item In the ordinary case, $H$ contains at least $bm^{L-1}$
distinct $(u,v)$-walks of length $L$ for every $u,v\in V(H)$,
including $u=v$. In the bipartite case, the same bound holds
for every $(u,v)\in A\times B$.
\item For every $v\in V(H)$,
$|\mathcal P(v)|=d_H(v)$ and $\Delta-|\mathcal P(v)|\le\eta m$.
\item In the ordinary case, each color $c\in[\Delta]$ belongs
to an even number of palettes and is absent from at most
$\eta m$ palettes. In the bipartite case, each color belongs
to equally many palettes on $A$ and $B$ and is absent from
at most $\eta m$ palettes on each side.
\end{enumerate}
Then $H$ has an edge $\Delta$-coloring whose set of colors at
each vertex $v$ is exactly $\mathcal P(v)$.
\end{theorem}

We prove Theorems~\ref{boundary:thm:main} and~\ref{palette:thm:main}
in Sections~\ref{sec:boundary} and~\ref{sec:palettes}, respectively,
using the boundary and palette arguments in~\cite{Folz2026}.
The boundary theorem retains the original endpoints of the exterior
edges, while the palette theorem permits unequal bipartite sides.
Both statements therefore allow the degree deficiencies that remain
after the partition.

We also use the following average-degree bound for critical graphs.

\begin{lemma}[Cao and Chen {\cite[Theorem~3]{CaoChen2020}}]
\label{bip:part:cao-chen}
For every $0<\varepsilon<1$, there exists a positive integer
$\Delta_0=\Delta_0(\varepsilon)$ such that every $\Delta$-critical
graph $G$ of order $n$ with $\Delta\ge\Delta_0$ satisfies
\[
 \overline d(G):=\frac1n\sum_{v\in V(G)}d_G(v)
 \ge(1-\varepsilon)\Delta.
\]
\end{lemma}

\begin{proof}[Proof of Theorem~\ref{main:quarter}]
Let $\alpha,\Delta_1$ be as in Theorem~\ref{boundary:thm:main},
and choose $n_0$ sufficiently large for all requirements below.
A $\Delta$-overfull set prevents an edge $\Delta$-coloring.
For the converse, let $G$ have order $n\ge n_0$, maximum degree
$\Delta\ge2\lfloor(n-1)/8\rfloor$, and no $\Delta$-overfull set.
Suppose that $G$ is not $\Delta$-edge-colorable. Choose an
edge-minimal subgraph $H$ that is not $\Delta$-edge-colorable,
and discard its isolated vertices. Its maximum degree is
$\Delta$: otherwise Vizing's Theorem gives an edge
$\Delta$-coloring. Thus $H$ is $\Delta$-critical and contains
no $\Delta$-overfull set. Set $m=|V(H)|$. Then
\[
 \Delta+1\le m\le n\le4\Delta+8,
 \qquad \Delta\ge2\lfloor(m-1)/8\rfloor.
\]
Since $\Delta\ge(n-8)/4$ and $m\ge\Delta+1$, both $m$ and
$\Delta$ tend to infinity with $n$. Increasing $n_0$ therefore meets
all the order and degree requirements below, including
$\Delta\ge\Delta_1$.

\medskip
\noindent\textbf{Step 1. Finding the components.}
Let $\eta(b,L)$ be the tolerance in
Theorem~\ref{palette:thm:main}. Use Lemma~\ref{elem:partition} with
$\varepsilon(b,L)=\min\{1,\alpha,\eta(b,L)\}/200$ to obtain
the finite list of parameter pairs and the deficiency tolerance $\zeta$.
The average-degree bound of Lemma~\ref{bip:part:cao-chen}, with error $\zeta$, gives
$\sum_v(\Delta-d_H(v))\le\zeta m\Delta\le\zeta m^2$ for the critical
graph $H$, so the partition lemma applies to $H$.
Denote its resulting components by $H_1,\ldots,H_{q+\ell}$ and
their vertex sets by $V_i$, listing the ordinary components first. Write
$A_j,B_j$ for the sides of part $V_{q+j}$. Set $s=q+\ell$ and let $m_i$
be the order of $H_i$ in the ordinary case and its larger side size in
the bipartite case. Since $m\le4\Delta+8$ and $h\le m/200$,
we have $\Delta-h>m/5$ for large $m$. Each ordinary component and
each individual bipartite side has more than $m/5$ vertices,
so $q+2\ell\le4$. In particular, $m/5<m_i\le m$.
Our choice of tolerance, together with
$\Delta>m/5$, gives $6h\le\alpha\Delta$ and
$12h\le3\eta(b,L)m/50$. Since the list of pairs is finite,
we can take $m$ large enough that every $m_i$ meets the palette
theorem's order requirement and $12h+10\le\eta(b,L)m_i$. Each required walk count is at least
$bm^{L-1}\ge bm_i^{L-1}$, so the same pair $(b,L)$ works in every part.

\medskip

\noindent\textbf{Step 2. Coloring the boundary.}
We color the edges outside the $H_i$ first. The missing colors
will then give the palettes for coloring the $H_i$.
Attach $\Delta-d_H(v)$ private leaves to every original vertex
$v$, let $D$ be the resulting graph, and let $W$ be its set of
leaves. Every original vertex has degree $\Delta$ in $D$.
Set $J=D-\bigcup_{i=1}^s E(H[V_i])$.
Choose an edge $e$ of $H_1$. Criticality gives an edge
$\Delta$-coloring of $H-e$. Since $s\le4$, $|U|\le h\le\alpha\Delta$,
and $\delta(H[V_i])\ge\Delta-h\ge(1-\alpha)\Delta$ for every
$i\in[s]$, Theorem~\ref{boundary:thm:main} gives an edge
$\Delta$-coloring of $J$ in which every color is missing at
an even number of vertices of every $V_i$.
Every color class covers $U$, since $d_J(u)=\Delta$ for $u\in U$.
It covers a number of vertices congruent to $m$ modulo two in
$V(H)$, and hence also in $W$, because a matching covers an even
number of vertices altogether.

For each $j\in[\ell]$, insert the edges within $A_j$ and within
$B_j$ greedily. At each endpoint, at most $h$ edges
lie outside the specified crossing graph, so fewer than
$2h<\Delta$ colors are excluded. Every inserted edge contributes
two incidences to its part and preserves all required parities.
Write $J'$ for the graph now colored. It consists of all edges
of $D$ outside the graphs $H_i$.
Every original vertex outside $U$ has at most $h$ incident
edges in $J'$. Every edge has an original endpoint, so summing
degrees over original vertices gives
\[
 e(J')\le |U|\Delta+mh\le6h\Delta.
\]

\medskip
\noindent\textbf{Step 3. Balancing and extending the coloring.}
We balance the coloring of $J'$ so that each color is
missing equally often on the two sides of every bipartite component
and every color class is small.
Let $M_c$ be its color class $c$. For $j\in[\ell]$, set
$s_j(c)=|A_j\cap V(M_c)|-|B_j\cap V(M_c)|$ and
$d_j=|A_j|-|B_j|$. Since each color is missing at an even number
of vertices of $V_{q+j}$, we have $s_j(c)\equiv d_j\pmod2$. Equality of the two degree sums
in the crossing graph gives
\[
 \sum_{c=1}^{\Delta}s_j(c)=\Delta d_j.
\]
Apply Lemma~\ref{pre:balance}, taking the vertex sets of the
$\ell$ bipartite components as signed parts, and those of the
$q$ ordinary components and, if nonempty, $W$ as unsigned parts.
It gives $s_j(c)=d_j$ for all $j,c$, while preserving all
part parities and coverage of $U$. The color classes differ
in size by at most $2\ell+q+1\le5$, and hence
\[
 |M_c|\le e(J')/\Delta+5\le6h+5.
\]

For $v\in V_i$, let $\mathcal P(v)$ be the set of colors
missing at $v$ in this coloring of $J'$.
Every original vertex has degree $\Delta$ in $D$, so
$|\mathcal P(v)|=d_{H_i}(v)$ and $\Delta-|\mathcal P(v)|\le h$.
Each color is absent from
the palettes of at most $2|M_c|\le12h+10$ vertices, even when
all parts are considered together. If $i\le q$, each color belongs to the palettes of an even number
of vertices of $H_i$. If $i=q+j$, each color belongs to equally
many palettes on $A_j$ and $B_j$, by $s_j(c)=d_j$.
Together with the walk counts and the bounds from Step~1,
these facts verify every hypothesis of
Theorem~\ref{palette:thm:main} for each $H_i$.
The resulting colorings use only colors missing in $J'$ at
each vertex, so they combine with the coloring of $J'$ to give
an edge $\Delta$-coloring of $D$. Deleting the deficiency edges
gives an edge $\Delta$-coloring of $H$, a contradiction.
\end{proof}

\subsection{Proof of Theorem~\ref{intro:thm:regular}}\label{sec:proof-regular}
For the regular case we need the following variant of
Theorem~\ref{boundary:thm:main}, proved in
Section~\ref{boundary:sec:regular}. It uses the order and regularity
of an ambient graph in place of the order bound on the subgraph
whose boundary is being colored.

\begin{lemma}\label{boundary:lem:regular-ambient}
There exist a constant $0<\alpha<1$ and a positive integer
$\Delta_1$ for which the following holds. Let $R$ be a simple
$\Delta$-regular graph of even order $N$, where
$\Delta\ge\max\{\Delta_1,2\lfloor N/10\rfloor\}$.
Let $G\subseteq R$ have a vertex partition
$V(G)=U\cup V_1\cup\cdots\cup V_k$, with $1\le k\le4$ and
$U$ possibly empty. Suppose that the following conditions hold.
\begin{enumerate}
\item $G$ has no $\Delta$-overfull set and
$|V(R)\setminus V(G)|\ge13$.
\item $|U|\le\alpha\Delta$ and
$\delta(G[V_i])\ge(1-\alpha)\Delta$ for every $i\in[k]$.
\item For some $i\in[k]$ and some $e\in E(G[V_i])$, the graph
$G-e$ has an edge $\Delta$-coloring.
\end{enumerate}
Let $D$ be obtained from $G$ by attaching $\Delta-d_G(v)$
distinct private leaves to each $v\in V(G)$, and set
\[
 J=D-\bigcup_{i=1}^k E(G[V_i]).
\]
Then $J$ has an edge $\Delta$-coloring in which every color is
missing at an even number of vertices of each $V_i$.
\end{lemma}

We delete the edges of one dense component, color the remaining
graph using at most four components, and then extend the coloring
over the deleted edges. We apply the boundary lemma to a critical
subgraph of the remaining graph, with the original regular graph
as its ambient graph. Any vertices of the chosen component that
belong to this critical subgraph lie in its exceptional set. The
other vertices of the component are absent from the critical
subgraph, as required by
Lemma~\ref{boundary:lem:regular-ambient}.

\begin{proof}[Proof of Theorem~\ref{intro:thm:regular}]
The necessity of the odd-set condition is immediate. For sufficiency,
let $G$ be a simple $\Delta$-regular graph of order $n$ with no
$\Delta$-overfull set and
$\Delta\ge2\lfloor n/10\rfloor$. The order $n$ is even, since
otherwise $V(G)$ itself would be overfull. In particular,
$n\le5\Delta+8\le6\Delta$ for large $n$.

Let $\alpha$ be the constant in
Lemma~\ref{boundary:lem:regular-ambient}, and let $\eta(b,L)$
be the palette tolerance in Theorem~\ref{palette:thm:main}.
Apply Lemma~\ref{elem:partition} with
\[
 \varepsilon(b,L)=\frac{\min\{1,\alpha,\eta(b,L)\}}{1000}.
\]
Fix the resulting finite list and deficiency tolerance $\zeta$.
Take $n_0$ large enough for this partition lemma, the average-degree
bound of Lemma~\ref{bip:part:cao-chen} with error $\zeta$, the
boundary lemma, and the palette theorem for every pair on the list.
All further inequalities hold by increasing $n_0$ if necessary.
Since $\Delta\ge(n-8)/5$, this also makes $\Delta$ sufficiently large.

\smallskip\noindent\textbf{Step 1. Deleting the edges of one component.}
The graph $G$ has zero degree deficiency, so the partition lemma
applies. Choose one of its components $H_0$,
set $S=V(H_0)$, and denote its parameters by $b,L,h$.
Write $m_S$ for its order in the ordinary case and for its larger
side size in the bipartite case. We have
\[
 |S|\ge\Delta-h+1,\qquad
 h\le\varepsilon(b,L)n<\Delta/100.
\]
Delete just the edges of $H_0$, and set $R=G-E(H_0)$.
Every vertex outside $S$ still has degree $\Delta$, whereas
$d_R(v)\le h$ for $v\in S$. In the bipartite case the edges
within either side of $H_0$ remain in $R$.

\smallskip\noindent\textbf{Step 2. Coloring the remaining graph.}
Suppose that $R$ has no edge $\Delta$-coloring. Choose an
edge-minimal subgraph $H$ with this property and discard its
isolated vertices. Vizing's theorem gives $\Delta(H)=\Delta$,
so $H$ is $\Delta$-critical. It has no $\Delta$-overfull set,
because it is a subgraph of $G$. Set $m=|V(H)|$; then
$\Delta+1\le m\le n\le5\Delta+8$.
Lemma~\ref{bip:part:cao-chen} gives
\[
 \sum_{v\in V(H)}(\Delta-d_H(v))
 \le\zeta m\Delta\le\zeta m^2.
\]
Apply the same partition lemma to $H$. This is a new partition;
the components of the first partition need not remain dense in $H$.
Let $U'$ be its exceptional
set, let $H_1,\ldots,H_s$ be its components on parts $V_i$,
and let $b',L',h'$ be the resulting parameters. Write
$s=q+\ell$, listing the $q$ ordinary components before the
$\ell$ bipartite components. Then
\[
 |U'|\le h',\qquad
 \delta(H_i)\ge\Delta-h',\qquad
 h'\le\varepsilon(b',L')m<\Delta/100.
\]

Every vertex in $S\cap V(H)$ has degree at most $h$ in $H$,
which is less than $\Delta-h'$. Thus all these vertices lie
in $U'$, and all the dense parts lie outside $S$.
Each ordinary component and each individual bipartite side has at
least $\Delta-h'$ vertices. If there were five such sets, then
\[
 5(\Delta-h')\le n-|S|\le4\Delta+h+7,
\]
contrary to $h,h'<\Delta/100$ for large $\Delta$.
Consequently $q+2\ell\le4$, and in particular $s\le4$.
Moreover,
\[
 |V(G)\setminus V(H)|
 \ge |S|-|U'|\ge\Delta-h-h'+1\ge13.
\]
This is why we retain $G$ as the ambient graph when coloring $H$: the vertices of $S$ outside $U'$ supply the required
vertices omitted from $H$.

Attach deficiency leaves to $H$, obtaining $D$, and set
$J=D-\bigcup_{i=1}^s E(H[V_i])$.
Choose an edge $e$ of $H_1$. Criticality gives an edge
$\Delta$-coloring of $H-e$. Our choice of tolerance gives
$h'\le\alpha\Delta$, so
Lemma~\ref{boundary:lem:regular-ambient}, with ambient graph $G$,
colors $J$ with every color missing at an even number of
vertices of every $V_i$.

We complete the coloring as in Steps~2--3 of the proof of
Theorem~\ref{main:quarter}. We give the estimates needed in
the present order range. Insert the edges within the sides
of each bipartite component greedily; fewer than $2h'<\Delta$ colors
are excluded at each edge. The resulting graph $J'$ consists
of all edges of $D$ outside the specified graphs $H_i$.
Every vertex of $U'$ has degree $\Delta$ in $J'$, and every
other original vertex has degree at most $h'$. Hence
$e(J')\le |U'|\Delta+mh'\le7h'\Delta$.

Apply Lemma~\ref{pre:balance}, taking $U'$ as the set that every
color class must cover. Take the vertex
sets of the bipartite components as signed parts, and those of the
ordinary components and, if nonempty, the set of deficiency leaves
as unsigned parts. Every color class has the same incidence parity
on the leaves: its incidence count on $V(H)$ is congruent to $m$
modulo two, and its total incidence count is even.
Write $M_c$ for the current color class $c$. For a bipartite
component with sides $A_j,B_j$, set
$s_j(c)=|A_j\cap V(M_c)|-|B_j\cap V(M_c)|$ and
$d_j=|A_j|-|B_j|$. The degree sums in its crossing graph give
$\sum_c s_j(c)=\Delta d_j$, and the even missing-color count
on $A_j\cup B_j$ gives $s_j(c)\equiv d_j\pmod2$.
Thus the hypotheses of the balancing lemma hold. It makes the
missing-color counts equal on the two sides of every bipartite
component and makes the color classes differ in size by at most
$q+2\ell+1\le5$.
Every color therefore covers at most $14h'+10$ vertices.
The palettes satisfy $|\mathcal P(v)|=d_{H_i}(v)$ and
$\Delta-|\mathcal P(v)|\le h'$. Each color occurs in an even number of
palettes in an ordinary component, and equally often on the two
sides of a bipartite component. The tolerance chosen above ensures that both
$h'$ and $14h'+10$ are at most $\eta(b',L')m_i$, where $m_i$
is the order or larger side size of $H_i$; indeed
$m_i\ge\Delta-h'$ and $m\le6\Delta$.
The walk counts also hold with $m_i$ in place of $m$,
because $m_i\le m$.
Theorem~\ref{palette:thm:main} colors every $H_i$, giving an
edge $\Delta$-coloring of $H$, a contradiction.
Thus $R$ has an edge $\Delta$-coloring.

\smallskip\noindent\textbf{Step 3. Coloring the deleted edges.}
Fix such a coloring of $R$. Every color class covers all vertices
outside $S$, since their degrees in $R$ equal $\Delta$.
If $M_c$ is color class $c$ and
$t_c=|V(M_c)\cap S|$, then
\[
 t_c\equiv n-|S|\equiv|S|\pmod2,
 \qquad t_c=2|M_c|-(n-|S|).
\]
We now make every $t_c$ small, so that each color is absent from
only a small number of the required palettes.
In the ordinary case, each color is therefore missing at an even
number of vertices of $S$. Apply Lemma~\ref{pre:balance}, taking
$V(G)\setminus S$ as the set that every color class must cover
and $S$ as one unsigned part.

In the bipartite case, write $S=A\cup B$ and set
$d=|A|-|B|$. For each color $c$, set
$s(c)=|A\cap V(M_c)|-|B\cap V(M_c)|$.
Since $s(c)$ and $t_c$ differ by $2|B\cap V(M_c)|$,
we have $s(c)\equiv t_c\equiv d\pmod2$. Since every deleted
edge has one end on each side,
\[
 \sum_c s(c)
 =\sum_{v\in A}d_R(v)-\sum_{v\in B}d_R(v)
 =\Delta d.
\]
Apply Lemma~\ref{pre:balance}, again requiring every color class
to cover $V(G)\setminus S$, and taking $A,B$ as the sides of
one signed part. It gives $s(c)=d$,
so every color is missing equally often on $A$ and $B$.

In either case, the color classes now differ in size by at most
two, and thus the numbers $t_c$ differ by at most four. Their
mean is at most $|S|h/\Delta$, giving
\[
 t_c\le\frac{1}{\Delta}\sum_{v\in S}d_R(v)+4
      \le\frac{|S|h}{\Delta}+4.
\]
Let $\mathcal P(v)$ be the colors missing at $v\in S$.
Then $|\mathcal P(v)|=d_{H_0}(v)$ and
$\Delta-|\mathcal P(v)|\le h$. A color is absent from at most $t_c$
of these palettes. Since $|S|\le2m_S$,
$h\le6\varepsilon(b,L)\Delta$, and $m_S\ge\Delta-h$,
our choices give
\[
 h\le\eta(b,L)m_S,\qquad
 t_c\le12\varepsilon(b,L)m_S+4\le\eta(b,L)m_S.
\]
All hypotheses of Theorem~\ref{palette:thm:main} now hold for
$H_0$. Its palette coloring extends the coloring of $R$ to
an edge $\Delta$-coloring of $G$. Since $G$ is regular, its
color classes form a $1$-factorization.
\end{proof}

\section{Preliminaries}\label{sec:preliminaries}
We collect the auxiliary results used in the proof and in the
algorithmic consequence.

We use Hall's Theorem in the following form. A bipartite
multigraph with sides $A,B$ has a matching covering $A$ if and
only if $|N(S)|\ge|S|$ for every $S\subseteq A$, where $N(S)$
is the set of neighbors of $S$. In particular, a simple bipartite
graph with $|A|=|B|$ and minimum degree at least $|A|/2$ has
a perfect matching.

\subsection{Detecting overfull sets}
We detect overfull sets by expressing the odd-set inequality as
a cut-capacity inequality. We first introduce the cut terminology.
For an undirected graph $H$ with edge capacities $\operatorname{cp}(e)\ge0$, set
$\operatorname{cp}(B)=\sum_{e\in B}\operatorname{cp}(e)$ for $B\subseteq E(H)$.
If $T\subseteq V(H)$ has even cardinality, a cut $\partial_H(X)$
is \emph{$T$-odd} if $|X\cap T|$ is odd. The two sides of a
cut have the same parity of intersection with $T$. If
$T=\emptyset$, there is no $T$-odd cut.

\begin{theorem}[Padberg and Rao {\cite{PadbergRao1982}}]
\label{pre:odd-cut}
Let $H$ be an undirected graph with nonnegative rational edge
capacities, and let $T\subseteq V(H)$ have even cardinality.
A minimum-capacity $T$-odd cut can be found, or its nonexistence
determined, in time polynomial in the input length.
\end{theorem}

\begin{corollary}\label{pre:overfull-detection}
Let $k$ be a positive integer and let $G$ be a graph with
$\Delta(G)\le k$. A $k$-overfull set can be found, or its
nonexistence determined, in time polynomial in
$|V(G)|+e(G)+\log(k+1)$.
\end{corollary}
\begin{proof}
Add a vertex $z$ to $G$. Give each original edge capacity one
and each new edge $zv$, for $v\in V(G)$, capacity $k-d_G(v)$.
These capacities are nonnegative. Let $H$ be the resulting
graph, and set $T=V(G)$ when $|V(G)|$ is even and
$T=V(G)\cup\{z\}$ when $|V(G)|$ is odd. Then $|T|$ is even.

For any $T$-odd cut, choose its side $X$ not containing $z$.
Since complementation preserves $T$-parity and $T\cap X=X$,
the set $X\subseteq V(G)$ has odd order. Conversely, every
odd set $X\subseteq V(G)$ gives a $T$-odd cut of $H$. Moreover,
\[
 \operatorname{cp}(\partial_H(X))
 =|\partial_G(X)|+\sum_{v\in X}(k-d_G(v))
 =k|X|-2e_G(X).
\]
Thus $X$ is $k$-overfull if and only if
$\operatorname{cp}(\partial_H(X))<k$. Apply Theorem~\ref{pre:odd-cut}.
If the minimum capacity is less than $k$, the side avoiding $z$
of a minimum $T$-odd cut is the required set. If the minimum is at least
$k$, or no $T$-odd cut exists, then no $k$-overfull set exists.
The construction adds one vertex and $|V(G)|$ edges, and all its capacities are
integers between zero and $k$, so its input length is polynomial in
$|V(G)|+e(G)+\log(k+1)$.
\end{proof}

\subsection{Balancing color classes}
To extend the boundary coloring, we need both the prescribed
parity of the number of vertices covered in each part and, for each
color, equal numbers of vertices missing that color on the two
sides of each bipartite component. The signed parts below
record these two sides.
A part with a specified partition into sides $A_i,B_i$ is
called signed; all other parts are unsigned.
The \emph{incidence parity}, or \emph{part parity}, of a color
class $M$ on a part $V$ is $|V(M)\cap V|\pmod2$.

\begin{lemma}\label{pre:balance}
Let $J$ have an edge $\Delta$-coloring with color classes
$M_c$, $c\in[\Delta]$, and let $U\subseteq V(J)$ be covered
by every color class. Partition $V(J)\setminus U$ into
$k_{\mathrm{s}}$ signed parts $A_i\cup B_i$ and
$k_{\mathrm{u}}$ unsigned parts. Suppose that the incidence
parity on each part is the same for every color.
For $i\in[k_{\mathrm{s}}]$ and $c\in[\Delta]$, set
$s_i(c)=|A_i\cap V(M_c)|-|B_i\cap V(M_c)|$.
Suppose that integers $d_i$, $i\in[k_{\mathrm{s}}]$, satisfy,
for every $i\in[k_{\mathrm{s}}]$,
\[
 s_i(c)\equiv d_i\pmod2\quad(c\in[\Delta]),
 \qquad \sum_{c=1}^{\Delta}s_i(c)=\Delta d_i.
\]
Then $J$ has an edge $\Delta$-coloring with color classes
$M'_c$, $c\in[\Delta]$, satisfying the following properties.
\begin{enumerate}
\item Every $M'_c$ covers $U$ and has the same incidence
parity on each part as $M_c$.
\item For every $i\in[k_{\mathrm{s}}]$ and $c\in[\Delta]$,
$|A_i\cap V(M'_c)|-|B_i\cap V(M'_c)|=d_i$.
\item Any two color classes differ in size by at most
$2k_{\mathrm{s}}+k_{\mathrm{u}}$.
\end{enumerate}
\end{lemma}
\begin{proof}
We first make the signed differences equal to their prescribed values,
then balance the class sizes. We use the alternating-path argument
of~\cite[Lemma A1.1]{Folz2026}, with the details given below.
Fix two colors $c_1,c_2$. The subgraph with edge set
$M_{c_1}\cup M_{c_2}$ has nontrivial components that are
$(c_1,c_2)$-alternating paths or even cycles; along each of
them the edge colors alternate between $c_1$ and $c_2$. Every endpoint of a nontrivial path is outside
$U$. Represent each part by an auxiliary vertex, called its label,
and replace each such path by an edge joining its endpoint labels.
The resulting auxiliary multigraph, in which we allow loops, is the quotient
used in this proof. Switching a path means interchanging its two
colors. Switching a set of quotient edges means switching the
corresponding vertex-disjoint paths in $J$, and therefore preserves a proper
edge-coloring and coverage of $U$.
Every quotient degree is even, by the common part parities:
vertices covered by both colors contribute no path end, and
the number covered by exactly one color is even on each part.
At an endpoint in a signed part, record the side sign, positive
on $A_i$ and negative on $B_i$, multiplied by $+1$ if the
endpoint is covered by $c_1$ and by $-1$ if it is covered by $c_2$.
The sum of these signs at label $i$ is $s_i(c_1)-s_i(c_2)$.
Switching a path subtracts each recorded endpoint sign from the
corresponding $s_i(c_1)$ and adds it to $s_i(c_2)$.
Contributions from the two endpoints are added if they lie in
the same signed part.

First make the signed differences exact in the order $i=1,\ldots,k_{\mathrm{s}}$.
At part $i$, choose a coloring minimizing $\sum_c(s_i(c)-d_i)^2$,
subject to the part parities, coverage of $U$, and the differences
already fixed at the earlier parts.
The degree sums preserve $\sum_c s_i(c)=\Delta d_i$ under
every recoloring. If the minimum is positive, the mean and parity conditions give
colors $c_1,c_2$ with $s_i(c_1)\ge d_i+2$ and $s_i(c_2)\le d_i-2$.
At each earlier signed label, pair opposite signs; their
numbers agree. At label $i$, pair opposite signs as far as
possible. At least four positive ends remain and no negative
end remains. Pair incidences arbitrarily at every other label.
Start at one unpaired incidence and follow the pairings until
another unpaired incidence is reached. This gives a trail whose
two end incidences are positive and lie at label $i$. Switch
this trail. It transfers two units from $s_i(c_1)$
to $s_i(c_2)$, transfers zero at the earlier signed labels, and
changes each part incidence by an even number. The sum of
squares strictly decreases, a contradiction. Thus $s_i(c)=d_i$
can be achieved for every signed part and color.

With these differences fixed, minimize the sum of squared class
sizes. For two colors, pair opposite signs at every signed label
and pair incidences arbitrarily at every unsigned label. All
resulting trails are closed. Give each one a direction of traversal.
We split a trail whenever an unsigned label repeats, or a signed
label has two visits with the same incoming sign. To specify the
split, write $(e_1,f_1)$ and $(e_2,f_2)$ for the incoming and
outgoing incidences at the two visits. Replace these two pairs by
$(e_1,f_2)$ and $(e_2,f_1)$. The segment from $f_1$ to $e_2$
and the complementary segment from $f_2$ to $e_1$ now close
separately, partitioning the edges of the trail into two shorter
closed trails. At an unsigned label there is no sign restriction.
At a signed label, $e_1,e_2$ have the same sign, and $f_1,f_2$
both have the opposite sign. Thus each new closing pair still
joins opposite signs; all other pairs are unchanged.

Repeat these splits until none is possible. Each resulting trail
visits an unsigned label at most once and a signed label at most
twice, since three visits would include two with the same incoming
sign. It therefore has at most $2k_{\mathrm{s}}+k_{\mathrm{u}}$
quotient edges. Every visit uses two incidences, and at a signed
label their signs cancel. Switching the trail thus preserves all
part parities and signed differences, as well as coverage of $U$.
It changes either class size by at most $2k_{\mathrm{s}}+k_{\mathrm{u}}$,
since switching one alternating path changes a class size by at most one. The transfers over all
these trails sum to $|M_{c_1}|-|M_{c_2}|$, because alternating
cycles of $J$ contribute zero. If this difference exceeded
$2k_{\mathrm{s}}+k_{\mathrm{u}}$, some trail would have a positive transfer smaller than
the difference. Switching it would strictly decrease the sum of
squared class sizes, a contradiction. This proves the last assertion.
\end{proof}

\subsection{Walks avoiding a sparse set of edges or vertices}
We will use the following count when corrections of the palette weights
must avoid unavailable colors or previously used edges.

\begin{lemma}\label{pre:walk-stability}
Let $H$ be a graph in either setting of
Theorem~\ref{palette:thm:main}. Let $m$ be its order in the
ordinary case and the size of its larger side in the bipartite
case. Fix an integer $L\ge3$ and an admissible ordered pair
$(u,v)$. The following bounds count distinct $(u,v)$-walks of
length $L$ in $H$.
\begin{enumerate}
\item For every subgraph $D\subseteq H$, at most
$L\Delta(D)m^{L-2}$ such walks use an edge of $D$.
\item For every $W\subseteq V(H)$, at most
$(L-1)|W|m^{L-2}$ such walks have a vertex of $W$ at an
internal position.
\item The number of such walks with a repeated vertex, apart
from the prescribed equality of the endpoints when $u=v$,
is $O_L(m^{L-2})$.
\end{enumerate}
\end{lemma}
\begin{proof}
Fix a position at which the walk uses an edge of $D$. Choose the
vertices at all unspecified positions except one endpoint of this
edge. Once its other endpoint is fixed, there are at most
$\Delta(D)$ choices for the last endpoint. If the edge of $D$ is the first or last edge of the walk,
leave its internal endpoint to be chosen last.
Summing over the $L$ positions proves the first bound.
For the second, prescribe an internal position in $W$ and
choose the vertices at the other unspecified positions. For the last assertion, choose
two positions that are not both endpoints and require their
vertices to agree. This removes at least one free
choice, and there are $O(L^2)$ pairs of positions.
In the bipartite case the side of every position is fixed, so
each free position has at most $m$ choices.
\end{proof}

\subsection{A concentration inequality}
We use the following Chernoff bound.

\begin{lemma}[{\cite[Corollary~2.3 and Theorem~2.8]{JLR}}]
\label{bip:pre:chernoff}
Let $Z$ be a sum of independent indicator random variables, and
set $\mu=\mathbb E Z$.
For every $0<\varepsilon\le1$,
\[
 \mathbb P(|Z-\mu|\ge\varepsilon\mu)\le2\exp(-\varepsilon^2\mu/3).
\]
\end{lemma}

\addtocontents{toc}{\protect\columnbreak}
\section{Proof of Theorem~\ref{boundary:thm:main}}\label{sec:boundary}

We contract the parts while recording which incidences share an
original endpoint.
For some edges we temporarily give the incidences distinct endpoints,
which a final matching will restore. The circuit argument and
the choice of color classes below follow the method of
Folz~\cite{Folz2026}. We first give the results on the contracted
graph that will be used in the proof.

\subsection{Compatible circuit decompositions}

After contracting a part, we must still distinguish edges that
shared an original endpoint, since they cannot receive the same
color. Forbidden parts retain this information, and compatible
circuits will give matchings in the original graph.
Let $J$ be a graph whose vertex set is partitioned into parts and
whose edges join distinct parts. Contract each part to a vertex,
called its \emph{label}. We retain
parallel edges. At a label, the incidences belonging to one vertex
of $J$ form a \emph{forbidden part}. More generally, a
\emph{forbidden system} on a loopless multigraph consists of a
partition of the incidences at each vertex into forbidden parts.
A \emph{circuit} is a connected $2$-regular subgraph; a circuit
of length two is called a \emph{digon}. A circuit is \emph{compatible}
if its two incidences at each label belong to distinct forbidden
parts. For the system obtained from $J$, the edges corresponding
to a compatible circuit form a matching in $J$.
We say that a quotient subgraph uses a vertex of $J$ if one of
its edges corresponds to an edge incident with that vertex before
contraction. When edges are deleted, each forbidden part is restricted
to the remaining incidences, and empty forbidden parts are omitted.

In a compatible circuit decomposition, each cut edge from a
forbidden part must be accompanied by another cut edge outside
that part. This gives the following necessary condition.
At a fixed label, we identify an incidence with its edge when
intersecting a forbidden part with a cut. A forbidden system is
\emph{admissible} if
\[
 2|F\cap B|\le |B|
\]
for every edge cut $B$ and every forbidden part $F$. A cut is
nontrivial if both its sides are nonempty. A label is
\emph{unrestricted} if all its forbidden parts are singletons.
We use the following result.

\begin{theorem}[Fleischner and Frank {\cite[Theorem~2.1]{FF1990}}]
\label{boundary:thm:planar}
Every admissible forbidden system on a planar Eulerian loopless
multigraph has a compatible circuit decomposition.
\end{theorem}

Here and below, an Eulerian graph may be disconnected: every vertex
has even degree. The theorem then applies to each component containing
an edge.

In particular, this applies on at most five labels whenever at least
one pair of labels is nonadjacent, and always on at most four labels.
Indeed, $K_5$ with one edge removed is planar: draw a triangle and
let one of the two nonadjacent vertices lie inside the triangle and the other outside. Parallel copies
preserve planarity. Fan and Zhang~\cite{FZ} extend the planar theorem
to multigraphs with no $K_5$ minor. For a set $S$ of labels, write
$\partial_J(S)$ for the set of edges with exactly one label in $S$.
The quantity $d_B(x)$ continues to count incidences at the
original vertex $x$, before contraction.

On five labels, the following configuration may remain after
compatible circuits have been removed.
We call a forbidden system on $K_5$ \emph{exceptional} if its
edges are partitioned into two Hamilton cycles, every forbidden
part is contained in one of the incidence pairs defined by these
cycles, and every label has a forbidden part of size two. We call
the two Hamilton cycles complementary. Thus each label has either two forbidden parts
of size two or one such part and two singleton parts.

\begin{lemma}\label{boundary:lem:five-labels}
Every Eulerian loopless multigraph on at most five labels,
equipped with an admissible forbidden system, has an edge
decomposition into compatible circuits and at most one
exceptional system.
\end{lemma}
\begin{proof}
Repeatedly remove a compatible circuit whose removal preserves
admissibility. We show that any nonempty remaining graph $H$
is exceptional. Every compatible circuit is admissible, so
$H$ has no compatible circuit decomposition: any circuit in
such a decomposition could be removed. By
Theorem~\ref{boundary:thm:planar} and the preceding observation, $H$
has five labels with every pair joined by an edge.

We first show that a cut tight for a forbidden part must be the
singleton cut at its label. We then determine the forbidden
parts at each label and, finally, the edge multiplicities.

\smallskip\noindent\emph{Tight cuts.}
For a cut $B$ and a forbidden part $F$, call $|B|-2|F\cap B|$
their \emph{slack}; the cut is \emph{tight for $F$} when this is zero.
To rule out tight cuts other than singleton cuts, we first record an identity. Let $F$ be a forbidden part at a label
$v$ with $|F|=d_H(v)/2$. Its complement $F'$ in the incidences at
$v$ also satisfies all cut inequalities. Indeed, if $v\in T$,
then the slack $|\partial_H(T)|-2|F'\cap\partial_H(T)|$ equals
$|\partial_H(T\setminus\{v\})|
-2|F\cap\partial_H(T\setminus\{v\})|$.
At every label having a part of half its degree, combine all
other parts into its complement. Denote the resulting system by
$H'$. It is admissible. Its compatible circuits are compatible
in $H$, and a removal preserving admissibility in $H'$ also
preserves it in $H$.

We claim that a nontrivial cut tight for a part $F$ at a label $w$ must be
the singleton cut at $w$. Otherwise, choose the side $S$ not
containing $w$. Since the cut is not the singleton cut at $w$,
we have $1\le |S|\le3$. Set
\[
 A=F\cap\partial_{H'}(S),\qquad
 B=\partial_{H'}(S)\setminus A,\qquad |A|=|B|=h.
\]
Contract the complement of $S$ to a label $z$ and give $z$ the
two forbidden parts $A,B$. The inequalities for $A$ follow from
those for $F$, and the inequalities for $B$ follow from the
identity for complementary forbidden parts above. All other inequalities lift to $H'$.
This graph has at most four labels, so it has a compatible
circuit decomposition. Its circuits through $z$ give paths
through $S$, each joining an $A$-edge from $w$ to a $B$-edge
whose other end is an outside label $u$.

Replace each such path by an edge $wu$, keeping its $w$-incidence
in $F$ and its $u$-incidence in the part of its last edge. If
$u=w$, retain the lifted compatible circuit separately. Keep
the edges wholly outside $S$. The resulting outer graph is
Eulerian and admissible. To see the latter assertion, consider
a cut with side $T$ containing $w$. Let $t$ count the old
crossing edges avoiding $S$, and let $b$ count the $B$-edges
whose outside endpoint lies outside $T$. The new cut has size
$t+b$, equal to the old cut with side $T\cup S$; every part
except $F$ inherits its inequality from that cut. If $p$ counts
the old crossing $F$-edges avoiding $S$, the old $F$-inequality
on the cut with side $T$ gives
\[
 t+2h-b\ge2(p+h),
 \qquad\text{and hence}\qquad t+b\ge2(p+b),
\]
which is the required new inequality for $F$.
The outer graph also has at most four labels and therefore
has a compatible circuit decomposition. Every new edge belongs
to $F$ at $w$, so an outer compatible circuit uses at most one
new edge. As each replacement path has no repeated vertex, its lift cannot
repeat a label of $S$. Together with the circuits lying inside
$S$ and those retained separately,
these circuits decompose $H'$, a contradiction. This proves
the claim.

At a label of $H'$, either there are two equal forbidden parts
or every part has size strictly smaller than half the degree.
For a forbidden part $F$ at a label $v$, the singleton-cut
inequality is $2|F|\le d_{H'}(v)$. Removing a compatible circuit
preserves this inequality: at a label visited by the circuit,
the degree decreases by two. If there
are two equal parts, each loses one incidence. Otherwise, each
part had size at most half the old degree minus one, so even
an unchanged part satisfies the new bound. At an unvisited
label nothing changes. Every other nontrivial cut slack is an
even positive integer by the claim, and hence is at least two. A digon or triangle crosses a cut at most twice, so its
removal preserves all these inequalities as well. Consequently
$H'$ has no compatible digon or triangle.

\smallskip\noindent\emph{The forbidden parts at each label.}
For distinct labels $u,v$, call the set of edges joining them
the \emph{parallel class} $uv$. We say that this class varies at
an endpoint if its incidences there belong to more than one
forbidden part; otherwise it is constant at that endpoint.
All copies in this class lie
in one forbidden part at $u$ or in one forbidden part at $v$.
Indeed, if two copies belong to different parts at $u$, they
must belong to the same part at $v$; every third copy must then
belong to that same $v$-part. Thus a parallel class can vary
between forbidden parts at at most one endpoint.
We now choose one edge copy from each class of a triangle.
Call a vertex of this triangle a \emph{blocked corner} if both
entire incident classes lie in the same forbidden part. We show
that, if no corner is blocked, compatible edge copies can be chosen. At a corner where neither class varies, their
fixed parts are distinct, so any copies are compatible there.
If just one class varies, choose a copy from a part different
from the fixed part of the other class. This is possible because
a varying class meets at least two parts. If both classes vary,
choose a copy of one arbitrarily, then choose a copy of the other
from a different part.

These choices can be made at all three corners simultaneously.
Each class varies at at most one end, so its copy is chosen at
most once by these rules. At its other end, every copy belongs
to the same part, and the choice therefore cannot change
compatibility there. Finally choose arbitrary copies for classes
constant at both ends. The three chosen edges form the required
compatible triangle.

At each label, group the other four labels according to the
forbidden parts containing their entire parallel classes,
leaving a class that varies at that label as a singleton.
Every triangle is blocked by one of these groups. No group
has size four, since its forbidden part would contain every
incidence. Write $m_{uv}$ for the positive multiplicity of the
parallel class $uv$.

Suppose a group at $v$ contains three labels $x,y,z$ and
excludes $u$. The whole class $uv$ lies in a forbidden part
at one of its ends. Apply the inequality on the cut $\{v\}$
to this part and to the part containing $vx,vy,vz$. Together
they give $m_{uv}=m_{vx}+m_{vy}+m_{vz}$.
If the whole class $uv$ lay in one part at $u$, this cut
would be tight for that part, contrary to the claim. Hence
$uv$ is constant at $v$ and varies at $u$.
The excluded label $u$ has no group of size three. Otherwise,
its triple would have to exclude $v$,
forcing $uv$ to vary at $v$ as well.

Let $t$ be the number of labels with a group of size three. Two such
labels $v,w$ cannot exclude the same label $u$: the triangle $uvw$
would be blocked nowhere, since both classes at $u$ vary there.
Each excluded label blocks at most one triangle, because it has a
varying class and no group of size three. The $t$ labels with triples
and their $t$ distinct excluded labels are disjoint. A triple blocks
three triangles; an excluded label blocks at most one; and each of
the remaining $5-2t$ labels blocks at most two, one for each possible
pair. There are therefore at most $3t+t+2(5-2t)=10$ blocked corners.
All ten triangles have a
blocked corner, so each has exactly one.

If $v$ has a triple and excludes $u$, each triangle $uvx$ must be
blocked at $x$. The group containing $u,v$ there cannot also contain
$y\notin\{u,v,x\}$, since $vxy$ is already blocked at $v$. Thus for
any two of the other three labels $x,y$, neither can block $uxy$.
All three such triangles would have to be blocked at $u$, which
blocks at most one. This proves $t=0$. To attain ten blocked corners,
every label must now block two triangles. Its four neighboring labels
therefore form two pairs, and no parallel class varies there.
Moreover, each triangle has exactly one blocked corner. Give the two pairs
at each label signs $+1,-1$, and give an edge the product of
its endpoint signs. Every triangle has edge-sign product $+1$,
since it is blocked at exactly one corner. Fix one label and
reverse both signs at each other label whose edge to it is
negative. These reversals preserve the triangle products and
make every edge at the fixed label positive. Each remaining
edge is therefore positive as well. Thus the two incidence
signs of every edge agree. The parallel classes are then
colored red or blue, with two classes of each color at
every label. They form two complementary Hamilton cycles.
All copies of a class belong to its specified forbidden part
at both ends, and the two part inequalities at each label
give equal sums of multiplicities for the red and blue classes.

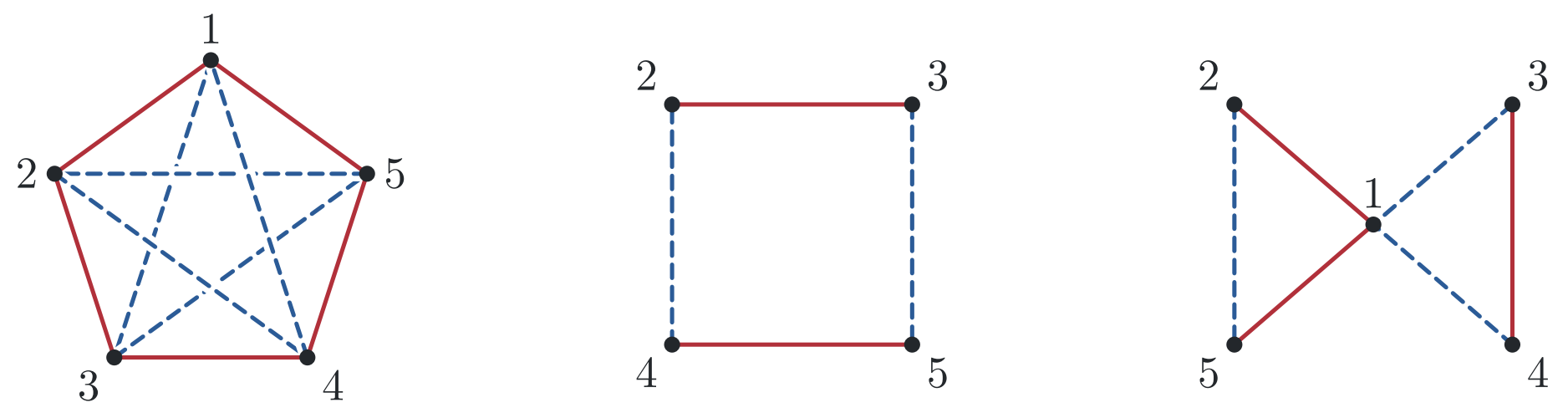
\begin{figure}[!htbp]
\centering
\begin{tikzpicture}[of diagram]
  \begin{scope}
    \foreach \i/\a in {1/90,2/162,3/234,4/306,5/18}
      \coordinate (v\i) at (\a:1.3);
    \foreach \i/\j in {1/3,3/5,5/2,2/4,4/1}
      \draw[of blue,of gap] (v\i)--(v\j);
    \draw[of red] (v1)--(v2)--(v3)--(v4)--(v5)--cycle;
    \node[of vertex,label=above:$1$] at (v1) {};
    \node[of vertex,label=left:$2$] at (v2) {};
    \node[of vertex,label=below left:$3$] at (v3) {};
    \node[of vertex,label=below right:$4$] at (v4) {};
    \node[of vertex,label=right:$5$] at (v5) {};
    \node[of panel] at (0,-1.95) {(a) $R$};
  \end{scope}
  \begin{scope}[xshift=4.6cm]
    % The same four-cycle is drawn without a crossing.
    \coordinate (u2) at (-.95,.95);
    \coordinate (u3) at (.95,.95);
    \coordinate (u5) at (.95,-.95);
    \coordinate (u4) at (-.95,-.95);
    \draw[of red] (u2)--(u3) (u5)--(u4);
    \draw[of blue] (u3)--(u5) (u4)--(u2);
    \node[of vertex,label=above left:$2$] at (u2) {};
    \node[of vertex,label=above right:$3$] at (u3) {};
    \node[of vertex,label=below right:$5$] at (u5) {};
    \node[of vertex,label=below left:$4$] at (u4) {};
    \node[of panel] at (0,-1.95) {(b) $C_1$};
  \end{scope}
  \begin{scope}[xshift=9.2cm]
    \coordinate (w1) at (0,0);
    \coordinate (w2) at (-1.1,.95);
    \coordinate (w5) at (-1.1,-.95);
    \coordinate (w3) at (1.1,.95);
    \coordinate (w4) at (1.1,-.95);
    \draw[of red] (w2)--(w1)--(w5) (w3)--(w4);
    \draw[of blue] (w2)--(w5) (w3)--(w1)--(w4);
    \node[of vertex,label=above:$1$] at (w1) {};
    \node[of vertex,label=above left:$2$] at (w2) {};
    \node[of vertex,label=below left:$5$] at (w5) {};
    \node[of vertex,label=above right:$3$] at (w3) {};
    \node[of vertex,label=below right:$4$] at (w4) {};
    \node[of panel] at (0,-1.95) {(c) $R-E(C_1)$};
  \end{scope}
\end{tikzpicture}
\caption{An exceptional system $R$ and its alternating four-cycle
$C_1=23542$ between $A_1=\{2,5\}$ and $B_1=\{3,4\}$.
In (a), solid red and dashed blue edges form the two
complementary Hamilton cycles. Deleting $C_1$ leaves the
triangles $1251$ and $1341$, which meet at label $1$.
In (a), crossings without a marked vertex are not vertices.}
\label{boundary:fig:exceptional}
\end{figure}

\smallskip\noindent\emph{The edge multiplicities.}
Let $A_i,B_i$ be the red and blue neighbor pairs at label $i$,
and let $r_i$ be the total multiplicity in either color at that label. The classes
between $A_i$ and $B_i$ form a compatible alternating four-cycle
$C_i$ as in Figure~\ref{boundary:fig:exceptional}(a),(b).
Write $m(A_i,B_i)$ for their total multiplicity.

Deleting a copy of $C_i$ must violate some cut inequality.
It removes one red and one blue incidence at each label it visits,
so it preserves the singleton inequalities. Every other cut
inequality has even slack at least two. Hence a
failure requires original slack two, all four edges of $C_i$
crossing the cut, and none belonging to the tested part.
At each label on $C_i$ its two edges use both forbidden parts.
Thus the tested part is at $i$, and the cut side avoiding $i$ is
$A_i$ or $B_i$. The part with no edge to that side has slack at
least six, since all ten classes are present. The other part has
slack $m(A_i,B_i)-r_i$, which must therefore equal two.

Every class occurs in two of the five four-cycles. Indeed, with
the red pentagon numbered cyclically from $0$ to $4$, the class
$01$ occurs in $C_2,C_4$ and $02$ in $C_3,C_4$. Rotations cover
the remaining classes. Since $r_i=d_{H'}(i)/2$, summing the five
slacks gives
\[
 10=\sum_i\bigl(m(A_i,B_i)-r_i\bigr)
    =2e(H')-e(H')=e(H').
\]
Each of the ten classes therefore has multiplicity one.
Consequently $H$ is a copy of $K_5$, and its original forbidden
parts refine the two cycle pairs. An unrestricted original
label could not have been coarsened and could not have two
pairs in $H'$. Thus every original label has a forbidden part of size two, and
$H$ is exceptional. Together with the circuits removed at the
start, this gives the required decomposition.
\end{proof}

\begin{lemma}\label{boundary:lem:absorb-five}
Let $R$ be an exceptional system on five labels, and let $C$
be a compatible circuit with $V(C)\subseteq V(R)$ and
$E(C)\cap E(R)=\emptyset$. Both inherit their forbidden parts
from the same system. Suppose that at least one of the following
conditions holds.
\begin{enumerate}
\item $C$ has length two or three.
\item $C$ contains a label $w$ at which neither forbidden part
used by $C$ contains an incidence of $R$.
\end{enumerate}
Then $R\cup C$ has a compatible circuit decomposition.
\end{lemma}
\begin{proof}
Color the two complementary Hamilton cycles of $R$ red and
blue. For a label $i$, let $A_i,B_i$ be its two red and its
two blue neighbors. The four edges between $A_i$ and $B_i$
form a compatible alternating four-cycle $C_i$. Deleting
$C_i$ from $R$ leaves the two triangles on $\{i\}\cup A_i$
and $\{i\}\cup B_i$; see Figure~\ref{boundary:fig:exceptional}(c).

Choose $i$ so that $C$ crosses both cuts with sides $A_i,B_i$.
If $e(C)\ge4$, take $i=w$; a circuit through $i$ on at least
four labels visits both neighbor pairs. If $C$ is a triangle,
choose $i\notin V(C)$, since its three labels meet both pairs.
If $C$ is a digon on $u,v$, choose $i$ so that $u,v$ belong
to different neighbor pairs. Such a label exists because every
parallel class belongs to two of the cycles $C_i$.

Set $H=(R\cup C)-E(C_i)$, which is Eulerian.
In $R$, temporarily combine all red incidences into one forbidden
part at each label, and all blue incidences into another. With these two parts, every cut inequality other
than the singleton inequality at the label of the part has slack
at least two. Indeed, a singleton cut at another label meets the
part at most once, and a two-label cut has six edges. The
four-cycle calculation in the preceding proof now shows that deleting
$C_i$ can violate only the red inequality on the cut with side
$A_i$ and the blue inequality on the cut with side $B_i$, both at
$i$. Their new slack is $-2$. Refining the parts of $R$ can only
increase the corresponding slacks. The compatible circuit $C$ has
nonnegative slack everywhere, and crosses both of these cuts at
least twice without using a part at $i$ that meets $R$. Thus $H$
is admissible for the original forbidden system.

At least one parallel class of $C_i$ remains absent from $H$.
A digon or triangle restores at most three such classes. When
$i\in V(C)$, at most $e(C)-2\le3$ edges of $C$ avoid $i$.
Thus $H$ is planar. Apply Theorem~\ref{boundary:thm:planar} to $H$
and restore $C_i$.
\end{proof}

\subsection{Correcting the boundary coloring}
We first prove the following more general recoloring statement.
Let $J$ be a graph with a partition $U,V_1,\ldots,V_k$ such that
every edge meets $U$ or joins two distinct sets $V_i$.
The set $U$ and a part consisting of private leaves are allowed
to be empty.
The vertices of $U$ are called exceptional. Every color class must cover $U$; on each remaining part we
prescribe the parity of the number of vertices it covers. For
$p=(p_1,\ldots,p_k)\in\{0,1\}^k$, a matching $M$ is
\emph{compatible with $p$} if it covers $U$ and
$|V(M)\cap V_i|\equiv p_i\pmod2$ for every $i$.
A coloring is compatible with $p$ if each color class is compatible.
For $I\subseteq[k]$, write $\mathcal V_I=\bigcup_{i\in I}V_i$.
For $T\subseteq U$, define $p(T,I)\in\{0,1\}$ by
\[
 p(T,I)\equiv |T|+\sum_{i\in I}p_i\pmod2.
\]
Keeping the index set $I$ in $p(T,I)$ avoids ambiguity when
one of the parts is empty.

\begin{proposition}\label{boundary:prop:recolor}
There exist a constant $0<\alpha<1$ and an integer
$\Delta_1\ge3$ for which the following holds.
Let $\Delta\ge\Delta_1$ be an integer, and let $J$ be a graph
with a vertex partition $U,V_1,\ldots,V_k$, where $1\le k\le5$.
Suppose that every edge has an endpoint in $U$ or joins two
distinct parts $V_i$. The set $U$ and a part consisting of
private leaves may be empty. Fix
$p=(p_1,\ldots,p_k)\in\{0,1\}^k$, and assume the following.
\begin{enumerate}
\item The degrees and sizes satisfy
\[
 \begin{gathered}
 d_J(u)=\Delta\quad(u\in U),\qquad |U|\le\alpha\Delta,\\
 d_J(v)\le\alpha\Delta\quad(v\in V(J)\setminus U),\qquad
 e(J)\le\alpha\Delta^2.
 \end{gathered}
\]
\item For each $i\in[k]$,
$\sum_{v\in V_i}d_J(v)\equiv\Delta p_i\pmod2$.
\item $J$ has an edge $\Delta$-coloring with at most two color
classes incompatible with $p$.
\item For every $T\subseteq U$ and $I\subseteq[k]$, set
$B=\partial_J(T\cup\mathcal V_I)$. Then
\begin{equation}\label{boundary:eq:cut-conditions}
 \begin{gathered}
 p(T,I)=1\quad\Longrightarrow\quad |B|\ge\Delta,\\
 p(T,I)=0\quad\Longrightarrow\quad
 |B|\ge2d_B(v)\quad(v\in V(J)).
 \end{gathered}
\end{equation}
\item If $k=5$, there is a simple graph $H$ of order $n$ and
maximum degree at most $\Delta$, with
$\Delta\ge2\lfloor(n-1)/8\rfloor$, and a vertex partition
$U,V_1,V_2,V_3,V_4$ such that
$\delta(H[V_i])\ge(1-\alpha)\Delta$ for every $i\in[4]$.
The graph $J$ is obtained by attaching $\Delta-d_H(v)$ distinct
private leaves to each $v\in V(H)$ and deleting the edges
of the four graphs $H[V_i]$. The part $V_5$ consists of these
leaves, and the targets satisfy
$p_i\equiv|V_i|\pmod2$ for $i\in[4]$ and
$p_5\equiv|V(H)|\pmod2$.
\end{enumerate}
Then $J$ has an edge $\Delta$-coloring compatible with $p$:
every color class covers $U$ and covers a number of vertices
congruent to $p_i$ modulo two in each $V_i$.
\end{proposition}

The cut conditions imply $|U|+\sum_i p_i\equiv0\pmod2$ by
taking $T=U$ and $I=[k]$. They also imply $p_i=0$ whenever
$V_i=\emptyset$.

\begin{proof}[Proof of Proposition~\ref{boundary:prop:recolor}]
Fix
\[
 \begin{gathered}
 N_{\mathrm{move}}=5\cdot10^5,\qquad \lambda_{\mathrm{sel}}=1/100,\\
 \lambda_{\mathrm{size}}=\frac{(\lambda_{\mathrm{sel}}/10)^{N_{\mathrm{move}}}}6,
 \qquad 0<\alpha<\frac{\lambda_{\mathrm{size}}\lambda_{\mathrm{sel}}}{10002}.
 \end{gathered}
\]
Then choose $\Delta_1$ sufficiently large. The constant
$N_{\mathrm{move}}$ bounds the number of vertices reassigned in Step~2, and
$\lambda_{\mathrm{size}}\Delta$ will be a lower bound on the number of selected colors.

Every color class covers $U$, since every vertex of $U$ has degree
$\Delta$.
We omit empty parts when applying Lemma~\ref{pre:balance} below.
If every color is compatible, there is nothing to prove.
Exactly one color cannot be incompatible, by the degree-sum congruence at each
part. We may therefore mark the two incompatible colors.
We keep them marked even if they become compatible when
parts are combined.

The difficulty is to correct the parity of each color at every part
without giving the same color to two edges incident with a vertex
before contraction. The
quotient graph records the parities, while its forbidden parts record
these possible conflicts. We use both descriptions throughout the proof.
Steps~1--2 prepare switches and a sufficiently large family of colors
whose crossing edges have small maximum degree. Steps~3--5 recolor
those crossing edges with the required parities on the combined parts.
Steps~6--7 restore the omitted edges and separate the combined parts.

We first replace certain endpoints by private leaves, so that the
alternating paths used for switches have bounded length.
For $u\in U$, call its edges to $V_i$ a large link if there are
at least $8\alpha\Delta$ of them. Give each edge in a large link
a private leaf in $V_i$ in place of its original end there, reusing
an end that already has degree one. Keep the edge colors and the prescribed part parities.
Every cut in~\eqref{boundary:eq:cut-conditions} has the same edges,
and the incidences of an original vertex have only been separated,
so all its cut inequalities persist. Throughout the proof,
labels record membership in a part, while forbidden parts retain
the incidences of individual vertices. Combining labels later
does not combine these forbidden parts.
Continue to write $J$ for this graph. At each $u\in U$ at most $(8k+1)\alpha\Delta\le41\alpha\Delta$
edges lie outside the large links. We will restore the original
ends after correcting the parities of the original parts.

\smallskip\noindent\textbf{Step 1. A switch between two parts.}
For $u\in U$ and $i\in[k]$, set
\[
 a_i(u)=\frac{|N_J(u)\cap V_i|}{\Delta}.
\]
For a real number $t$, write $t^+=\max\{t,0\}$. Define a
graph $G_{\mathrm{sw}}$ on $[k]$ by joining $i,j$ when
\begin{equation}\label{boundary:eq:switch-pairs}
 \sum_{u\in U}(a_i(u)a_j(u)-100\alpha)^+>4k^2.
\end{equation}
This graph depends only on $J$ and its original parts.

If $ij\in E(G_{\mathrm{sw}})$, we claim that in any proper edge $\Delta$-coloring two
specified colors $c_1,c_2$ can have their $i,j$ parities
toggled, with all other parities preserved. At any $u\in U$ with
$a_i(u)a_j(u)>100\alpha$, both links are large. Choose distinct
helper colors $c_3,c_4$ on these links, different from
the specified colors. Their other ends $x\in V_i$ and $y\in V_j$
are private leaves. For $t=1,2$, let $z$ be the other end of the
$c_t$-edge at $u$. If $z\notin U$, require $c_{t+2}$ to be
absent at $z$. If $z\in U$, require its $c_{t+2}$-edge to belong
to a large link. The first test excludes at most $\alpha\Delta$
choices and the second at most $41\alpha\Delta$.

There are initially at least $a_i(u)a_j(u)\Delta^2-4\Delta$
ordered pairs with the required locations. The two tests together
cost at most $82\alpha\Delta^2$ pairs. For large $\Delta$ the
proportion of all $(\Delta-2)(\Delta-3)$ helper pairs that succeed
at $u$ is therefore at least $(a_i(u)a_j(u)-100\alpha)^+$.
By~\eqref{boundary:eq:switch-pairs}, some fixed pair succeeds at
more than $3k^2$ vertices.

At a successful $u$, the $c_1,c_3$ component is a path
from $x$ through $u$ to the vertex $z$ in the first test or to its
private leaf in the second. It contains at most two vertices of $U$.
The $c_2,c_4$ path from $y$ has the same property.
There are at most $k^2$ ordered pairs of labels at the other ends of these paths, so four
successful vertices agree on this pair. Choose one of them.
Its two paths contain at most two other vertices of $U$, so
one of the other three successful vertices shares neither path.
For each color pair, switch the two distinct path components
through the chosen vertices. These two paths have the same pair
of endpoint labels, so their combined switches preserve every
part parity. At the first chosen vertex the
specified colors now form the isolated path $xuy$. Switching it
toggles exactly the required parities. Every path end lies outside $U$.

\smallskip\noindent\textbf{Step 2. Combining parts and selecting color classes.}
Combine the original parts in each connected component of
$G_{\mathrm{sw}}$. Delete edges between original parts belonging to
the same combined part, and call the resulting graph $J_0$.
Keep all edges meeting $U$. The target at a combined part is
the sum of its original targets modulo two.
Deleting an internal edge removes two incidences from that
part. Consequently the degree-sum congruences and the given
coloring pass to $J_0$, with the same two marked colors.
A cut formed from combined parts and a subset of $U$ is an
original cut with the same target and edge set. Thus
\eqref{boundary:eq:cut-conditions} passes to $J_0$ as well.

Once the vertices of $U$ are assigned to combined parts, call
an edge \emph{crossing} if its endpoints have different assigned
labels. The \emph{crossing degree} of a vertex is the number of
crossing edges incident with it. We seek a large family of unmarked
colors whose crossing edges have small maximum degree. More precisely, we need a positive
fraction of all colors, but their crossing degree must be a small
fraction of the number retained. This will leave enough colors to
insert whole circuits greedily in Step~5. We obtain this in three
stages: bound the total crossing degree at vertices where it is large,
retain colors with a common short record there, and move the vertices
specified by that record so that their recorded edges become internal.
Initially assign each $u\in U$ to the combined part containing
an original part $V_{i(u)}$ in which it has as many neighbors
as possible. Since at most $|U|$ of its neighbors lie in $U$,
\[
 a_{i(u)}(u)\ge\frac{1-\alpha}{k}\ge\frac1{2k}.
\]
Let $t(u)$ be the largest $a_j(u)$ for a label $j$ outside
this component of $G_{\mathrm{sw}}$, taking $t(u)=0$ if none exists.
Set $U_1=\{u\in U:t(u)\ge2000\alpha\}$.
For $u\in U_1$, choose $j$ attaining $t(u)$. The pair $i(u),j$
is not an edge of $G_{\mathrm{sw}}$, and
\[
 a_{i(u)}(u)a_j(u)-100\alpha
 \ge\frac{t(u)}{2k}-100\alpha
 \ge\frac{t(u)}{4k}.
\]
There are at most ten pairs of original labels. Summing
\eqref{boundary:eq:switch-pairs} over those pairs gives
\[
 \sum_{u\in U_1}t(u)\le4k\cdot10\cdot4k^2\le20000.
\]
Initially the crossing degree of $u$ is at most
$kt(u)\Delta+|U|$. Since $t(u)\ge 2000\alpha$
on $U_1$, the sum of these crossing degrees over $U_1$ is at most
$120000\Delta$. Every exceptional vertex outside $U_1$ has
crossing degree at most $10001\alpha\Delta$.

For each color record, at every $u\in U_1$ where its edge crosses,
the initial combined label of the other endpoint of that edge. The total number of
entries over the records is at most $120000\Delta$. Discard the
marked colors and those with more than $N_{\mathrm{move}}$ entries.
At least $\Delta/2$ color classes remain for large $\Delta$.

A record consists of pairs $(u,j)$: the edge of that color at $u$
crosses to a vertex with initial label $j$. A set of entries shared by several records therefore specifies
vertex moves that work simultaneously for all their colors.
The following maximization balances the number of retained colors
against the number of moves imposed on them.
We retain many records containing the same chosen entries; their
other entries need not agree. Maximality will ensure that each additional
entry occurs in only a small fraction of those records.

For a set $\mathcal A$ of entries contained in at least one remaining
record, let $n(\mathcal A)$ be the number of remaining colors whose
records contain $\mathcal A$. Choose $\mathcal A$ maximizing
$n(\mathcal A)(8/\lambda_{\mathrm{sel}})^{|\mathcal A|}$, and let $U_2$ be its
set of vertices. A record has at most one entry at each vertex,
so $|U_2|=|\mathcal A|\le N_{\mathrm{move}}$.
Comparison with the empty record gives
$n(\mathcal A)\ge(\Delta/2)(\lambda_{\mathrm{sel}}/8)^{N_{\mathrm{move}}}$.
At each $u\in U_1\setminus U_2$, at most
$5\lambda_{\mathrm{sel}}n(\mathcal A)/8$ of these colors cross. Otherwise
one of the at most five possible labels at $u$ would increase the
quantity maximized on adjoining that entry to $\mathcal A$. Indeed, for any
entry $(u,j)$ not already in $\mathcal A$, maximality gives
$n(\mathcal A\cup\{(u,j)\})\le(\lambda_{\mathrm{sel}}/8)n(\mathcal A)$.
Summing over the at most five labels gives the asserted bound.
Thus outside the vertices specified by $\mathcal A$, few retained colors
cross at any one vertex.

Among the $n(\mathcal A)$ colors discard those using an edge within $U_2$.
By simplicity this costs at most $\binom{N_{\mathrm{move}}}{2}$
colors. In each retained color class, the neighbor matched to a
vertex of $U_2$ now lies outside $U_2$. Its assigned label therefore
remains fixed when the vertices of $U_2$ are reassigned. For large $\Delta$ the number $s_0$ left satisfies
\[
 s_0\ge 5n(\mathcal A)/8,\qquad
 s_0\ge(\Delta/2)(\lambda_{\mathrm{sel}}/10)^{N_{\mathrm{move}}}.
\]
Reassign each vertex of $U_2$ to its recorded label. In every
retained color class, its matching edge now has both ends assigned
to that label. Its crossing degree in this family is therefore zero. At any
other vertex, only edges to $U_2$ can change their crossing status.
Since $J_0$ is simple, there are at most $|U_2|\le N_{\mathrm{move}}$
such edges. The crossing degree in the
family is therefore at most $\lambda_{\mathrm{sel}}s_0+N_{\mathrm{move}}$
on $U_1\setminus U_2$ and at most $10001\alpha\Delta+N_{\mathrm{move}}$
outside $U_1$.

Fix these assignments and set aside all edges of $J_0$ that
are now internal. Let ${p'}$ be the parity target on the
assigned combined parts. Its coordinate is the sum of their
original targets plus the number of assigned exceptional
vertices, modulo two. Each original color covers each
exceptional vertex once, and setting aside an internal edge
removes two incidences. Thus every unmarked crossing class
has target ${p'}$. The sum of its coordinates is even, since
it is congruent to $|U|+\sum_i p_i$. The crossing degree of
each assigned label is congruent to $\Delta$ times its target:
assigned exceptional vertices contribute $\Delta$ each, and
setting aside an internal edge removes two incidences.
Every union of assigned parts has the form $T\cup\mathcal V_I$, with
target $p(T,I)$ and the same cut as before. Hence
\eqref{boundary:eq:cut-conditions} holds for the crossing
graph with these adjusted targets.

Decompose the quotient of each surviving crossing color class
into edge-disjoint circuits and paths whose endpoints pair the
labels of odd degree. For example, remove circuits until a forest
remains, and then pair its incidences at each label except for
one incidence at each label of odd degree.
There are $r=|\{j:{p'}_j=1\}|/2\le2$ paths. If $r=2$, retain a most
frequent one of the three pairings of the four odd labels.
For $r\le1$, the pairing is unique. Let $\Pi$ be the chosen
pairing, $\mathcal R$ the retained color set, and $s$ its size.
Then
\[
 s\ge s_0/3\ge \lambda_{\mathrm{size}}\Delta.
\]
By our choice of $\alpha$ and since
$s_0\le3s$, the crossing edges of the full classes in
$\mathcal R$ have maximum degree at most $3\lambda_{\mathrm{sel}} s+N_{\mathrm{move}}$.

\smallskip\noindent\textbf{Step 3. Making the crossing graph Eulerian.}
For each odd label, introduce one auxiliary vertex assigned to
that label. For each pair in $\Pi$, join its two auxiliary vertices by $\Delta$
parallel auxiliary edges and assign one to each original
color. There are $r$ auxiliary edge types. An unmarked
augmented class is a matching before contraction and has even
degree at every label. Decompose its quotient into circuits.
The matching property makes every such circuit compatible. For a color of $\mathcal R$, close each
recorded path with its auxiliary edge and retain its recorded
circuits. Each type has $s$ such circuits containing no other
auxiliary type.

The full augmented quotient is Eulerian. Indeed, a label of
target ${p'}_j$ has original crossing degree congruent to
$\Delta{p'}_j$ and receives $\Delta{p'}_j$ auxiliary
incidences, so its new degree is even for either parity of
$\Delta$.

It is also admissible. Fix a cut, write $B$ for its original
crossing edges, and let $a\in\{0,1,2\}$ be the number of
auxiliary types crossing it. Since these types pair the odd
labels, $a\equiv p(T,I)\pmod2$ for the corresponding
union $T\cup\mathcal V_I$ of assigned parts. The augmented cut has
$|B|+a\Delta$ edges. If $a=0$, its target is even, and
\eqref{boundary:eq:cut-conditions} gives
$|B|\ge2d_B(v)$ for every nonauxiliary vertex before contraction;
no auxiliary vertex
has a crossing incidence. If $a=1$, the target is odd, so
$|B|\ge\Delta$ and the augmented cut has at least
$2\Delta$ edges. If $a=2$, its auxiliary edges alone give
at least $2\Delta$ edges. In the latter two cases, every
forbidden part has at most $\Delta$ incidences in the cut:
a nonauxiliary vertex has degree at most $\Delta$, and all edges
at an auxiliary vertex belong to one auxiliary type. Thus all
forbidden-part cut inequalities hold.

To apply Lemma~\ref{boundary:lem:five-labels}, we make the
union of the two full marked augmented classes admissible.
This union is Eulerian, since the full graph and all the unmarked
classes are Eulerian. A vertex before contraction contributes
at most two incidences. A failed inequality must therefore
have cut size two and slack $-2$, with both cut edges in the
same forbidden part. Two edges have at most two common
endpoints before contraction, so at most two vertices can fail
on one cut. There are at most fifteen nontrivial cuts up to
complementation. Thus at most thirty inequalities
fail.

For each failed inequality choose an entire unmarked augmented
color class of positive slack there. Such a class exists because
the full graph is admissible. An unmarked class splits into
compatible circuits, so all its slacks are nonnegative and even.
The chosen class thus contributes at least two to the failed
inequality. Add the distinct chosen classes to the two marked ones.
Their union $J^*$ is admissible and Eulerian, and uses at most
thirty classes besides the two marked ones.

Select every color of $J^*$ and every color of $\mathcal R$.
Selecting a color retains its entire crossing class, even when
only one of its circuits is needed below.
Lemma~\ref{boundary:lem:five-labels} decomposes $J^*$ into
compatible circuits and possibly one exceptional system $R$.
Classes outside $J^*$ keep their circuit decompositions.

\smallskip\noindent\textbf{Step 4. Removing an exceptional system.}
If $R$ remains and $r>0$, fix an auxiliary type. At least $s-30$
of its recorded circuits remain outside $J^*$, since a class uses
only one edge of that type. The incidences of the exceptional system use at most twenty
nonauxiliary vertices before contraction. Among the recorded circuits, at most
$20(3\lambda_{\mathrm{sel}}s+N_{\mathrm{move}})<s-30$ meet any
of these vertices, by the crossing degree of $\mathcal R$ and the
choice of $\lambda_{\mathrm{sel}}$. Choose a circuit whose
nonauxiliary vertices all avoid those of $R$. Auxiliary vertices
may be shared; the internal label used next has only
nonauxiliary incidences on this circuit. A digon or triangle satisfies
Lemma~\ref{boundary:lem:absorb-five} directly. In a longer circuit, deleting its unique auxiliary edge leaves
a path with an internal label. The two incidences at that label
come from nonauxiliary vertices not used by $R$, so the same
lemma applies.

Now suppose $r=0$. The remaining exceptional system has five
labels, so no original labels have been combined. Thus every
component of $G_{\mathrm{sw}}$ is a singleton and
the private-leaf hypothesis of the proposition applies. Write $W$
for $V_5$ together with the exceptional vertices assigned to it. If the
assigned part $W$ contains only leaves, it is unrestricted,
contrary to Lemma~\ref{boundary:lem:five-labels}. Thus $W$
contains an assigned original vertex.

Measure the five assigned parts in the original graph, before the
large links were given private ends. Each of the four parts other
than $W$ has even order, since its adjusted target is zero. Choose
a smallest such original set $S$, of order $a$. At least one
original vertex lies in $W$, so $4a\le|V(H)|-1$. Since $a$ is even,
\begin{equation}\label{boundary:eq:small-even-order}
 a\le2\lfloor(|V(H)|-1)/8\rfloor\le\Delta.
\end{equation}
The part $S$ contains an original part $V_i$, with $i\in[4]$,
so $a\ge(1-\alpha)\Delta$.

Let $Y\subseteq S$ consist of the original ends represented by the
four incidences of $R$ at this label, and set $n_Y=|Y|$. In the
modified graph the four incidences use two or three vertices, so
$n_Y\le3$. Let $D$ be $H$ with its original deficiency leaves attached.
Every original vertex has degree $\Delta$ in $D$, and
$\partial_D(S)$ has exactly the edges of the modified label cut.
Counting incidences at their original ends gives
\begin{align}
 |\partial_D(S)|-2d_{\partial_D(S)}(Y)
 &=\Delta(a-2n_Y)-2e_H(S\setminus Y)+2e_H(Y)\notag\\
 &\ge\Delta(a-2n_Y)-(a-n_Y)(a-n_Y-1)\notag\\
 &=(\Delta-a)(a-2n_Y)+a-n_Y(n_Y+1)\notag\\
 &\ge a-12>0.\label{boundary:eq:small-even-part}
\end{align}
Simplicity gives the second line. For the last line, use
\eqref{boundary:eq:small-even-order}, $n_Y\le3$, and
$a\ge(1-\alpha)\Delta>12$ for sufficiently large $\Delta$.

We use this inequality to find a circuit that can be combined with
$R$. The set $Y$ records the original endpoints already used by $R$
at the label representing $S$; we need a circuit whose two
incidences at this label have original endpoints outside $Y$.
The expression on the left of~\eqref{boundary:eq:small-even-part}
counts incidences outside $Y$ positively and incidences in $Y$
negatively. Its positivity forces the desired circuit, as follows.

The exceptional system contributes $-4$ to the left side of
\eqref{boundary:eq:small-even-part}. A compatible circuit
avoiding $S$ contributes zero. A circuit visiting $S$
contributes $2-2j$, where $j$ counts incidences at $S$ whose
original ends belong to $Y$. Some visiting circuit therefore has
both original ends there outside $Y$. Neither forbidden part used
by this circuit at the label contains an incidence of $R$, so it satisfies
Lemma~\ref{boundary:lem:absorb-five}. If this is a retained
circuit outside $J^*$, also select its original color class if necessary.
A circuit from the decomposition of $J^*$ already consists
entirely of selected edges.

In either case, $R$ and one additional circuit have a
compatible circuit decomposition. All other circuits remain
in place, so no exceptional system remains.

\smallskip\noindent\textbf{Step 5. Coloring the selected crossing edges.}
We now turn the compatible circuits into color classes. A compatible
circuit corresponds to a matching before contraction. All edges
of that matching can therefore receive one color. Circuits carrying auxiliary edges must be colored first:
these edges prescribe the odd labels of each color. The remaining
circuits change every label degree by an even number.

Let $q$ be the number of selected colors and let $J_1$ be their
crossing graph after auxiliary edges are omitted. Set
$d=\Delta(J_1)$. Each additionally selected full class is a
matching, so
\begin{equation}\label{boundary:eq:selected-degree}
 s\le q\le s+33,\qquad d\le3\lambda_{\mathrm{sel}} s+N_{\mathrm{move}}+33.
\end{equation}
The additional classes comprise the two marked classes,
at most thirty classes chosen to correct the cut inequalities,
and at most one further class in Step~4. For each auxiliary type,
at least $s-31$ recorded circuits containing only that type
remain intact. A compatible circuit uses at most one edge
of a type, since all copies have the same auxiliary endpoints.

If $r=1$, exactly $q$ circuits contain an auxiliary edge.
Give them distinct selected colors. If $r=2$, first give
distinct colors to circuits containing both types. Each type
has $q$ edges, and each circuit just colored uses one edge of
each type. Hence the two families containing only one type
have the same size $t$, with $t\ge s-31$.
Delete their auxiliary edges and view each remaining piece as a
matching in $J_1$. Each such matching covers at most eight vertices.
At any one of its vertices, at most $d$ pieces of the other family
can have an incident edge, since those pieces are edge-disjoint
subgraphs of $J_1$. Thus it meets at most $8d$ pieces in that family. By
\eqref{boundary:eq:selected-degree}, $16d<t$ for large
$\Delta$. Form a bipartite graph whose parts are these two
families of matchings, joining two matchings when they are
vertex-disjoint. Its minimum degree is greater than half its
part size. Hall's Theorem pairs the families so that each pair
has a union that is a matching. Give these pairs the remaining
selected colors.

In either case, each selected color uses exactly one edge of
each auxiliary type. Deleting the auxiliary edges leaves
target ${p'}$ in every selected color. If $r=0$, start with
all selected colors empty, which already gives target zero.
In all cases insert the remaining circuits greedily. Each
corresponds to a matching covering at most ten vertices and
hence excludes at most
$10d<q$ selected colors. A color is available for the entire
circuit, and insertion preserves every label parity.
Unselected crossing classes keep their original colors.
Thus the crossing graph is properly colored and every color class
has incidence parity ${p'}_j$ at each assigned label $j$.

\smallskip\noindent\textbf{Step 6. Restoring the internal edges.}
Every internal edge of $J_0$ meets $U$. First restore all such
edges outside the large links, using selected colors. Their total
degree is at most $41\alpha\Delta$ at a vertex of $U$ and at most
$\alpha\Delta$ elsewhere. The selected crossing edges have maximum
degree $d$; unselected crossing edges use none of the selected colors. Thus $2d+82\alpha\Delta<q$
by~\eqref{boundary:eq:selected-degree} and the choice of $\alpha$,
and greedy insertion suffices. Each remaining internal edge has a
private leaf at its nonexceptional end. At each $u\in U$, all other incident edges have been restored,
so the number of remaining edges is exactly the number of colors
missing at $u$. Their private ends impose no further restriction.
Use these missing colors in any order. Every restored edge has
both ends in one assigned part, so all assigned part parities
are preserved.

Every exceptional vertex is now covered by every color class.
Remove these vertices from the assigned combined parts and again
treat them as the separate set $U$; no graph edge is changed.
Each such vertex contributes one incidence per color, so this
subtracts exactly the target shift introduced by its assignment.
Every color now has the required target on the combined
parts. Restore the edges deleted when original parts were
combined. Each endpoint lies outside $U$ and has degree at
most $\alpha\Delta$ in $J$, so fewer than
$2\alpha\Delta<\Delta$ colors are excluded at any edge.
Greedy insertion preserves combined parity.

\smallskip\noindent\textbf{Step 7. Restoring the original parities.}
At this stage the coloring is proper and has the right parity on
each combined part. We must now correct the parity at each of
its original parts. The switch in Step~1 toggles two specified
colors at the two labels joined by an edge of $G_{\mathrm{sw}}$ and
preserves all other color--label parities.

Fix a rooted spanning tree in each component of $G_{\mathrm{sw}}$ and
process its nonroot labels from the leaves toward the root.
At the current label, the number of colors with wrong parity
is even, by the original degree-sum congruence. Pair those
colors and apply the switch in Step~1 along the edge to the
parent. This corrects the current label, changes only its
parent's parity in the specified colors, and preserves the sum
over the component. Every other color keeps its parity at
every original label, including when it is used as a helper color in Step~1.

The correct component parity now forces the parity at each root.
Apply Lemma~\ref{pre:balance} to the original parts, all unsigned,
requiring every color class to cover $U$. Each color has at most $e(J)/\Delta+k\le\alpha\Delta+k$
edges. Restore the large links one at a time. Remove the edges of the current link. Keep its set of $s_{\mathrm{link}}$ colors
at $u$ and assign them to the $s_{\mathrm{link}}$ distinct
original neighbors. Each neighbor excludes at most $\alpha\Delta$
colors, and each color occurs at at most $2(\alpha\Delta+k)$
vertices. Both numbers are smaller than $s_{\mathrm{link}}/2$ for large $\Delta$.
The bipartite graph of permissible color--neighbor pairs therefore
has minimum degree greater than $s_{\mathrm{link}}/2$. Hall's Theorem gives the
required bijection. It preserves every class size and the parity
on every original part, so the same bounds apply to the next link.
The resulting coloring of the original graph is compatible.
\end{proof}

\begin{proof}[Proof of Theorem~\ref{boundary:thm:main}]
Let $\alpha_0,\Delta_0$ be the constants in
Proposition~\ref{boundary:prop:recolor}. Take $\alpha=\alpha_0/6$
and choose $\Delta_1\ge\max\{\Delta_0,8,1/\alpha_0\}$.
Let $G,D,J,U,V_1,\ldots,V_k$ be as in the theorem, and let
$W$ be the set of private leaves. The degree condition gives
$n\le4\Delta+8\le5\Delta$. Every vertex of $U$ has degree
$\Delta$ in $J$, every other original vertex has degree at
most $\alpha\Delta$, and each leaf has degree one. Moreover,
\[
 e(J)\le |U|\Delta+n\alpha\Delta
          \le6\alpha\Delta^2=\alpha_0\Delta^2.
\]
Thus the size and degree hypotheses of the proposition hold.

Set $V_{k+1}=W$, take $p_i\equiv|V_i|\pmod2$ for $i\in[k]$,
and take $p_{k+1}\equiv n\pmod2$.
Then $|U|+\sum_{i=1}^{k+1}p_i\equiv0\pmod2$.
The required degree congruences follow from
\[
 \sum_{v\in V_i}d_J(v)=\Delta|V_i|-2e_G(V_i)\quad(i\in[k]),
 \qquad |W|=\Delta n-2e(G).
\]

We verify~\eqref{boundary:eq:cut-conditions}. For
$T\subseteq U$ and $I\subseteq[k+1]$, set
$B=\partial_J(T\cup\mathcal V_I)$. Replacing $T$ by $U\setminus T$
and $I$ by $[k+1]\setminus I$ preserves both $B$ and $p(T,I)$,
so we may assume $k+1\notin I$. Then $T\cup\mathcal V_I\subseteq V(G)$,
$p(T,I)\equiv|T\cup\mathcal V_I|\pmod2$, and
\[
 B=\partial_D(T\cup\mathcal V_I),\qquad |B|=\Delta|T\cup\mathcal V_I|-2e_G(T\cup\mathcal V_I).
\]
For odd $|T\cup\mathcal V_I|$, the absence of an overfull set gives
$|B|\ge\Delta$. If $|T\cup\mathcal V_I|$ is even, add or remove
an original vertex $v$ to obtain the odd set
$(T\cup\mathcal V_I)\oplus\{v\}$. Since $d_D(v)=\Delta$,
\[
 |B|-2d_B(v)
 =|\partial_D((T\cup\mathcal V_I)\oplus\{v\})|-\Delta\ge0.
\]
Indeed, changing the side of $v$ adds its $\Delta-d_B(v)$
noncrossing edges to the cut and removes its $d_B(v)$ crossing
edges. The resulting odd cut has size at least $\Delta$ by the
same odd-set inequality.
A private leaf has at most one incidence in $B$, and $|B|$
is even, so its inequality also holds.

Write $e=xy$ for the edge in the hypothesis. Extend an edge
$\Delta$-coloring of $G-e$ to the deficiency edges by giving them
distinct missing colors at their original ends. Every original
vertex is then covered in every color, except that $x,y$
each miss one color. Since $e$ lies inside a part, restriction
to $J$ colors all its edges and every color still covers $U$.
For any color missing at neither $x$ nor $y$, deleting internal
edges preserves its parity $|V_i|$ on each $V_i$, and matching
parity gives target $n$ on $W$. Thus at most two classes are
incompatible. When $k=4$, the additional
five-label hypothesis of Proposition~\ref{boundary:prop:recolor}
holds by the construction of $J$ and the minimum-degree
assumption. Applying the proposition to the $k+1$ parts gives
$|V(M_c)\cap V_i|\equiv|V_i|\pmod2$ for every color class
$M_c$ and $i\in[k]$. Thus each color is missing at an even
number of vertices of each $V_i$, as required.
\end{proof}

\subsection{Proof of Lemma~\ref{boundary:lem:regular-ambient}}
\label{boundary:sec:regular}
The omitted vertices in Lemma~\ref{boundary:lem:regular-ambient}
are retained when we measure the fifth part. This is where
regularity of the original graph enters. Requiring at least
thirteen omitted vertices ensures that this part is larger than
the loss of twelve in \eqref{boundary:eq:small-even-part}.

\begin{proof}[Proof of Lemma~\ref{boundary:lem:regular-ambient}]
Let $\alpha_0,\Delta_0$ be the constants in
Proposition~\ref{boundary:prop:recolor}. Take $\alpha=\alpha_0/7$
and choose $\Delta_1$ sufficiently large. Attach the deficiency
leaves to $G$, obtaining $D$, and denote their set by $W$.
Set $J=D-\bigcup_{i=1}^k E(G[V_i])$.
The degree bound gives $N\le5\Delta+8\le6\Delta$. Consequently,
$e(J)\le |U|\Delta+|V(G)|\alpha\Delta\le7\alpha\Delta^2$.
The other degree bounds, the targets
$p_i\equiv|V_i|\pmod2$ and $p_{k+1}\equiv|V(G)|\pmod2$,
the cut inequalities, and the coloring with at most two
incompatible classes follow exactly as in the proof of
Theorem~\ref{boundary:thm:main}. Their verification uses the
odd-set inequalities in $G$, not an order bound on $G$.

When $k\le3$, Proposition~\ref{boundary:prop:recolor} applies
directly. Suppose $k=4$. We use its proof, replacing only the
case $r=0$ in Step~4. This is the only place where the extra
five-label hypothesis uses the order bound on the original graph.
All other steps, including the case $r>0$, remain unchanged.

Suppose that an exceptional system $R_0$ remains and that $r=0$.
There are five labels, so no original parts have been combined.
For $i\in[4]$, let $U_i$ be the exceptional vertices assigned to
$V_i$, set $X_i=V_i\cup U_i$, and set
\[
 X_5=V(R)\setminus\bigcup_{i=1}^4X_i.
\]
The sets $X_1,\ldots,X_4$ have even orders, because their adjusted
targets are zero. Since $N$ is even, $X_5$ also has even order.
Each of the first four sets has more than twelve vertices by
its minimum-degree bound. Also $X_5$ contains
$V(R)\setminus V(G)$, so $|X_5|\ge13$.
Choose a smallest $X_i$, of order $a$. All five orders are even,
and hence
\[
 14\le a\le2\lfloor N/10\rfloor\le\Delta.
\]

If $i\le4$, use $S=X_i$ in the calculation
\eqref{boundary:eq:small-even-part}, with the completion $D$ of
$G$. Every vertex of $X_i$ has degree $\Delta$ in $D$ and
$G[X_i]$ is simple. That calculation gives positive slack and
therefore supplies the compatible circuit needed to absorb $R_0$.

It remains to handle $i=5$. Let $U_5=U\setminus\bigcup_{i=1}^4U_i$,
and write $X=V(R)\setminus X_5=\bigcup_{i=1}^4X_i$.
The fifth label represents the part $U_5\cup W$ in the boundary
graph. Let $B$ be the cut of this part. Giving large links
private ends did not change the edges of this cut. Counting
the deficiency edges at vertices of $X$ gives
\begin{align}
 |B|&=e_G(X,U_5)+\sum_{v\in X}(\Delta-d_G(v))\notag\\
    &=|\partial_R(X_5)|+2\bigl(e_R(X)-e_G(X)\bigr)
      \ge |\partial_R(X_5)|.\label{boundary:eq:ambient-cut}
\end{align}
Here $X_5$ is a set in the original regular graph, whereas the
fifth label in the boundary graph is $U_5\cup W$.
Every edge of $R$ between $X$ and $X_5$ contributes one edge
to $B$, either an edge of $G$ or a deficiency edge at its end
in $X$. Each edge of $R[X]$ missing from $G$ contributes two
deficiency edges to $B$. This explains both the equality and
the nonnegative extra term in~\eqref{boundary:eq:ambient-cut}.

The four incidences of $R_0$ at the fifth label use at most three
vertices before contraction. Let $Y\subseteq U_5$ be the original vertices
among them, set $t=|Y|$, and let $n_{\mathrm{leaf}}$ be the number of leaves
among them. Then $t+n_{\mathrm{leaf}}\le3$. Every such leaf has one incidence
in $B$, and $d_B(Y)\le d_{\partial_R(X_5)}(Y)$.
Using~\eqref{boundary:eq:ambient-cut} and counting in the simple
regular graph $R$, we obtain
\begin{align*}
 |B|-2d_B(Y)-2n_{\mathrm{leaf}}
 &\ge |\partial_R(X_5)|-2d_{\partial_R(X_5)}(Y)-2n_{\mathrm{leaf}}\\
 &=\Delta(a-2t)-2e_R(X_5\setminus Y)+2e_R(Y)-2n_{\mathrm{leaf}}\\
 &\ge (\Delta-a)(a-2t)+a-t(t+1)-2n_{\mathrm{leaf}}\\
 &\ge a-12>0.
\end{align*}
For the last inequality, use $a\le\Delta$, $a\ge14$,
and $t(t+1)+2n_{\mathrm{leaf}}\le12$ when $t+n_{\mathrm{leaf}}\le3$.

This is the same positive slack used in Step~4, now counting
the original vertices and the leaves used by $R_0$ at the fifth
label. The exceptional system contributes $-4$. A compatible
circuit through the label contributes $2-2j$, where $j$ is the
number of its incidences using one of those vertices or leaves.
Some circuit must therefore avoid all these vertices and leaves
at this label.
Lemma~\ref{boundary:lem:absorb-five} absorbs $R_0$ with that circuit.
Steps~5--7 of the proposition now complete the coloring, and
the conclusion follows as in Theorem~\ref{boundary:thm:main}.
\end{proof}

\section{Proof of Theorem~\ref{palette:thm:main}}\label{sec:palettes}
Following the absorption approach in~\cite{Folz2026}, we set aside
a sparse graph that can be decomposed together with the eventual
uncovered remainder. We specify the required degree and parity
conditions on that remainder below. A weighted random greedy
procedure and a final greedy step cover the other edges by
edge-disjoint allowed triangles.

Let $H$ and the palettes $\mathcal P(v)$ be as in
Theorem~\ref{palette:thm:main}. If a bipartite side is empty,
then $H$ has no edges and $|\mathcal P(v)|=d_H(v)$ makes every palette
empty, so the conclusion is immediate. We may therefore assume
that both sides are nonempty in the bipartite case.
Introduce a set $\mathcal C$ of color vertices, disjoint from
$V(H)$ and indexed by $[\Delta]$. Let $G$ be obtained from $H$
by adding the edge $vc$
for each $c\in \mathcal P(v)$. We call the edges of $H$ \emph{original edges}
and the other edges \emph{color edges}. An \emph{allowed triangle}
is a triangle $uvc$ with $uv\in E(H)$ and
$c\in \mathcal P(u)\cap \mathcal P(v)$.
It suffices to decompose $G$ into allowed triangles: the triangle
$uvc$ assigns color $c$ to $uv$, and each color edge $vc$ is used
exactly once. The resulting coloring is proper and has palette
$\mathcal P(v)$ at every vertex $v$.

Every allowed triangle uses one original edge and two color
edges. The following degree and parity conditions are therefore
necessary for a decomposition into allowed triangles.
A subgraph of $G$ is \emph{balanced} if, within that subgraph,
the original-edge and color-edge degrees agree at every vertex
of $H$, and every color vertex has even degree in the ordinary
case, or equal degrees into the two sides in the bipartite case. The palette assumptions imply
that $G$ is balanced. An allowed triangle is balanced.
Edge-disjoint unions of balanced graphs, and differences of nested
balanced graphs, are also balanced.

Counting walks by their first edge gives
$\delta(H)\ge bm$. The palette conditions give
$bm\le\Delta\le(1+\eta)m$. In the bipartite case, both sides have
order at least $bm$. Thus, if
$\eta\le1/2$, every side under consideration and $\mathcal C$ have order
$\Theta(m)$. Constants implicit in this section
may depend on $b$ and $L$ and on other fixed constants introduced
before $\eta$ is chosen.

\subsection{A family of triangles and digons}
The binary cube $\{0,1\}^k$ has an edge between two points when
they differ in exactly one coordinate. The following operation
redistributes sets of labels by exchanges along these edges. We use it first to
balance the triangles selected for a given remainder, and then to pair
all unused arc copies between each pair of auxiliary vertices.

\begin{lemma}\label{palette:lem:cube-exchange}
Let $k,r$ be nonnegative integers. Assign an initial set and a
target set of $r$ labels to each vertex of the binary cube
$\{0,1\}^k$. Suppose that, for every label, its total number of
initial occurrences equals its total number of target occurrences.
Then the initial sets can be transformed into the target sets by
exchanges with the following properties.
\begin{enumerate}
\item There are $2k$ steps, using the coordinates in the order
$1,\ldots,k,k,\ldots,1$. Each step uses only edges in its
specified coordinate.
\item After every step, each vertex carries exactly $r$
distinct labels.
\item At each step, equally many labels move in the two
directions along every edge that is used.
\end{enumerate}
\end{lemma}
\begin{proof}
We prove the assertion by induction on the dimension. The main task
is to divide the occurrences of every label between the two halves
so that each half has the same initial and target counts. Once this
is done, induction rearranges the labels within the halves.

For $k=0$ the assertion follows from the occurrence counts. For
$k\ge1$, pair the points differing in their first coordinate. We
divide the labels in each initial and target pair between the two
halves of the cube. Use a directed graph with three vertex classes:
one vertex for each initial pair, one for each label, and one for
each target pair. Join an initial-pair vertex to each label occurring
in that pair, and join each label to every target-pair vertex where
it is required. Orient the arcs from initial pairs to labels and
from labels to target pairs. Give each arc weight $1$ or $2$,
equal to the number of occurrences of its label in its pair. Every pair has total weight
$2r$, and the weights into and out of each label agree.

Split the weight vector into two nonnegative integer vectors as follows. The edges of
odd weight form an undirected graph in which every degree is even: a pair vertex has total
weight $2r$, and at a label vertex the total incident weight is twice
its incoming weight. In either case the number of odd incident
weights is even. Orient them along Euler
tours and let $\gamma$ be $+1$ or $-1$ according as this orientation
agrees with the given direction or opposes it. Set $\gamma=0$ on the other
edges. At every vertex, the sum of $\gamma$ on outgoing edges
equals its sum on incoming edges. For the weight
vector $\omega$ the two vectors $(\omega+\gamma)/2$ and $(\omega-\gamma)/2$ have entries $0$
or $1$. Each gives total $r$ at every pair and equal incoming and
outgoing weights at every label.

Use the two vectors for the two values of the first coordinate.
In each initial pair move the labels to their assigned points. The
number crossing each way is the same because both points keep $r$
labels. A label present at both points stays at both. Make the same
preparatory exchanges for the target pairs. Within each half, all
label counts now agree between the initial and target assignments.
Apply induction there in parallel, then undo the preparatory exchanges
for the targets. This gives the stated sequence of coordinates.
\end{proof}

For a directed multigraph, \emph{balanced} means that the indegree
equals the outdegree at every vertex. Its total degree counts both
incoming and outgoing arcs, with multiplicity; $\Delta$ denotes the
maximum of these total degrees. A directed triangle
is a directed cycle of length three; a directed digon consists of
two opposite arcs between the same pair of vertices.

\begin{lemma}\label{palette:lem:circulation}
Fix $0<\beta<1$. For all sufficiently large $m$, every oriented
graph $D$ on at most $4m$ vertices with
$\Delta(D)\le m^{1-\beta}$ has the following extension.
There are a set $W$ disjoint from $V(D)$, a directed multigraph
$R$ on $V(D)\cup W$, and a family $\mathcal F$ of directed
triangles and digons in $D\cup R$ satisfying these properties.
\begin{enumerate}
\item $|W|=o(m)$, every arc of $R$ has an endpoint in $W$, and
$\Delta(D\cup R)=O_\beta(m^{1-\beta/4})$.
\item Every arc of $D\cup R$ lies in at most four members of
$\mathcal F$. Each arc of $D$ lies in exactly one triangle
of $\mathcal F$, whose third vertex belongs to $W$. Every
member of $\mathcal F$ other than these triangles contains at
most one vertex of $D$.
\item For every balanced directed subgraph $J\subseteq D$,
some members of $\mathcal F$ form an edge decomposition of
$R\cup J$.
\end{enumerate}
\end{lemma}
\begin{proof}
Call the collection of auxiliary arcs between a fixed pair of
vertices a \emph{link}. There are two requirements on the construction.
For any balanced $J$, the triangles using its arcs must be completed by auxiliary
triangles so that every auxiliary link is used equally often in its
two directions. The unused copies on each link must then be paired
into digons chosen from a fixed list. The first requirement is a
redistribution problem for labels; the second is a matching problem
for arc copies. The cube lemma supplies both constructions.

\smallskip\noindent\textbf{Constructing the triangles.}
Choose a binary cube $Q_k$ on new vertices, of order $N=2^k$ with
$N\le m^{1-\beta/2}<2N$, and distribute $V(D)$ among its points as
evenly as possible by a map $g$. For every arc $u\to v$ make a
request from $g(v)$ to $g(u)$. The bipartite multigraph of requests
on source and destination copies of $Q_k$ has maximum degree at most
$r_0=\lceil4m/N\rceil\lceil m^{1-\beta}\rceil=O(N)$. Properly color
its edges with $r_0$ labels. This follows by padding to a $r_0$-regular
bipartite multigraph on the same equal sides and successively
removing perfect matchings by Hall's theorem.

Introduce $r_0$ additional labels, with a fixed ordering, and
$2r_0$ new vertices $w_i$, one for each original or additional label.
For an arc $u\to v$ with request label $i$, and for each arc
$x\to y$ in the following list, include the directed triangle
$x\to y\to w_i\to x$:
\[
 u\to v,\qquad v\to g(v),\qquad g(u)\to u.
\]
For every label $i$ and every directed cube edge $x\to y$, list two
copies of $x\to y\to w_i\to x$. Use private arc copies in every
triangle, except that a triangle on $u\to v$ uses that arc of $D$.
Thus the listed triangles are edge-disjoint before the alternative
digon decompositions are introduced. The private arc copies distinguish
different uses of the same ordered pair of vertices.

\smallskip\noindent\textbf{Selecting triangles for a given remainder.}
Consider first which triangles to select for a balanced $J\subseteq D$.
Select its three triangles per arc and inspect the remaining imbalance.
At a vertex of $D$, the arcs to and from $w_i$ cancel within each
group of three. The arcs between $v$ and $g(v)$ balance because $J$
is balanced at $v$. At the cube points, a request from $x$ to $y$
has used one arc $x\to w_i$ and one arc $w_i\to y$. Moving label $i$ from $x$ to $y$
will cancel these two contributions. If $x=y$, they already
cancel and no move is needed.

At each point, let the labels of its selected outgoing requests be the initial labels,
and let the labels of its incoming requests be the target labels. The sets have no repeated label because the request edges were
properly colored. Each label has equal global counts,
and the two sets at a point have equal size by summing the balances of
$J$ over the set $g^{-1}(x)$ of vertices assigned to that point $x$. If this size is $t$, add the first $r_0-t$ additional
labels to both sets. Every point now has $r_0$ labels. Apply
Lemma~\ref{palette:lem:cube-exchange}. For a move of $i$ along
$x\to y$ select the triangle $x\to y\to w_i\to x$. Two copies suffice:
the coordinate of the edge is used twice and its tail has at most one
occurrence of $i$ at each step. Across each cube edge the opposite
moves balance. At a point $x$ the moves of label $i$ contribute net
usage of $w_i\to x$ equal to the indicator that $i$ is initially
present at $x$ minus the indicator that it is required there. This cancels exactly the arcs contributed by the selected requests, including
the zero difference for the additional labels. Thus on every auxiliary
link the selected triangles use its two directions equally often.

\smallskip\noindent\textbf{Covering the unused arc copies.}
Equality of the two directional counts would let us pair the unused
copies arbitrarily. We cannot list every possible pair, since a
single arc would then belong to too many listed digons. Instead we
use successive cube exchanges. Each position has only two possible
next positions, so each arc copy belongs to only a bounded number
of digons in the fixed list.

It remains to prepare, before $J$ is chosen, the digons that cover
the unused copies on a link. Pad its two directions to equal size $s$
and then to $n=2^h$, where $h=\max\{1,\lceil\log_2s\rceil\}$, so
$n\le2s$. Let the copies of one direction be the inputs of a binary
$h$-cube and the opposite copies be its outputs. Between them use the
$2h$ coordinate steps of Lemma~\ref{palette:lem:cube-exchange}.
At each intermediate time and point add two arc copies. The slot names
refer to successive time steps on this one link: there is an incoming
slot of the output direction and an outgoing slot of the input
direction. List the digon joining these two slots, called the
\emph{local digon}. Also list a \emph{transition digon}
for every permitted transition, joining its outgoing slot to the
incoming slot one step later. A transition either keeps the cube point fixed or changes its
current coordinate. The positions $(t,x)$, with
$t\in\{0,\ldots,2h\}$ and $x\in\{0,1\}^h$, are the vertices of a
layered directed graph; the permitted transitions from time $t$
to time $t+1$ are its arcs. Inputs supply the first outgoing slots and
outputs the last incoming slots.
After removing the copies used by the selected triangles, mark the
surviving inputs and outputs active and the others inactive. The
active counts agree because the selected triangles used equally
many copies in the two directions on this link.
Figure~\ref{palette:fig:slots} shows these pairings. Its nodes
represent arc copies on one fixed auxiliary link, not new vertices
of the directed multigraph. Each drawn edge represents a listed
digon, pairing copies with opposite directions.

\begin{figure}[!htbp]
\centering
\begin{tikzpicture}[of diagram]
\foreach \x/\t in {0/t,4/{t+1}} {
 \node at (\x+.6,2.2) {time $\t$};
 \foreach \y/\q in {1.2/0,0/1} {
  \node[of vertex,label=above:{$\mathrm{in}$}] (i\x\q) at (\x,\y) {};
  \node[of vertex,label=above:{$\mathrm{out}$}] (o\x\q) at (\x+1.2,\y) {};
  \draw[of red] (i\x\q)--(o\x\q);
 }
}
\draw[of blue,solid] (o00)--(i40);
\draw[of blue,solid] (o01)--(i41);
\draw[of blue,dashed] (o00)--(i41);
\draw[of blue,dashed] (o01)--(i40);
\node[anchor=east] at (-.3,1.2) {$0$};
\node[anchor=east] at (-.3,0) {$1$};
\node at (2.6,-.8) {red: local digons; blue: transition digons};
\end{tikzpicture}
\caption{Two positions differing in the current cube coordinate, at
consecutive intermediate times. A transition may stay at its position
(solid) or flip that coordinate (dashed). At each intermediate
position on a selected directed path, use the arriving and departing
transition digons. At every other intermediate position, use the
local digon. In either case both slots are covered exactly once.}
\label{palette:fig:slots}
\end{figure}
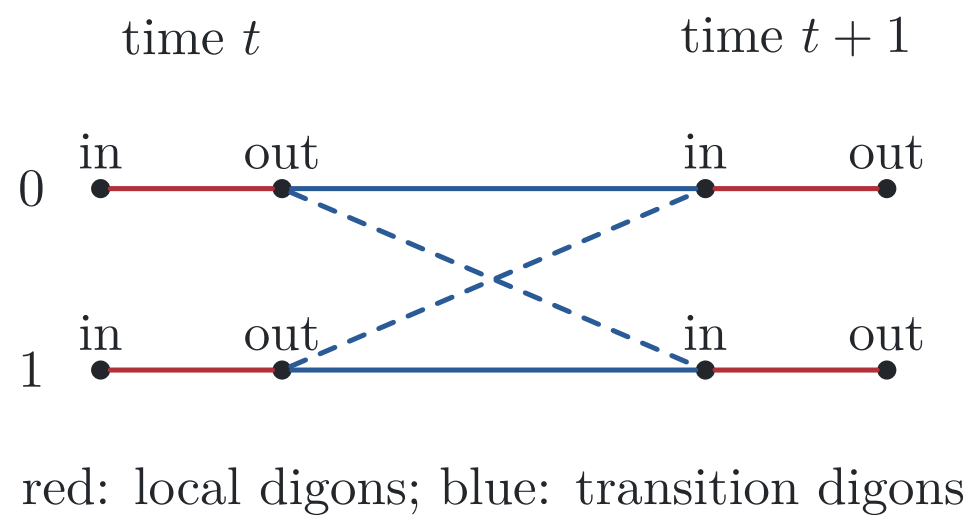

Apply the same cube lemma with $r=1$ and labels active and inactive.
Each step permutes the points within each pair. The active inputs
therefore determine vertex-disjoint directed paths in the layered
graph, ending at precisely the active outputs. Select the transition
digons corresponding to the arcs of these paths. At each unused
intermediate position select the digon between its two slots. To check exact coverage, consider each kind of slot separately.
Every surviving input is the initial vertex of one selected
directed path, and every surviving output is the terminal vertex
of one. At an intermediate position on such a path, the
incoming slot is covered by the arriving transition digon and the
outgoing slot by the departing transition digon. At every other
intermediate position, the local digon covers its two slots.
Deleted inputs and outputs were already covered by the selected
triangles. Thus every arc copy is covered exactly once. Applying the
construction on all links proves the partition assertion.

\smallskip\noindent\textbf{Size and multiplicity bounds.}
The number of new vertices is $N+2r_0=O(N)=o(m)$. Before adding the
arc copies used for digons,
an old vertex has $O(\Delta(D))$ incidences. A cube point has $O(r_0)$
request incidences and $O(r_0k)$ incidences from cube triangles. A label occurs on at most $N$ requests, and its vertex $w_i$
has $O(Nk)$ incidences in all. The slot construction multiplies the number of arc
copies on a link by $O(\log(2s))$ and adds no vertices. Thus the
maximum degree is $O(N\log^2m)=O_\beta(m^{1-\beta/4})$. Each slot
is in at most two transition digons and one digon at its position,
and each arc is in at most one triangle. Only triangles on arcs of
$D$ can contain two old vertices. This proves the remaining claims.
\end{proof}

\subsection{An absorber for the prescribed palettes}
The uncovered edges will lie in a sparse graph $D$, but the
precise remainder is not known when the construction starts.
We therefore reserve a graph that can complete the decomposition
for every balanced remainder.
Let $D\subseteq G$. An \emph{absorber for $D$} is a balanced
subgraph $R\subseteq G-E(D)$ such that $R\cup J$ has a
decomposition into allowed triangles for every balanced
$J\subseteq D$.

\begin{lemma}\label{palette:lem:absorber}
Fix $b>0$, an odd integer $L\ge3$, $0<\alpha<1$, and $\lambda>0$.
There exist a constant $\eta_0>0$ and a positive integer $m_0$,
depending only on $b,L,\alpha,\lambda$, for which the following holds. Let $H$ and its palettes satisfy
conditions~(1)--(3) of Theorem~\ref{palette:thm:main} with these
$b,L$ and with $0<\eta\le\eta_0$ in the palette bounds.
Let $G$ be the augmented graph defined above, and suppose that
$m\ge m_0$.
For every subgraph $D\subseteq G$ with
$\Delta(D)\le \lambda m^{1-\alpha}$, there is a subgraph $R\subseteq G$
satisfying the following properties.
\begin{enumerate}
\item $R$ is balanced, $E(R)\cap E(D)=\emptyset$, and
$\Delta(R)=o(m)$.
\item For every balanced subgraph $J\subseteq D$, the graph
$R\cup J$ has a decomposition into allowed triangles.
\end{enumerate}
The graph $R$ is chosen independently of $J$ and is an absorber
for $D$.
\end{lemma}
\begin{proof}
The absorber must be chosen before the remainder $J$ is known.
We therefore first express every possible remainder using members
of one fixed family of small balanced graphs. For each member we
build a graph with two triangle decompositions: one that includes
that member and one that does not. Once $J$ is known, we choose
between these two decompositions according to which members occur
in the chosen partition.

We construct $D_{\mathrm{aux}}\subseteq G-E(D)$ and a fixed family of balanced
graphs $J_f\subseteq D_{\mathrm{aux}}\cup D$, each of bounded order, so that
some of the $J_f$ form an edge decomposition of $D_{\mathrm{aux}}\cup J$ for every balanced
$J\subseteq D$. For each $J_f$, we then construct a graph $R_f$
such that both $R_f$ and $R_f\cup J_f$ decompose into allowed
triangles. We choose the graphs $R_f$ to be pairwise edge-disjoint,
with no edges in $D_{\mathrm{aux}}\cup D$, and keep $\Delta(D_{\mathrm{aux}}\cup\bigcup_f R_f)=o(m)$.
The required absorber will be $D_{\mathrm{aux}}\cup\bigcup_f R_f$.
We verify its decomposition property after constructing these graphs.

\smallskip\noindent\textbf{Step 1. Representing balanced remainders by trails.}
To construct the $J_f$, we apply Lemma~\ref{palette:lem:circulation}
using three formal classes $A,B,\mathcal C$ in cyclic order. In the
bipartite case, take the two sides and the colors as these classes,
and orient the edges of $G$ as $A\to B\to \mathcal C\to A$. In the ordinary case,
replace each vertex of $H$ by two formal copies. The edge $uv$
has alternatives $u_Av_B$ and $v_Au_B$, directed from $A$ to $B$.
The color edge $vc$ has alternatives $v_Bc$ and $cv_A$, with
the displayed directions. Let $D^*$ consist of all alternatives
of the edges of $D$.

Every balanced $J\subseteq D$ has a balanced lift
$J^*\subseteq D^*$ using one alternative of each edge. In the
bipartite case, its prescribed orientation is already balanced:
the two types of degree agree at every vertex of $H$, and each
color has equal degrees into $A$ and $B$. In the ordinary case,
pair each original-edge incidence with a color-edge incidence at
the same vertex of $H$, and pair the incidences at each color vertex. The resulting closed
trails consist successively of one original edge and two color
edges. Orient each trail. Lift the tail and head of each original
edge to $A$ and $B$, respectively; the two following color edges
then run from $B$ through $\mathcal C$ to $A$. At each formal
vertex, the pairing matches an incoming edge with an outgoing
edge. The lift is therefore balanced and uses one alternative
of every edge of $J$, including its color edges.

The graph $D^*$ has at most $2m+\Delta\le7m/2$ vertices
and maximum degree $O(m^{1-\alpha})\le m^{1-\alpha/2}$ for
sufficiently large $m$. Apply Lemma~\ref{palette:lem:circulation}
with $\beta=\alpha/2$, and denote its auxiliary directed multigraph
by $D_{\mathrm{aux}}^*$. Identify the two copies of each vertex
of $H$, and map the new vertices injectively into $\mathcal C$. This is
possible because their number is $o(m)$ and $|\mathcal C|=\Theta(m)$.
At most two formal vertices are mapped to any one vertex of $G$: a vertex
of $H$ has its two copies, and a color may also be the image of
one new vertex.

\smallskip\noindent\textbf{Step 2. Replacing auxiliary arcs by trails.}
The auxiliary graph gives the desired decompositions formally; its
arcs need not be edges of $G$. We replace each auxiliary arc by a
trail of bounded length in $G$, using the same trail in every
member of the auxiliary family containing that arc. There are two points to verify: enough trails
are available, and they can be chosen without accumulating large
degree at any vertex.

Every auxiliary arc has at least one end in $\mathcal C$, since
it touches a new vertex, now mapped to a color. For an end in $H$,
retain its formal class $A$ or $B$ before identification. We use
the following trails, all following the cyclic order
$A\to B\to\mathcal C\to A$:
\[
\begin{array}{c@{\qquad}l@{\qquad}c}
\text{Prescribed ends} & \text{Trail} & \text{Length}\\
c,c'\in\mathcal C & c\to x\to y\to c' & 3\\
u\in A,\ c'\in\mathcal C & u\to v\to c\to x\to y\to c' & 5\\
u\in B,\ c'\in\mathcal C & u\to c\to x\to y\to c' & 4\\
c\in\mathcal C,\ v\in A & c\to x\to y\to c'\to v & 4\\
c\in\mathcal C,\ v\in B & c\to x\to y\to c'\to u\to v & 5
\end{array}
\]
Only the two ends in the first column are prescribed. All other
vertices in a row are chosen subject to the adjacency and
distinctness requirements below. In all rows $x\in A$, $y\in B$,
and $xy$ is an original edge; $c$ must belong to $\mathcal P(x)$ and
$c'$ to $\mathcal P(y)$. In the ordinary case, the displayed vertices of
$A,B$ are interpreted in $H$ after identifying their formal copies.

Choose the extension at an end in $H$ first. In the second row,
choose $v\in N_H(u)$ and then $c\in \mathcal P(v)$; in the third,
choose $c\in \mathcal P(u)$. In the fourth row choose $c'\in \mathcal P(v)$;
in the fifth choose $u\in N_H(v)$ and then $c'\in \mathcal P(u)$.
Each such choice has $\Omega(m)$ possibilities, since both
the minimum degree and every palette size are at least $bm$.
With $c,c'$ fixed, choose $x$ with $c\in \mathcal P(x)$ and then a neighbor
$y$ with $c'\in \mathcal P(y)$. We say that a vertex \emph{admits} a color
when that color belongs to its palette. All but $O(\eta m)$ vertices of $A$ admit
$c$, and each has at least $bm-O(\eta m)$ neighbors admitting
$c'$. Since $|A|=\Theta(m)$, there are $\Omega(m^2)$ choices
for this central edge.

Choose each free color different from the prescribed colors,
and choose the vertices of $H$ on the trail to be distinct.
These exclusions remove only $O(m^{j-2})$ assignments from a
trail of length $j$. Thus each row has $\Omega(m^{j-1})$
realizations, with $3\le j\le5$. The prescribed colors in the
first row may coincide: when $c=c'$, the trail is the closed trail
$cxyc$ of length three. Otherwise each chosen trail is a path. The further restrictions
needed for different trails and listed pieces are imposed below.

We control the degrees in this step and Step~3 by keeping a
count at each vertex that bounds the number of incidences already
used or reserved there. We call each unit of this count a charge.
Fix a constant larger than the number of incidences in any one
local construction; it may depend on $L$. Whenever we choose a free
vertex, reserve this many units of degree at that vertex. For a
vertex placed on a trail, reserve units also for every listed piece
using that trail and its later graph $R_f$. There are only boundedly
many such pieces per arc. In Step~3, charge the freely chosen common
color for its entire local construction at once, and charge each
new internal path vertex when chosen. At prescribed vertices, the
required units are bounded by a constant times the number of their
incidences in the listed pieces. Thus every later edge is either
already charged at its prescribed end or charged when its free end
is chosen; there is no uncounted later contribution.

We now bound these charges. Lemma~\ref{palette:lem:circulation}, with
$\beta=\alpha/2$, gives $O(m^{2-\alpha/8})$ charges in total
and $O(m^{1-\alpha/8})$ prescribed charges at each vertex.
At most two formal vertices are mapped to any one vertex, and
each arc is used only boundedly many times. Thus these estimates
remain valid after the mapping.
Stop choosing a vertex freely when its charge reaches
$m^{1-\alpha/16}$. The number of vertices excluded by this rule is at most the total
charge divided by the threshold, namely
$O(m^{2-\alpha/8}/m^{1-\alpha/16})=O(m^{1-\alpha/16})=o(m)$.
The prescribed charge at a vertex is smaller than the threshold
in order of magnitude. A free choice adds only boundedly many
charges, so its overshoot is bounded as well. Hence the total
charge, and therefore the final degree, is $o(m)$ at every vertex.

We next check that these degree restrictions leave a valid choice.
For a block of $k=O(1)$ free choices with $\Omega(m^k)$
admissible assignments, these excluded vertices rule out only
$o(m^k)$ assignments. The same holds for edges of $D$ or
edges already used, since they have maximum degree $o(m)$
and every new edge has a free end. More explicitly, an excluded value in one coordinate rules out
$O(m^{k-1})$ assignments. A forbidden adjacency with two free
ends rules out $o(m^2)$ choices for those ends and hence
$o(m^k)$ assignments. With one prescribed end, its forbidden
degree excludes only $o(m)$ choices for the free end.
Requiring the free vertices to be distinct from each other and
from the prescribed vertices removes $O(m^{k-1})$ further assignments.
We can therefore choose the auxiliary trails to be pairwise
edge-disjoint and to have no edge in $D$. For every listed directed
triangle or digon, require all freely chosen vertices of $H$ in
its replacement trails to be distinct
and to avoid its prescribed vertices of $H$. These are only
boundedly many restrictions at each choice, because an arc lies
in boundedly many listed triangles or digons. The only members
with two old vertices
are the triangles on edges of $D^*$. If both vertices lie in $H$, they
are the distinct ends of an original edge. Every other piece has at
most one prescribed vertex of $H$. Thus the required distinctions
are possible. They are imposed within each listed piece;
different pieces may share vertices.

Let $D_{\mathrm{aux}}$ be the union of the chosen trails in $G$;
$D_{\mathrm{aux}}^*$ remains the auxiliary directed multigraph.
In each listed directed triangle or digon, replace every auxiliary
arc by its chosen trail and any arc of $D^*$ by its underlying
edge of $D$. Traversing these edges in cyclic order gives a closed
trail in $G$. Its edge set defines a balanced graph
$J_f\subseteq D_{\mathrm{aux}}\cup D$. The distinctness requirements
ensure that every vertex of $H$ in $J_f$ has one original edge
and one color edge. Hence the original edges of $J_f$ form a matching.
A color vertex may occur several times on the
trail; each occurrence has one incoming and one outgoing color edge.
It consequently has even degree, and in the bipartite case its
degrees into the two sides agree. Each edge of
$D_{\mathrm{aux}}\cup D$ lies in boundedly many $J_f$, since the auxiliary
trails are edge-disjoint and an edge of $D$ has at most two
alternatives. A listed formal piece contains at most one edge of $D^*$,
so no $J_f$ uses both alternatives of one edge of $D$.
For every balanced $J\subseteq D$, choose a balanced lift $J^*$.
The selected directed triangles and digons decompose $D_{\mathrm{aux}}^*\cup J^*$.
Replacing their auxiliary arcs by the chosen trails uses each
auxiliary trail once and one alternative of each edge of $J$.
The resulting subgraphs therefore decompose $D_{\mathrm{aux}}\cup J$.

\smallskip\noindent\textbf{Step 3. Absorbing each balanced piece.}
The original edges of $J_f$ form a matching. We will color all of
them with one common color $c_0$. To make that color available,
we join the paired vertices at each color occurrence by an odd path
and interchange two colors along it. The odd length ensures that
the path edge at each endpoint changes from $c_0$ to the color
of the incident color edge in $J_f$.

Fix $f$. We construct $R_f$ so that both $R_f$ and $R_f\cup J_f$
decompose into allowed triangles. Choose a color $c_0$ belonging
to every palette $\mathcal P(v)$ for $v\in V(J_f)\cap V(H)$, distinct from
every color vertex in $J_f$, so that
all edges $vc_0$ avoid $D_{\mathrm{aux}}\cup D$ and the graphs $R_{f'}$
already constructed. Reserve these edges. There are at least a fixed positive multiple of $m$ choices
before the $o(m)$ exclusions, since $J_f$ has
bounded order and the palette omissions are sufficiently small.
At each occurrence of a color vertex $c$ on the closed trail
defining $J_f$, pair the vertex immediately before it with the vertex immediately
after it. Every vertex of $V(J_f)\cap V(H)$ belongs to exactly
one such pair, because it has one original edge and one color
edge in $J_f$. The paired vertices lie on opposite sides in
the bipartite case.

For each such pair of distinct vertices, choose an ordering $x,y$,
with $x\in A$ and $y\in B$ in the bipartite case. Choose an
$(x,y)$-path of odd length $L$ in $H$. Require its internal
vertices to admit both $c$ and $c_0$, to lie outside
$V(J_f)\cap V(H)$, and to avoid the internal vertices of the paths
already chosen for this $f$. Since the pairs partition
$V(J_f)\cap V(H)$, the chosen paths for this $f$ are pairwise
vertex-disjoint.
We also require the path edges and the internal color edges to avoid
$D_{\mathrm{aux}}\cup D$, the earlier graphs $R_{f'}$, and the edges already
reserved or used for this $f$. Their union has maximum degree $o(m)$
by the degree accounting in Step~2. For the internal color edges,
it suffices to avoid the neighbors of $c$ and $c_0$ in this union,
which excludes $o(m)$ vertices. The vertices that have reached the
charge threshold also form a set of size $o(m)$.
The palette restrictions exclude $O(\eta m)$ further vertices.
Lemma~\ref{pre:walk-stability} therefore removes only
$O_{b,L}(\eta m^{L-1})+o(m^{L-1})$ of the at least
$bm^{L-1}$ distinct $(x,y)$-walks of length $L$, including all
walks with a repeated vertex. For sufficiently small $\eta$ and
large $m$, there remain $\Omega(m^{L-1})$ distinct $(x,y)$-paths
of length $L$ satisfying all these requirements.
Color each path alternately with $c_0,c$, starting and ending
with $c_0$. Let $R_f$ consist of the path edges and the color
edges used by this coloring. This is a decomposition of $R_f$
into allowed triangles. For a second decomposition, interchange
$c_0$ and $c$ along every path and give each original edge of
$J_f$ color $c_0$. Those edges form a matching. At each path
endpoint $v$, the path now uses color $c$ and the edge of $J_f$
uses color $c_0$. These use exactly the color edge $vc$ of
$J_f$ and the reserved edge $vc_0$ of $R_f$. Each internal
vertex still uses one edge of each color. Thus the second
coloring decomposes $R_f\cup J_f$ into allowed triangles.
Figure~\ref{palette:fig:absorber-switch} shows the interchange on one path.

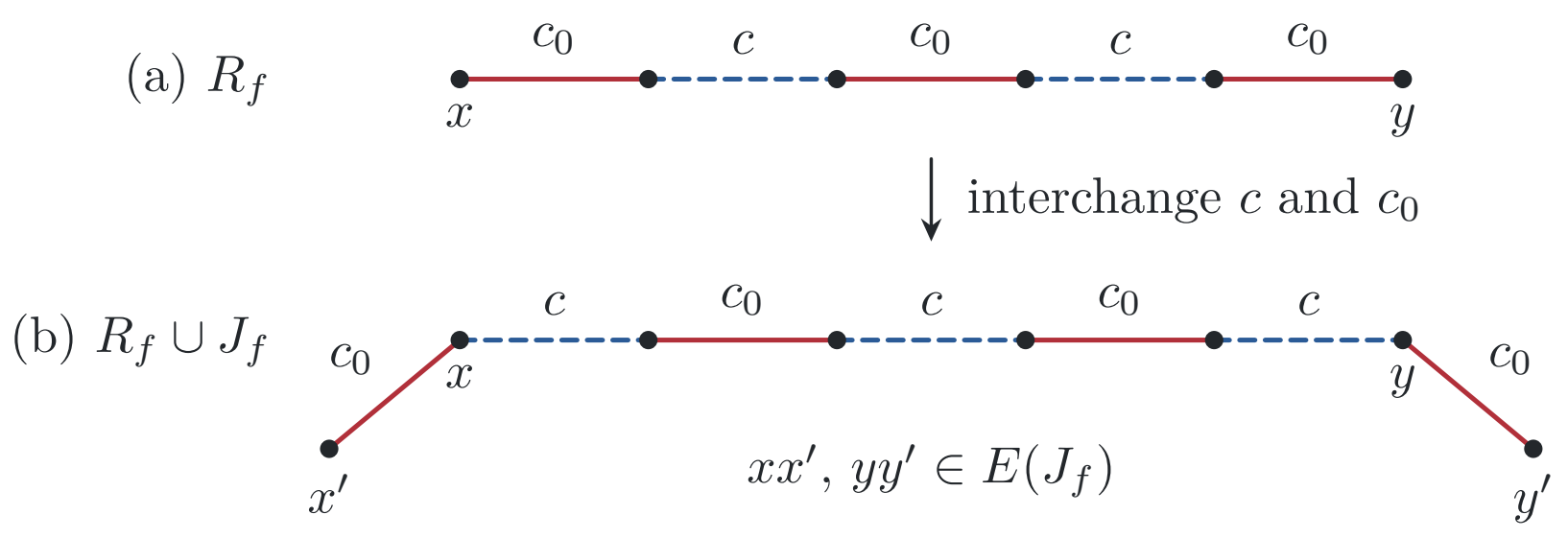
\begin{figure}[!htbp]
\centering
\begin{tikzpicture}[of diagram]
  \begin{scope}
    \foreach \i in {0,...,5} \coordinate (v\i) at (1.3*\i,0);
    \foreach \i/\j in {0/1,2/3,4/5}
      \draw[of red] (v\i)--(v\j) node[midway,above=3pt] {$c_0$};
    \foreach \i/\j in {1/2,3/4}
      \draw[of blue] (v\i)--(v\j) node[midway,above=3pt] {$c$};
    \foreach \i in {0,...,5} \node[of vertex] at (v\i) {};
    \node[below=3pt] at (v0) {$x$};
    \node[below=3pt] at (v5) {$y$};
    \node[of panel,anchor=east] at (-1.3,0) {(a) $R_f$};
  \end{scope}
  \draw[of arrow] (3.25,-.55)--(3.25,-1.12)
    node[midway,right=5pt] {interchange $c$ and $c_0$};
  \begin{scope}[yshift=-1.8cm]
    \foreach \i in {0,...,5} \coordinate (u\i) at (1.3*\i,0);
    \foreach \i/\j in {0/1,2/3,4/5}
      \draw[of blue] (u\i)--(u\j) node[midway,above=3pt] {$c$};
    \foreach \i/\j in {1/2,3/4}
      \draw[of red] (u\i)--(u\j) node[midway,above=3pt] {$c_0$};
    \coordinate (xp) at (-.9,-.75);
    \coordinate (yp) at (7.4,-.75);
    \draw[of red] (xp)--(u0) node[midway,above left=2pt] {$c_0$};
    \draw[of red] (u5)--(yp) node[midway,above right=2pt] {$c_0$};
    \foreach \i in {0,...,5} \node[of vertex] at (u\i) {};
    \node[of vertex] at (xp) {};
    \node[of vertex] at (yp) {};
    \node[below=3pt] at (u0) {$x$};
    \node[below=3pt] at (u5) {$y$};
    \node[below=3pt] at (xp) {$x'$};
    \node[below=3pt] at (yp) {$y'$};
    \node[of panel,anchor=east] at (-1.3,0) {(b) $R_f\cup J_f$};
    \node at (3.25,-.9) {$xx',\,yy'\in E(J_f)$};
  \end{scope}
\end{tikzpicture}
\caption{One odd path in the two colorings used by the absorber,
illustrated with length five. Only original edges are drawn.
The vertices $x,y$ are paired for
color $c$. After the interchange, each original edge of the
matching in $J_f$ receives color $c_0$. Here $xx'$ and $yy'$
are distinct matching edges; the proof also permits $xy$ to be
a matching edge.}
\label{palette:fig:absorber-switch}
\end{figure}

The common color $c_0$ is chosen freely and charged for its
entire graph $R_f$ at that choice. After the edges $vc_0$ are
reserved, every additional edge has a free internal end. Also,
for each color vertex, its prescribed uses are bounded by the
number of incidences at that vertex in the graphs $J_f$. Thus the degree argument applies
throughout and the union of all these graphs has maximum degree
$o(m)$.

Set $R=D_{\mathrm{aux}}\cup\bigcup_f R_f$. For a balanced $J\subseteq D$,
choose the subgraphs $J_f$ that decompose $D_{\mathrm{aux}}\cup J$.
Use the decomposition of $R_f\cup J_f$ at each chosen index
and the decomposition of $R_f$ at every other index. These
decompositions cover each edge of $R\cup J$ exactly once.
Taking $J=\emptyset$ shows that $R$ is balanced. The construction
also gives $R\subseteq G-E(D)$ and $\Delta(R)=o(m)$.
\end{proof}

\subsection{Completing the decomposition}
Our first task is to cover most edges by edge-disjoint allowed
triangles. A sparse graph reserved for the final step will finish
the covering provided few edges remain
\emph{at each vertex}, including each color vertex.

\begin{lemma}\label{palette:lem:greedy}
Fix positive constants $a,\lambda,\mu$. There exist a constant
$\kappa>0$ and a positive integer $m_0$, depending only on
$a,\lambda,\mu$, for which the following holds.
Let $m\ge m_0$, and let $D$ be a simple graph on at most $\mu m$
vertices with $N=e(D)\ge am^2$. Let $\mathcal T$ be a family
of distinct triangles of $D$, with a weight $\omega(C)$ assigned to each
$C\in\mathcal T$. Suppose that
\begin{enumerate}
\item $0\le \omega(C)\le \lambda/m$ for every $C\in\mathcal T$;
\item $\sum_{C\in\mathcal T:\,e\in E(C)}\omega(C)=1$ for every
$e\in E(D)$.
\end{enumerate}
Then some pairwise edge-disjoint triangles from $\mathcal T$
can be deleted from $D$ so that the remaining graph has maximum
degree $O_{a,\lambda,\mu}(m^{1-\kappa})$.
\end{lemma}
\begin{proof}
We use a random greedy procedure and show that it leaves few edges
at every vertex. Counting the total number of uncovered edges would
not suffice: the final absorption step needs a bound on the maximum
degree. We therefore track both the available triangle weight at
each edge and the remaining degree at each vertex.

\smallskip\noindent\textbf{The process and the quantities to track.}
A triangle is available if all three of its edges remain in the
graph. Choose an available triangle with probability proportional
to its fixed weight, and remove its three edges. After $i$ choices the fraction
of edges still present is $s=s_i=1-3i/N$. We track two quantities.
For every edge $e$ of $D$, let $\omega_e$ be the sum of the weights
of all triangles containing $e$ whose \emph{other two} edges remain,
even if $e$ itself has been removed. This convention avoids a sudden change in the tracked
quantity when $e$ is removed. If edges survived independently with
probability $s$, the other two edges of a triangle would survive
with probability $s^2$. This motivates the comparison with $s^2$;
the argument below proves the needed estimate without assuming
independence. For every vertex $v$, let $f_v$ be its remaining degree
divided by $\mu m$, and write $f_v^0$ for its initial value. We compare
$\omega_e$ with $s^2$ and $f_v$ with $s f_v^0$. Tracking $s f_v$
instead makes both target values constant multiples of $s^2$. Set
\[
 \varepsilon_i=\max\bigl\{\max_e|\omega_e-s^2|,\ \max_v|s f_v-s^2f_v^0|\bigr\}.
\]
Stop if $\varepsilon_i>s^2/2$, or when $i=i_{\mathrm{end}}=\lfloor N(1-m^{-\kappa})/3\rfloor$.
We will choose $\kappa$ small enough that, with positive probability,
the error never triggers a stop before $i_{\mathrm{end}}$. At time zero the error is
zero, since every edge has total triangle weight one.

\smallskip\noindent\textbf{Conditional changes in one step.}
At a step before stopping, the available weight $\omega_{\mathrm{tot}}$ is one third of the
sum of $\omega_e$ over the $Ns$ remaining edges, so $\omega_{\mathrm{tot}}\ge Ns^3/6>0$.
Thus a triangle of positive weight is available at every step
before stopping. In fact $\omega_{\mathrm{tot}}=Ns^3/3+O(Ns\varepsilon_i)$, because each of the $Ns$
remaining edges has available weight $s^2+O(\varepsilon_i)$.
Dividing its available weight by $\omega_{\mathrm{tot}}$, a remaining edge is removed with probability
\[
 \frac{\omega_e}{\omega_{\mathrm{tot}}}=\frac{3}{Ns}+O\!\left(\frac{\varepsilon_i}{Ns^3}\right).
\]
Two distinct edges are removed together with probability at most
$\lambda/(m\omega_{\mathrm{tot}})$, since they determine at most one triangle.
Write $\mathbb E_i$ for expectation given the current graph.
For any quantity $Z$ in the process, let $\operatorname{inc}_i(Z)$
denote its value after step $i+1$ minus its value after step $i$. A triangle counted by $\omega_e$ is lost if either of its
other two edges is removed. Summing these probabilities and subtracting
their intersection gives
\[
 \mathbb E_i\operatorname{inc}_i(\omega_e)
 =-\frac{6\omega_e}{Ns}
    +O\!\left(\frac{\omega_e \varepsilon_i}{Ns^3}+\frac{\lambda \omega_e}{mNs^3}\right).
\]
For the edges incident with a vertex $v$, sum the single-edge probabilities, obtaining
$\mathbb E_i\operatorname{inc}_i(f_v)=-3f_v/(Ns)+O(f_v \varepsilon_i/(Ns^3))$.
Since $\omega_e\le 3s^2/2$ and $f_v\le3s/2$, both estimates give the same
error bound. More precisely, let $Z$ denote $\omega_e$ or $s f_v$, with
$\gamma=1$ or $f_v^0$, respectively. Then
\[
 \left|\mathbb E_i\operatorname{inc}_i(Z)-\gamma\,\operatorname{inc}_i(s^2)\right|
 \le \lambda_0\left(\frac{\varepsilon_i}{Ns}+\frac{1}{mNs}+\frac{1}{N^2}\right)
\]
for a fixed $\lambda_0$. To obtain this bound for $\omega_e$, replace
$\omega_e$ in its leading term by $s^2+O(\varepsilon_i)$.
The leading change is then $-6s/N$, which agrees with
$\operatorname{inc}_i(s^2)=-6s/N+9/N^2$ up to $O(N^{-2})$.
The two error terms above are bounded by
$O(\varepsilon_i/(Ns)+1/(mNs))$ using $\omega_e\le3s^2/2$.
For the vertex quantity, use
$\operatorname{inc}_i(s f_v)=(s-3/N)\operatorname{inc}_i(f_v)-3 f_v/N$ and
$|f_v-s f_v^0|\le \varepsilon_i/s$. The term $1/N^2$ also accounts
for the quadratic term in $\operatorname{inc}_i(s^2)$. Substitution in this
identity gives leading term $-6f_v/N$; replacing $f_v$ by
$sf_v^0+O(\varepsilon_i/s)$ gives the same bound.

\smallskip\noindent\textbf{Accumulated random error.}
The conditional estimates describe the mean change at each step.
Subtracting that mean gives increments with conditional mean
zero, which we call centered increments. We bound their
partial sums simultaneously for every tracked edge and vertex;
this is what allows the estimates to hold throughout the process.

Each removal decreases $\omega_e$ by at most $3\lambda/m$:
at most three removed edges differ from $e$, and each belongs
to at most one triangle with $e$. It decreases $f_v$
by at most $2/(\mu m)$. We claim that each centered increment
$\xi=\operatorname{inc}_i(Z)-\mathbb E_i\operatorname{inc}_i(Z)$ satisfies
$|\xi|=O(1/m)$ and $\mathbb E_i\xi^2=O(1/(mN))$.
For the variance bound, the second moment of a nonnegative
decrease is at most its upper bound times its mean, and centering
cannot increase variance. This applies directly to the decrease in
$\omega_e$. For $Z=sf_v$, the only random term in its increment is
$(s-3/N)\operatorname{inc}_i(f_v)$; the other term is determined by the current
state. Its conditional variance is therefore at most that of
$\operatorname{inc}_i(f_v)$. The expected decrease in $\omega_e$ is $O(s/N)$ and
that in $f_v$ is $O(1/N)$, since their respective current values
are $O(s^2)$ and $O(s)$ and each relevant edge is removed with
probability $O(1/(Ns))$. These observations give the stated
$O(1/(mN))$ variance bound.

Take centered increments to be zero after stopping.
For $\theta=\sqrt m$, we have $|\theta \xi|\le1$ for large $m$.
The inequality $e^t\le1+t+t^2$ for $|t|\le1$, together with
$\mathbb E_i\xi=0$ and the variance bound, gives
$\mathbb E_i e^{\pm\theta \xi}\le\exp(O(1/N))$.
There are at most $N/3$ steps. Iteration and Markov's inequality bound the probability that any fixed centered sum
exceeds $m^{-1/4}$ in absolute value by $2\exp(O(1)-m^{1/4})$.
No independence between successive increments is needed: at each
step we condition on the current graph and then average the
conditional exponential bound.
The exponential bound is smaller than any fixed negative power of
$m$. There are $O(m^2)$ tracked quantities and $O(m^2)$ times.
A union bound therefore shows that, with probability $1-o(1)$,
all their centered sums satisfy this bound simultaneously.

\smallskip\noindent\textbf{Keeping the error below the stopping threshold.}
Fix an outcome for which all these centered sums are small. The
remaining error comes from the conditional means, and their bounds
depend on the errors accumulated so far. We now control this
feedback by an elementary product estimate.

On this event, sum the errors in the conditional means up to
any $i$, including the stopping step. The accumulated deterministic error other than the
terms containing $\varepsilon_j$ is
$O(m^{-1}\log m+N^{-1})=o(m^{-1/4})$. Since
$\sum_{j<i}1/(Ns_j)\le\frac13\log(1/s_i)$, we obtain
\[
 \varepsilon_i\le 2m^{-1/4}+\sum_{j<i}\frac{\lambda_0 \varepsilon_j}{Ns_j}
 \le 2m^{-1/4}\prod_{j<i}\left(1+\frac{\lambda_0}{Ns_j}\right)
 \le 2m^{-1/4}s_i^{-\lambda_0/3}.
\]
To see the product bound, set $a_0=2m^{-1/4}$ and
$\beta_j=\lambda_0/(Ns_j)$. The sequence
$a_0\prod_{j<i}(1+\beta_j)$ equals
$a_0+\sum_{j<i}\beta_j a_0\prod_{h<j}(1+\beta_h)$,
by telescoping.
Induction therefore bounds any sequence satisfying the first
inequality by this product. For the last inequality, use
$1+t\le e^t$ and the bound above on $\sum_{j<i}1/(Ns_j)$.
Choose $\kappa>0$ with
$\kappa(2+\lambda_0/3)<1/4$. Since $s_i\ge m^{-\kappa}$ for $i\le i_{\mathrm{end}}$, this choice gives
$\varepsilon_i/s_i^2\le2m^{-1/4+\kappa(2+\lambda_0/3)}=o(1)$.
In particular, a first stop caused by excessive error would
contradict the estimate at that very step. The process therefore
reaches $i_{\mathrm{end}}$. Each remaining degree is at most
$(3/2)\mu m s_{i_{\mathrm{end}}}=O(m^{1-\kappa})$, as required.
\end{proof}

\begin{proof}[Proof of Theorem~\ref{palette:thm:main}]
We will choose a sparse graph $D$ and an absorber $R$ for $D$.
It then suffices to cover every edge of $G-E(D\cup R)$ by
edge-disjoint allowed triangles having no edge in $R$. The
uncovered edges of $G-E(R)$ will form a balanced subgraph of
$D$, which can be absorbed by $R$.

The original graph has at least $m$ vertices and minimum degree at least
$bm$, and the augmented graph has at most $4m$ vertices. Use
$a=b/4$, $\lambda=2/b$ and $\mu=4$ in Lemma~\ref{palette:lem:greedy}.
Fix $0<\alpha<\min\{\kappa/8,1/4\}$. Apply
Lemma~\ref{palette:lem:absorber} with this $\alpha$ and the fixed degree
constant from Step~1, and let $\eta_0$ be its error bound.
Choose $0<\eta\le\min\{\eta_0,b/8,1/2\}$ sufficiently small
in terms of these constants and $b,L$; finally, choose $m$
sufficiently large. All further restrictions on $\eta$ below
depend only on these fixed constants.

\smallskip\noindent\textbf{Step 1. The reserve and absorber.}
Set $p=m^{-\alpha}$, and include each edge of $G$ independently
with probability $p$ in a graph $D_0$. Every edge of $G$ belongs
to at least $(b/2)m$ allowed triangles, after decreasing $\eta$ if necessary. For an original edge, choose
a common available color at its ends. For a color edge $vc$,
choose a neighbor of $v$ in $H$ that admits $c$. For each fixed
edge $e$, the two-edge subgraphs obtained by deleting $e$ from
these triangles are pairwise edge-disjoint.
Also, $\Delta(G)\le2m$: each original vertex has degree
$2d_H(v)\le2m$, and each color has at most $|V(H)|\le2m$
neighbors. Lemma~\ref{bip:pre:chernoff} and a union bound show that, with
probability $1-o(1)$,
\[
 \Delta(D_0)\le4pm,\qquad
 \text{every edge has }\Omega(p^2m)\text{ allowed triangles
 whose other two edges lie in }D_0.
\]
Complete every edge of $D_0$ to an allowed triangle using two unused
edges outside $D_0$. To see that this can be done greedily, recall that
$\Delta(D_0)\le4pm$ and every edge has at least $b_1m$ possible third
vertices for a constant $b_1>0$ depending only on the walk parameters.
Indeed $\delta(H)\ge bm$, $\Delta\ge bm$, and the palette omissions are
at most $\eta m$ per vertex and per color on each relevant side.
Stop using a vertex as the third vertex once it has been chosen
$\lceil \lambda_1 pm\rceil$ times. There are $O(pm^2)$ edges in $D_0$, so at
most $O(m/\lambda_1)$ vertices are excluded in this way. Choose the fixed
constant $\lambda_1$ large enough that this is less than $b_1m/2$. At either
prescribed endpoint the used degree is at most
$2\Delta(D_0)+2\lceil \lambda_1 pm\rceil=o(m)$, and the edges of $D_0$ exclude
at most $4pm$ further neighbors there. A third vertex is available.

Let $D$ be the union of the chosen triangles. This graph is balanced,
has maximum degree at most $\lambda_{\mathrm{res}}pm$ for a fixed
$\lambda_{\mathrm{res}}>0$, and contains $D_0$.
We call $D$ the \emph{reserve}, since its edges will be used to cover
the remainder left by the random greedy procedure.
In particular, every edge still has $\Omega(p^2m)$ allowed triangles
whose two other edges lie in $D$. Apply Lemma~\ref{palette:lem:absorber}
with $\lambda=\lambda_{\mathrm{res}}$. This constant can be fixed before $\eta$ by using $b_1=b/2$
when $\eta\le b/4$. Let $R$ be the resulting balanced absorber for $D$.
Then $G_0=G-E(D\cup R)$ is balanced. Set $H_0=G_0[V(H)]$ and
write $\mathcal P_0(v)=\{c:vc\in E(G_0)\}$ for the remaining palette at $v$.
Balance gives $|\mathcal P_0(v)|=d_{H_0}(v)$. Since
$\Delta(D\cup R)=o(m)$, we have $\Delta-|\mathcal P_0(v)|\le2\eta m$,
and each color is absent from at most $2\eta m$ of the remaining
palettes on each relevant side.

\smallskip\noindent\textbf{Step 2. Fractional triangle weights.}
To apply the random greedy lemma, we need a fractional triangle
decomposition: nonnegative weights on allowed triangles such that
the weights of all triangles containing any fixed edge sum to one.
For each original edge, we can obtain this sum by distributing weight
one among its allowed triangles. The difficulty is to obtain sum
one at the color edges as well. We first
correct those sums while temporarily allowing all colors, then
remove the forbidden adjustments. The latter operation introduces
a smaller error, so the procedure can be repeated.

Initially distribute weight one on each remaining original edge $xy$
equally among the colors in $\mathcal P_0(x)\cap \mathcal P_0(y)$. There are
$\Delta-O(\eta m)$ such colors, so each weight is
$(1+O(\eta))/\Delta$.
At each remaining color edge $vc$, define its \emph{residual}
$\rho_{vc}$ to be one minus the sum of the weights of its triangles, and set $\rho_{vc}=0$ when $vc$ is absent. Initially
$\rho=\max_{v,c}|\rho_{vc}|=O(\eta)$. Indeed, for an available
$vc$, there are $\Delta-O(\eta m)$ neighbors of $v$ in $H_0$
that admit $c$, and each corresponding triangle has weight
$(1+O(\eta))/\Delta$. We describe one
correction that preserves total one on each original edge and decreases
$\rho$ by a factor of at least two. A positive residual requires additional weight, while a negative
residual requires a decrease.

\smallskip\noindent\emph{A correction with no palette restrictions.}

We do the calculation on a bipartite graph with sides $A,B$.
In the ordinary case use the bipartite double cover of the remaining
original graph. Give both copies of a vertex its current palette and
each lift of an edge its current weights. For $u,v\in V(H_0)$,
projection gives a bijection from the length-$L$ $(u^+,v^-)$-walks
in the double cover to the length-$L$ $(u,v)$-walks in $H_0$,
including when $u=v$.
Thus, in either case, for every $(x,y)\in A\times B$ there are
$\Omega(m^{L-1})$ distinct $(x,y)$-walks of length $L$ in the
bipartite graph under consideration, by Lemma~\ref{pre:walk-stability}. In the ordinary case, average the
adjustments on the two lifts of each edge after the correction below.
This preserves the total weight on that edge and makes the residual
at a color edge the average of the residuals at its two lifts.
Since the correction makes both lifted residuals zero, the
averaged adjustment corrects the original residual as well.

First allow every color on every original edge for the purpose of
calculating a signed adjustment. The residuals satisfy
\[
 \sum_{c\in\mathcal C}\rho_{vc}=0\quad(v\in A\cup B),
 \qquad
 \sum_{x\in A}\rho_{xc}=\sum_{y\in B}\rho_{yc}
 \quad(c\in\mathcal C).
\]
For the first identity, the number of color edges at a vertex equals
its original degree, and each incident original edge has total
weight one. For the second, each color has equal degrees into the
two sides and each triangle contributes the same weight on both.
In the double cover the residuals are identical on the two copies,
so the second identity also holds there.

These identities are precisely what the correction needs. The
zero sum over colors lets us change the color weights without
changing their sum on an original edge. Equality of the two side
sums lets us correct one side and then the other without undoing
the first correction. No equality of $|A|$ and $|B|$ is needed.

For $x\in A$ and each color $c$, distribute $\rho_{xc}$ equally among
the vertices $y\in B$, then distribute each share equally among the
distinct $(x,y)$-walks of length $L$. For each traversal of an edge
on such a walk, add the signed share to its color-$c$ weight with
alternating signs, starting positive. If an edge is traversed more
than once, add all its contributions. Consecutive contributions
cancel at every internal position, and the two ends both gain
the signed share, since $L$ is odd
(Figure~\ref{palette:fig:walk-cancellation}).

\begin{figure}[!htbp]
\centering
\begin{tikzpicture}[of diagram]
  \foreach \i in {0,...,5} \coordinate (v\i) at (1.3*\i,0);
  \foreach \i/\j in {0/1,2/3,4/5}
    \draw[of red] (v\i)--(v\j) node[midway,above=3pt] {$+1$};
  \foreach \i/\j in {1/2,3/4}
    \draw[of blue] (v\i)--(v\j) node[midway,above=3pt] {$-1$};
  \foreach \i in {0,...,5} \node[of vertex] at (v\i) {};
  \node[left=4pt] at (v0) {$x$};
  \node[right=4pt] at (v5) {$y$};
  \node[anchor=east] at (-.5,-.55) {sum};
  \foreach \i/\s in {0/1,1/0,2/0,3/0,4/0,5/1}
    \node at (1.3*\i,-.55) {$\s$};
\end{tikzpicture}
\caption{Alternating signs on an odd walk cancel at every
internal position and give sum one at each end. The example has length five. The same cancellation holds when vertices or edges repeat.}
\label{palette:fig:walk-cancellation}
\end{figure}
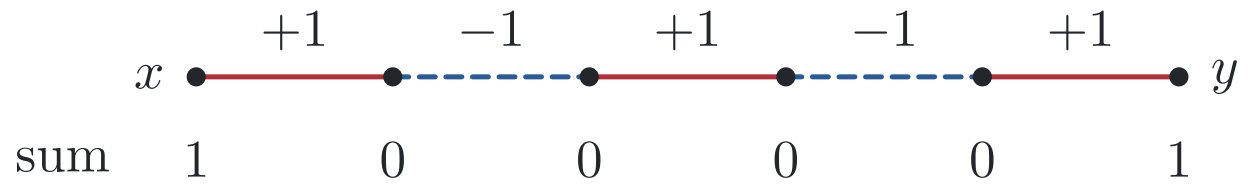

The total color-$c$ adjustment at each $x\in A$ is now $\rho_{xc}$.
Every $y\in B$ receives the same amount
$\mu_c=\sum_{x\in A}\rho_{xc}/|B|$. Apply the same distribution rule
with $A$ and $B$ interchanged, using the remaining amounts
$\rho_{yc}-\mu_c$. Their sum on $B$ is zero, so this
second operation supplies the adjustment at each $y$ and contributes
zero at every $x$. These cancellations count incidences with
multiplicity, so repeated vertices or edges in a walk cause no change.

Write $h_{ec}$ for the resulting adjustment. Its sum over colors
is zero on each edge $e$: the same walk distributions are used
for every color and the amounts sent from each vertex sum to zero.
For the second operation this also uses $\sum_c\mu_c=0$,
which follows by summing the original residuals over $A$ and all colors.
Both sides have order $\Theta(m)$, so
$|\mu_c|\le |A|\rho/|B|=O(\rho)$.
Each signed share in either operation therefore has absolute value
$O(\rho/m^L)$.
Fixing the position and direction of a traversal of an edge $e$
leaves at most $L-1$ free vertex positions. Thus the total number
of occurrences of $e$ among all walks of length $L$ is
$O_L(m^{L-1})$. Summing the corresponding contributions gives
$|h_{ec}|=O(\rho/m)$.

\smallskip\noindent\emph{Restoring the palette restrictions.}
At this point the signed adjustments correct every residual
exactly, but some adjustments may use forbidden colors. Erasing
them loses exact correction at the color edges. We show that the
resulting error is much smaller than the error we started with.

For $e=xy$, let $a_e$ be the sum of the adjustments $h_{ec}$ over
colors outside $\mathcal P_0(x)\cap\mathcal P_0(y)$. Erase these
adjustments by setting them to zero. Since the original sum
$\sum_c h_{ec}$ was zero, the retained adjustments sum to $-a_e$.
Add $a_e$ divided by the number of retained colors to each
retained adjustment. Their sum on $e$ is then zero again. Only $O(\eta m)$
colors are erased and $\Theta(m)$ remain. Since each erased
adjustment has size $O(\rho/m)$, their sum has size $O(\eta\rho)$;
sharing this sum among the surviving colors changes each retained
adjustment by $O(\eta\rho/m)$. For a color edge $vc\in E(G_0)$,
summing these changes over the $O(m)$ incident original edges
costs $O(\eta\rho)$. At most $O(\eta m)$ neighbors lack $c$;
erasing their terms costs a further $O(\eta\rho)$.
The new residual is therefore $O(\eta\rho)$, at most $\rho/2$
for sufficiently small $\eta$. Each triangle weight changes by
$O(\rho/m)$ in this round. The
correction preserves total weight one on every original edge and
uses only allowed triangles, so the balance identities remain valid
for the next round.

\smallskip\noindent\emph{Iteration and nonnegative limiting weights.}
After $j$ rounds the largest residual is at most $2^{-j}$ times
its initial value. The changes in any one triangle weight are
bounded by a constant times this residual divided by $m$, so their
sum converges. The limiting residual is zero at every color edge,
and the original-edge sums have remained one throughout.
Thus the limiting weights total one on every edge of $G_0$.
Summing the geometric series, their total change is
$O(\eta/m)$; since the initial weights are
$(1+O(\eta))/\Delta$ and $\Delta=\Theta(m)$, sufficiently small
$\eta$ keeps the limiting weights in $[0,2/\Delta]$.

\smallskip\noindent\textbf{Step 3. Covering the remaining edges.}
The weights of Step~2 total one at every edge of $G_0$ and are at most
$2/\Delta\le 2/(bm)$. Removing $D\cup R$ cost $o(m)$ edges per vertex,
so $G_0$ still contains at least $bm^2/4$ original edges.
Lemma~\ref{palette:lem:greedy} gives edge-disjoint allowed triangles
whose uncovered graph
$J$ in $G_0$ has maximum degree $d=O(m^{1-\kappa})=o(p^4m)$.
If $J$ has no edges, the uncovered edges of $G-E(R)$ already lie
in $D$. Otherwise, cover the edges of $J$ in any order. For each one, Step~1
gives at least $b_2 p^2m$ choices of a third vertex completing an allowed
triangle with two edges of $D$, for a fixed $b_2>0$. Use a vertex as
the third vertex at most $q_{\mathrm{cap}}=\lceil \lambda_2 d/p^2\rceil$ times. Since there
are $O(md)$ target edges, fewer than $b_2p^2m/3$ vertices reach their
cap if $\lambda_2$ is large enough. At each endpoint at most $d+2q_{\mathrm{cap}}$ reserve
edges have been used: one for each earlier target incident there, and
two for each use as a third vertex. The two endpoints together thus exclude at most
$2(d+2q_{\mathrm{cap}})=o(p^2m)$ more choices, since $d=o(p^4m)$. A third vertex
remains available at each step, giving the required edge-disjoint
triangles.

Let $J_D$ be the graph formed by the edges of $G-E(R)$ that
remain uncovered. Then $J_D\subseteq D$, and $J_D$ is balanced
because $G$, $R$ and every chosen triangle are balanced.
By Lemma~\ref{palette:lem:absorber}, the graph $R\cup J_D$
has a decomposition into allowed triangles. Together with the
triangles already chosen, these form a decomposition of $G$, which gives
the required edge coloring of $H$.
\end{proof}

\section{Sharpness constructions}\label{sharp:sec}
We give the sharpness constructions for Theorems~\ref{main:quarter}
and~\ref{intro:thm:regular}. Both use the same Petersen multigraph.

\subsection{The general bound}\label{sharp:sec:general}
The following construction shows that the degree bound in
Theorem~\ref{main:quarter} is best possible. The odd-degree
construction uses four cliques.
Replacing these cliques by slightly larger graphs gives the
even-degree construction of Folz~\cite{Folz2026}.

\begin{samepage}
\begin{proposition}\label{sharp:thm:family}
For every integer $\Delta\ge3$, there is a connected simple
graph $G$ of maximum degree $\Delta$ with the following properties.
\begin{enumerate}
\item $\chi'(G)=\Delta+1$ and $G$ has no $\Delta$-overfull set.
\item Its order is $4\Delta+9$ when $\Delta$ is even and
$4\Delta+5$ when $\Delta$ is odd.
\item $e(G)=\Delta(|V(G)|-1)/2$, and $G$ has a vertex of degree two.
\end{enumerate}
\end{proposition}
\end{samepage}
\begin{proof}
Set $q=2\lfloor\Delta/2\rfloor+2$. If $\Delta$ is odd, let
$D=K_{\Delta+1}$. If $\Delta$ is even, let $D$ be obtained from
$K_{\Delta+2}$ by deleting a perfect matching. In either case,
$D$ has order $q$, is $\Delta$-regular, and has a
$1$-factorization: take a $1$-factorization of the complete
graph of even order $q$, deleting one factor when $\Delta$ is even.
Fix this factorization of $D$.

Let $D_1,\ldots,D_4$ be pairwise vertex-disjoint copies of $D$.
For each $i\in[4]$, choose $z_i\in V(D_i)$ and set
$H_i=D_i-z_i$. Choose distinct $u_i,v_i\in N_{D_i}(z_i)$ and
let $W_i=N_{D_i}(z_i)\setminus\{u_i,v_i\}$. Thus
$|V(H_i)|=q-1$ and $|W_i|=\Delta-2$.
Add five vertices $w_0,\ldots,w_4$ and, in addition to the edges
inside the $H_i$, include precisely
\[
 \{w_iw_{i+2}:i\in\mathbb Z_5\},\qquad
 \bigcup_{i=1}^4\{w_ix:x\in W_i\},\qquad
 \{v_iu_{i+1}:i\in[3]\}.
\]
Here $\mathbb Z_5=\{0,1,2,3,4\}$ with addition modulo five;
the subscripts on the $w$-vertices are read in this way. The
resulting graph $G$ is simple. Each $H_i$ is connected and
$W_i\ne\emptyset$, so the edges from $w_i$ to $W_i$, together
with the cycle on the $w$-vertices, make $G$ connected. Its order is
$4(q-1)+5=8\lfloor\Delta/2\rfloor+9$.
Every vertex has degree $\Delta$, except that
$d_G(u_1)=d_G(v_4)=\Delta-1$ and $d_G(w_0)=2$.
Consequently, its total degree deficiency is $\Delta$ and
$e(G)=\Delta(|V(G)|-1)/2$.
For odd $\Delta$, each $H_i$ is $K_\Delta$; see
Figure~\ref{sharp:fig:construction}.

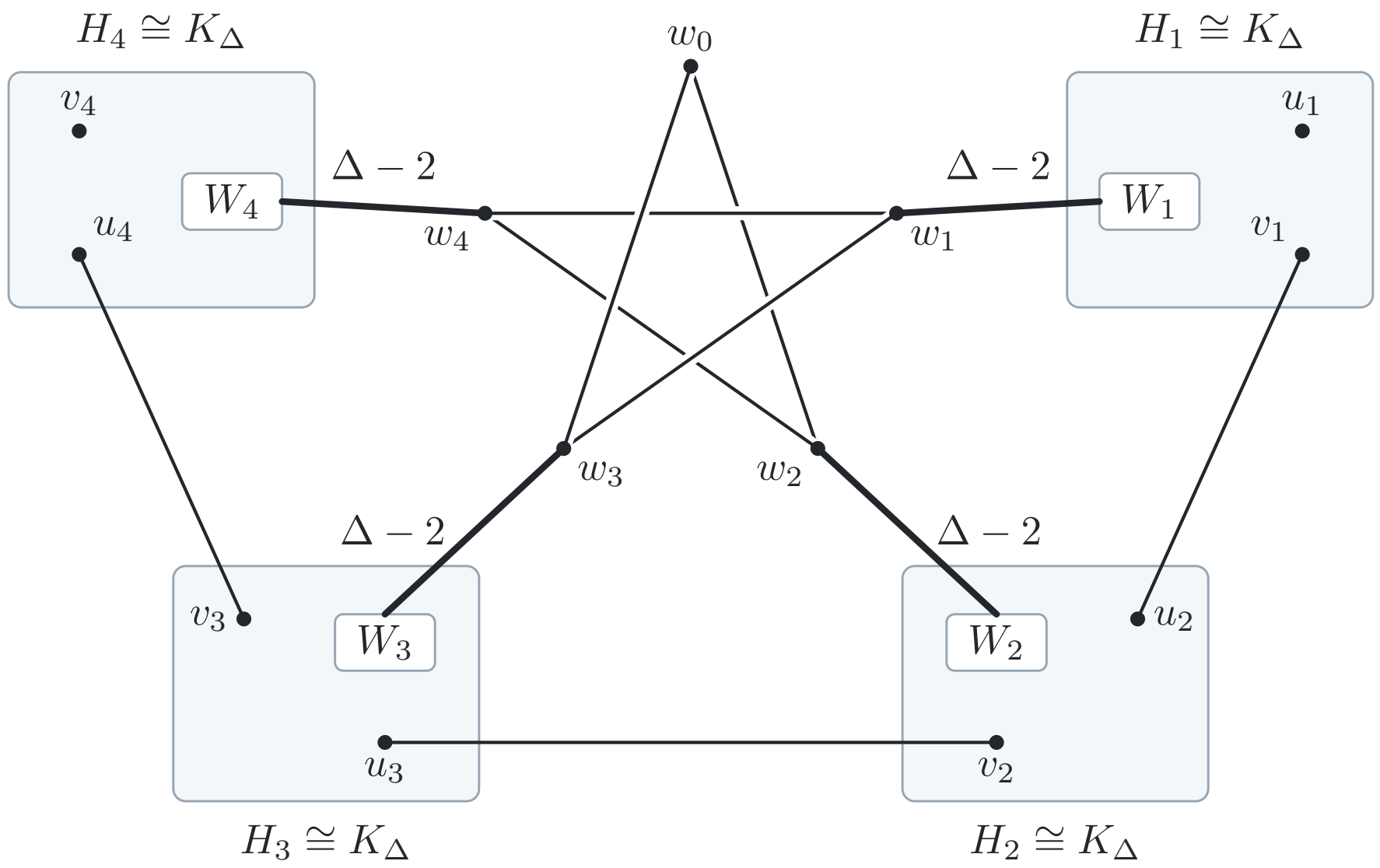
\begin{figure}[htbp]
\centering
\begin{tikzpicture}[of diagram]
  % Equal-sized clique boxes, with the same vertex and edge styles as above.
  \draw[of region] (3.2,0) rectangle (5.8,2);
  \draw[of region] (1.8,-4.2) rectangle (4.4,-2.2);
  \draw[of region] (-4.4,-4.2) rectangle (-1.8,-2.2);
  \draw[of region] (-5.8,0) rectangle (-3.2,2);
  \node at (4.5,2.35) {$H_1\cong K_{\Delta}$};
  \node at (3.1,-4.55) {$H_2\cong K_{\Delta}$};
  \node at (-3.1,-4.55) {$H_3\cong K_{\Delta}$};
  \node at (-4.5,2.35) {$H_4\cong K_{\Delta}$};
  \coordinate (u1) at (5.2,1.5);
  \coordinate (v1) at (5.2,.45);
  \coordinate (u2) at (3.8,-2.65);
  \coordinate (v2) at (2.6,-3.7);
  \coordinate (u3) at (-2.6,-3.7);
  \coordinate (v3) at (-3.8,-2.65);
  \coordinate (u4) at (-5.2,.45);
  \coordinate (v4) at (-5.2,1.5);
  \node[of set] (W1) at (3.9,.9) {$W_1$};
  \node[of set] (W2) at (2.6,-2.85) {$W_2$};
  \node[of set] (W3) at (-2.6,-2.85) {$W_3$};
  \node[of set] (W4) at (-3.9,.9) {$W_4$};
  \draw[of edge] (v1)--(u2) (v2)--(u3) (v3)--(u4);
  % The five central vertices retain the inner cycle of the Petersen graph.
  \coordinate (c0) at (0,2.05);
  \coordinate (c1) at (1.75,.8);
  \coordinate (c2) at (1.08,-1.2);
  \coordinate (c3) at (-1.08,-1.2);
  \coordinate (c4) at (-1.75,.8);
  \foreach \i/\j in {0/2,2/4,4/1,1/3,3/0}
    \draw[of edge,of gap] (c\i)--(c\j);
  \draw[of bundle] (c1)--node[above=4pt] {$\Delta-2$} (W1.west);
  \draw[of bundle] (c2)--node[right=5pt] {$\Delta-2$} (W2.north);
  \draw[of bundle] (c3)--node[left=5pt] {$\Delta-2$} (W3.north);
  \draw[of bundle] (c4)--node[above=4pt] {$\Delta-2$} (W4.east);
  \node[of vertex,label=above:$u_1$] at (u1) {};
  \node[of vertex,label=above left:$v_1$] at (v1) {};
  \node[of vertex,label=right:$u_2$] at (u2) {};
  \node[of vertex,label=below:$v_2$] at (v2) {};
  \node[of vertex,label=below:$u_3$] at (u3) {};
  \node[of vertex,label=left:$v_3$] at (v3) {};
  \node[of vertex,label=above right:$u_4$] at (u4) {};
  \node[of vertex,label=above:$v_4$] at (v4) {};
  \node[of vertex,label=above:$w_0$] at (c0) {};
  \node[of vertex,label=below right:$w_1$] at (c1) {};
  \node[of vertex,label=below left:$w_2$] at (c2) {};
  \node[of vertex,label=below right:$w_3$] at (c3) {};
  \node[of vertex,label=below left:$w_4$] at (c4) {};
\end{tikzpicture}
\caption{The construction when $\Delta$ is odd. Each shaded box represents a clique
$H_i\cong K_{\Delta}$; its internal edges are omitted. The set
$W_i=V(H_i)\setminus\{u_i,v_i\}$ has $\Delta-2$ vertices, and each
thick segment represents the $\Delta-2$ edges joining $w_i$ to $W_i$.
The thin edges between boxes are $v_1u_2$, $v_2u_3$, and $v_3u_4$.
The five central vertices form the cycle $w_0w_2w_4w_1w_3w_0$.
Gaps at crossings distinguish edges and do not indicate missing edges.}
\label{sharp:fig:construction}
\end{figure}

We show that $G$ is of class~2. Write the Petersen graph $R$ with
vertex set $\{x_i,w_i:i\in\mathbb Z_5\}$ and edges
$x_ix_{i+1}$, $w_iw_{i+2}$, and $x_iw_i$ for
$i\in\mathbb Z_5$. Let $R'$ be obtained from $R$ by replacing
each edge $x_iw_i$ by $\Delta-2$ parallel copies. Thus $R'$ is
$\Delta$-regular. We call the edges $x_iw_i$ and their parallel
copies \emph{spokes}; all other edges are \emph{nonspoke edges}.

Every perfect matching of $R$ contains either one spoke or five
spokes. Indeed, if it contains $t$ spokes, the remaining $5-t$
vertices on each of the two five-cycles must be matched internally,
so $t$ is odd. If $t=3$, the two remaining indices would have to
differ both by $\pm1$ and by $\pm2$ modulo five, which is
impossible. Hence every perfect matching of $R'$ contains either
four or zero nonspoke edges. A $\Delta$-edge-coloring of $R'$
would partition its ten nonspoke edges into sets of size four,
a contradiction.

Replace $x_i$ in $R'$, for each $i\in[4]$, by $H_i$, attaching
its $\Delta$ incident edge copies bijectively to
$N_{D_i}(z_i)$. Attach $x_{i-1}x_i$ to $u_i$, attach
$x_ix_{i+1}$ to $v_i$, and attach the spoke copies to $W_i$.
Denote the resulting
$\Delta$-regular multigraph by $R^*$. With these attachments,
$G=R^*-x_0$.

To obtain a contradiction, we show that a $\Delta$-edge-coloring
of $G$ would extend to $R^*$ and then contract to one of $R'$.
Suppose that $G$ has such a coloring. Its order is
odd and $e(G)=\Delta(|V(G)|-1)/2$, so every color class has
exactly $(|V(G)|-1)/2$ edges and misses exactly one vertex.
At any vertex $v$, the number of missing colors is
$\Delta-d_G(v)$, which is the number of deleted edge copies
joining $v$ to $x_0$ in $R^*$. Restore these copies using the
colors missing at $v$. Since each color is missing at only one
vertex, the restored edges have distinct colors at $x_0$.
This gives a $\Delta$-edge-coloring of $R^*$.

Every color class in this coloring of $R^*$ is a perfect
matching. Since $|V(H_i)|=q-1$ is odd, each color occurs on
at least one edge of $\partial_{R^*}(V(H_i))$. This cut has exactly
$\Delta$ edges, so each color occurs on exactly one of them.
Deleting the edges inside each $H_i$ and identifying its vertices
to $x_i$ therefore gives a $\Delta$-edge-coloring of $R'$, a contradiction.
Thus $\chi'(G)=\Delta+1$ by Vizing's theorem.

We verify the odd-set bound using the same perfect matchings of the
Petersen graph. For each $i\in\mathbb Z_5$ there is a unique perfect
matching whose only spoke is $x_iw_i$: deleting index $i$ leaves a
four-vertex path in each cycle. Use those five matchings and take
$2\Delta-5$ further copies of the all-spoke matching. These form
a family of $2\Delta$ indexed matchings in which repetitions are allowed.
They cover every nonspoke edge
twice: the five single-spoke matchings contribute ten incidences on
each five-cycle, equally among its five edges by rotation. Each spoke
is covered $1+(2\Delta-5)=2(\Delta-2)$ times. Assign these
occurrences, two apiece, to the $\Delta-2$ parallel copies in $R'$.

Each matching of $R'$ uses one edge at $x_i$. In $D_i$, take the unique
factor of the fixed $1$-factorization containing the edge at $z_i$
corresponding to that incidence. Delete $z_i$ from this factor and use
the chosen joining edge at its uncovered neighbor. Doing this at each
of $x_1,\ldots,x_4$ extends every indexed matching to a perfect
matching of $R^*$. Each factor of $D_i$ is used twice, because its corresponding edge
at $x_i$ is used twice. Thus every internal edge and every joining
edge of $R^*$ is covered twice.
Restricting to $G$ gives $2\Delta$ matchings covering each edge twice.
An odd set $X$ contains at most $(|X|-1)/2$ edges of any one matching.
Summing over the matchings gives $2e_G(X)\le\Delta(|X|-1)$.
\end{proof}

The construction identifies the obstruction excluded by
\eqref{boundary:eq:small-even-order}. Attach private leaves to
make every original vertex have degree $\Delta$, and set
$S_i=V(H_i)\cup\{w_i\}$ for $i\in[4]$. Let a fifth part $W$ consist of $w_0$ and all
private leaves. Each $S_i$ has even order
$q$, and $|V(G)|=4q+1$. Here $q=\Delta+1$ when $\Delta$ is odd
and $q=\Delta+2$ when $\Delta$ is even. Thus all four even
parts exceed $\Delta$, at the first excluded order in each case.

Delete edges within these five parts and contract each part,
labeling $W$ by $0$ and $S_i$ by $i$. The quotient is $K_5$:
its edges are the two complementary cycles $0\,1\,2\,3\,4\,0$ and
$0\,2\,4\,1\,3\,0$. At each label the two incidences of the second cycle
come from the same vertex $w_i$ before contraction, whereas those of the
first cycle belong to distinct vertices. This is an exceptional
forbidden system. More explicitly, set $Y_i=\{u_i,v_i,w_i\}$.
In the graph with private leaves, the cut $\partial(S_i)$ has
four edges, all incident with $Y_i$. Hence
\[
 |\partial(S_i)|-2d_{\partial(S_i)}(Y_i)=-4.
\]
Thus the positive inequality in
\eqref{boundary:eq:small-even-part} fails. All boundary edges
belong to this exceptional system, so no edge-disjoint
compatible circuit remains to absorb it.

The degree bound is also the exact threshold for each
sufficiently large order.

\begin{corollary}\label{sharp:cor:cutoff}
For every integer $n\ge n_0$, set
$\Delta=2\lfloor(n-1)/8\rfloor-1$.
There is a connected graph of order $n$ and maximum degree
$\Delta$ that is of class~2 and has no $\Delta$-overfull set.
\end{corollary}
\begin{proof}
The integer $\Delta$ is odd and $0\le n-(4\Delta+5)\le7$.
By Proposition~\ref{sharp:thm:family}, there is a connected
counterexample of order $4\Delta+5$ with a vertex of degree two.
If $n>4\Delta+5$, attach a path of length $n-(4\Delta+5)$ at
that vertex, using new vertices for the rest of the path. If equality
holds, make no change. The maximum degree remains $\Delta$, and it remains of class~2
since it contains the original graph.

Keep the $2\Delta$ indexed matchings from the proof of the proposition.
At a vertex of degree $d<\Delta$, exactly $2\Delta-2d\ge2$ of them
leave the vertex uncovered. Give a new pendant edge to two such
matchings. Each edge is again covered twice. Repeat at the new end
until the prescribed path is present. The same count on every odd set
then gives the required bound.
\end{proof}

\subsection{The regular bound}\label{sharp:sec:regular}

We now give the sharpness construction for the regular bound.
Starting from the Petersen multigraph in
Proposition~\ref{sharp:thm:family}, we replace all five vertices
$x_i$ of the outer cycle by the corresponding graphs $H_i$.
Keeping the fifth replacement makes the resulting graph regular.

\begin{proposition}\label{sharp:prop:regular}
For every integer $\Delta\ge3$, there is a connected simple
$\Delta$-regular graph $G$ with the following properties.
\begin{enumerate}
\item $\chi'(G)=\Delta+1$ and $G$ has no $\Delta$-overfull set.
\item Its order is $5\Delta+10$ when $\Delta$ is even and
$5\Delta+5$ when $\Delta$ is odd.
\end{enumerate}
\end{proposition}
\begin{proof}
Set $q=2\lfloor\Delta/2\rfloor+2$ and take the graph $D$ and its
fixed $1$-factorization from the proof of
Proposition~\ref{sharp:thm:family}. Thus $D=K_{\Delta+1}$ when
$\Delta$ is odd, and $D$ is obtained from $K_{\Delta+2}$ by deleting
one perfect matching when $\Delta$ is even. In both cases $D$ has
order $q$ and is $\Delta$-regular.

Take vertex-disjoint copies $D_i$ for $i\in\mathbb Z_5$. Choose
$z_i\in V(D_i)$, set $H_i=D_i-z_i$, and choose distinct
$u_i,v_i\in N_{D_i}(z_i)$. Let
$W_i=N_{D_i}(z_i)\setminus\{u_i,v_i\}$. Add vertices
$w_i$, $i\in\mathbb Z_5$, and take, besides the edges inside the
$H_i$, precisely the edges
\[
 \{v_i u_{i+1}:i\in\mathbb Z_5\},\qquad
 \{w_iw_{i+2}:i\in\mathbb Z_5\},\qquad
 \bigcup_{i\in\mathbb Z_5}\{w_i x:x\in W_i\}.
\]
All subscripts are read modulo five. Each neighbor of $z_i$ regains
its one lost edge, and every other vertex of $H_i$ retains degree
$\Delta$. Also $d_G(w_i)=2+|W_i|=\Delta$. The graph is simple and
connected, since every $H_i$ is connected and $W_i\ne\emptyset$.
Its order is $5(q-1)+5=5q$, as required.

Deleting the edges inside the five graphs $H_i$ and identifying
each $V(H_i)$ to a vertex $x_i$ gives the multigraph $R'$ used
in Proposition~\ref{sharp:thm:family}: the Petersen graph with each
spoke replaced by $\Delta-2$ parallel copies. Suppose that $G$
has an edge $\Delta$-coloring. Every color class is a perfect
matching. Since $|V(H_i)|=q-1$ is odd, every color meets
$\partial_G(V(H_i))$, which has exactly $\Delta$ edges. Thus each
color occurs on this cut exactly once, and contraction gives an
edge $\Delta$-coloring of $R'$. This is impossible: every perfect
matching of $R'$ contains either zero or four nonspoke edges,
whereas $R'$ has ten nonspoke edges. Vizing's theorem now gives
$\chi'(G)=\Delta+1$.

It remains to verify the odd-set bound. Use the $2\Delta$ indexed
perfect matchings of $R'$ constructed in the proof of
Proposition~\ref{sharp:thm:family}, which cover every edge twice.
For each matching and each $i\in\mathbb Z_5$, its edge at $x_i$
corresponds to an edge $z_i u$ of $D_i$. Inside $H_i$, use the
factor of $D_i$ containing $z_i u$, with $z_i$ deleted; the joining
edge covers the remaining vertex $u$. This extends the matching
to a perfect matching of $G$. Every factor of $D_i$ is used twice,
since its edge at $z_i$ corresponds to an edge of $R'$ covered
twice. Consequently these $2\Delta$ perfect matchings cover every
edge of $G$ twice. For every odd set $X$,
\[
 2e_G(X)\le 2\Delta\frac{|X|-1}{2}
          =\Delta(|X|-1),
\]
so $G$ has no $\Delta$-overfull set.
\end{proof}

For even $n$, the condition $\Delta\ge2\lfloor n/10\rfloor$ is
equivalent to $n\le5\Delta+8$ when $\Delta$ is even and
$n\le5\Delta+3$ when $\Delta$ is odd. The examples above have the
first even order beyond the respective bounds. In particular,
for every integer $r\ge1$ there is a connected regular graph of
order $n=10r+10$ and degree $\Delta=2r+1=2\lfloor n/10\rfloor-1$ that
has no $\Delta$-overfull set and no $1$-factorization. Thus the degree
bound in Theorem~\ref{intro:thm:regular} is best possible for
both parities of $\Delta$.

\section*{Acknowledgements}
While we were revising this manuscript, Dr.~Leilei Zhang informed
us by email that, with assistance from GPT, he had obtained
results corresponding to Theorems~\ref{main:quarter}
and~\ref{intro:thm:regular} for sufficiently large graphs under the
degree conditions $\Delta\ge n/4$ and $\Delta\ge n/5$, respectively.
His arguments build on the methods in the earlier arXiv submission
of Chen and Shan~\cite{ChenShan2026}. We thank him for sharing these
improvements with us.

\section*{AI statement}
ChatGPT 6 Astra was used to assist with technical details, and
ChatGPT Pro was subsequently used, through the chat interface,
to revise the manuscript. 
The authors take full responsibility
for the content of the manuscript.

\end{document}